\documentclass[10pt,a4paper]{article}
\usepackage[utf8]{inputenc}

\usepackage{mathrsfs,amsthm,xcolor,verbatim,bbm,amsmath,amsfonts,amssymb,nicefrac,enumitem,hyperref,bm,mathtools,xparse,etoolbox}

\numberwithin{equation}{section}
\usepackage{thmtools}
\usepackage[margin=1in]{geometry}
\usepackage[capitalise,nameinlink]{cleveref} %
\usepackage{tikz}
\usetikzlibrary{matrix,chains,positioning,decorations.pathreplacing,arrows,shapes,math,arrows.meta,calc, shapes.geometric}
\usepackage[most]{tcolorbox}

\makeatletter
\AtBeginEnvironment{tcolorbox}{%
   \def\@mpfn{footnote}%
}
\makeatother

\makeatletter
\usepackage{capt-of}
\usepackage{stmaryrd} %

\crefname{enumi}{item}{items}
\crefname{equation}{}{}
\crefname{subsection}{Subsection}{Subsections}
\crefname{figure}{Figure}{Figures}

\usepackage{xurl}

\theoremstyle{plain}
\newtheorem{theorem}{Theorem}[section]

\newtheorem{prop}[theorem]{Proposition}

\newtheorem{definition}[theorem]{Definition}



\theoremstyle{definition}

\DeclareMathAlphabet{\mathpzc}{OT1}{pzc}{m}{it}

\DeclareFontEncoding{LS1}{}{}
\DeclareFontSubstitution{LS1}{stix}{m}{n}
\DeclareMathAlphabet{\mathscr}{LS1}{stixscr}{m}{n}

\renewcommand{\P}{\mathbb{P}}

\newcommand{\R}{\mathbb{R}}
\newcommand{\N}{\mathbb{N}}

\newcommand{\smalll}{\mathbb{L}}

\makeatletter
\newcommand{\bvert}{%
  \mathord{%
    \mathchoice
      {\makebox[0.40em][c]{\scalebox{4.0}[1]{$\m@th\displaystyle\vert$}}}
      {\makebox[0.40em][c]{\scalebox{4.0}[1]{$\m@th\textstyle\vert$}}}
      {\makebox[0.38em][c]{\scalebox{4.0}[1]{$\m@th\scriptstyle\vert$}}}
      {\makebox[0.36em][c]{\scalebox{3.5}[1]{$\m@th\scriptscriptstyle\vert$}}}%
  }%
}
\newcommand{\bbvert}{%
  \mathord{%
    \mathchoice
      {\makebox[0.40em][c]{\scalebox{6.0}[1]{$\m@th\displaystyle\bigl\vert$}}}
      {\makebox[0.40em][c]{\scalebox{6.0}[1]{$\m@th\textstyle\big\vert$}}}
      {\makebox[0.38em][c]{\scalebox{6.0}[1]{$\m@th\scriptstyle\big\vert$}}}
      {\makebox[0.36em][c]{\scalebox{5.5}[1]{$\m@th\scriptscriptstyle\big\vert$}}}%
  }%
}
\makeatother

\newcommand{\bignnorm}[1]{\mathopen{\bbvert}#1\mathclose{\bbvert}}
\newcommand{\nnorm}[1]{\mathopen{\bvert}#1\mathclose{\bvert}}

\renewcommand{\d}{ \mathrm{d}}

\renewcommand{\c}[1]{\mathfrak{c}^{#1}}

\newcommand{\NewtonSchulzAlgorithm}[2]{\operatorname{NS}_{#1}^{#2}}
\newcommand{\TNewtonSchulzAlgorithm}[3]{\operatorname{NS}_{#1}^{#2}(#3)}
\newcommand{\TbigNewtonSchulzAlgorithm}[3]{\operatorname{NS}_{#1}^{#2}\bigl(#3\bigr)}
\newcommand{\nscoe}{b}

\newcommand{\cA}{\mathcal{A}}

\newcommand{\cL}{\mathcal{L}}

\newcommand{\cond }{ \chi }

\newcommand{\bfb}{\mathbf{b}}

\newcommand{\bfe}{\mathbf{e}}

\newcommand{\bfm}{\mathbf{m}}

\newcommand{\bfu}{\mathbf{u}}
\newcommand{\bfv}{\mathbf{v}}

\newcommand{\sgn}{\operatorname{sign}}

\newcommand{\scrC}{\mathscr{C}}
\newcommand{\scrc}{\mathscr{c}}

\newcommand{\convexcoe}{r}

\newcommand{\g}{\|\mathfrak{g}\|_{Lip}}

\newcommand{\fC}{\mathfrak{C}}
\newcommand{\fD}{\mathfrak{D}}

\newcommand{\bbC}{\mathbb{C}}

\newcommand{\fd}{\mathfrak{d}}

\newcommand{\cd}{\mathbf{X}}
\newcommand{\nb}{N}

\newcommand{\bbV}{\mathbb{V}}

\newcommand\restr[2]{{%
  \left.\kern-\nulldelimiterspace %
  #1 %
  \vphantom{|} %
  \right|_{#2} %
  }}

\DeclarePairedDelimiter{\abs}{\lvert}{\rvert}

\newcommand{\qandq}{\quad \text{and} \quad }
\newcommand{\qqandqq}{\qquad\text{and}\qquad}

\newcommand{\democrator}{\mathbb{V}}

\makeatletter
\newcommand{\doublefootnotemark}{%
  \refstepcounter{footnote}%
  \Hy@raisedlink{\hyper@anchorstart{Hfootnote.\theHfootnote}\hyper@anchorend}%
  \edef\firstfootnote{\thefootnote}%
  \edef\firstHfootnote{\theHfootnote}%
  \refstepcounter{footnote}%
  \Hy@raisedlink{\hyper@anchorstart{Hfootnote.\theHfootnote}\hyper@anchorend}%
  \edef\secondfootnote{\thefootnote}%
  \edef\secondHfootnote{\theHfootnote}%
  \textsuperscript{\normalfont
    \hyperlink{Hfootnote.\firstHfootnote}{\firstfootnote},%
    \hyperlink{Hfootnote.\secondHfootnote}{\secondfootnote}}%
}
\makeatother

\usepackage{cleveref}

\ExplSyntaxOn

\seq_new:N \g_cflist_loaded
\seq_new:N \g_cflist_pending

\NewDocumentCommand{\cfadd} { m } {
  \seq_if_in:NnF \g_cflist_loaded { #1 } {
    \seq_if_in:NnF \g_cflist_pending { #1 } {
      \seq_gput_right:Nn \g_cflist_pending { #1 }
    }
  }
}

\NewDocumentCommand{\cfconsiderloaded} { m } {
  \seq_gput_right:Nn \g_cflist_loaded {#1}
}

\NewDocumentCommand{\cfremove} { m } {
  \seq_gremove_all:Nn \g_cflist_pending { #1 }
}

\NewDocumentCommand{\cfload} { o } {
  \seq_if_empty:NTF \g_cflist_pending {
    \IfValueTF{#1}{\ignorespaces}{\unskip}
  } {
    (cf.\ \cref{\seq_use:Nn \g_cflist_pending {,}})\IfValueTF{#1}{#1~}{\unskip}
    \seq_gconcat:NNN \g_cflist_loaded \g_cflist_loaded \g_cflist_pending
    \seq_gclear:N \g_cflist_pending
    \IfValueT{#1}{\ignorespaces}
  }
}

\NewDocumentCommand{\cfclear} {} {
  \seq_gclear:N \g_cflist_loaded
  \seq_gclear:N \g_cflist_pending
}

\NewDocumentCommand{\cfout} { o } {
  \seq_if_empty:NTF \g_cflist_pending {\unskip\IfValueT{#1}{\ignorespaces}} {
    (cf.\ \cref{\seq_use:Nn \g_cflist_pending {,}})\IfValueTF{#1}{#1~}{\unskip}
    \seq_gclear:N \g_cflist_pending
    \IfValueT{#1}{\ignorespaces}
  }
}

\NewDocumentCommand{\ifnocf} { m } {
  \seq_if_empty:NT \g_cflist_pending { #1 }
}

\ExplSyntaxOff

\ExplSyntaxOn

\bool_new:N \g_noteobserve

\NewDocumentCommand{\setnote}{}{
  \bool_gset_true:N \g_noteobserve
}

\NewDocumentCommand{\setobserve}{}{
  \bool_gset_false:N \g_noteobserve
}

\NewDocumentCommand{\nobs}{ o }{
  \IfValueT{#1}{
    \str_if_eq:noTF {note} {#1} {
      \bool_gset_true:N \g_noteobserve
    } {
      \str_if_eq:noTF {Note} {#1} {
        \bool_gset_true:N \g_noteobserve
      } {
        \bool_gset_false:N \g_noteobserve
      }
    }
  }
  \bool_if:nTF { \g_noteobserve } {
    \bool_gset_false:N \g_noteobserve
    note
  } {
    \bool_gset_true:N \g_noteobserve
    observe
  }
  \IfValueF{#1}{~}
}

\NewDocumentCommand{\Nobs}{ o }{
  \IfValueT{#1}{
    \str_if_eq:noTF {note} {#1} {
      \bool_gset_true:N \g_noteobserve
    } {
      \str_if_eq:noTF {Note} {#1} {
        \bool_gset_true:N \g_noteobserve
      } {
        \bool_gset_false:N \g_noteobserve
      }
    }
  }
  \bool_if:nTF { \g_noteobserve } {
    \bool_gset_false:N \g_noteobserve
    Note
  } {
    \bool_gset_true:N \g_noteobserve
    Observe
  }
  \IfValueF{#1}{~}
}

\ExplSyntaxOff

\ExplSyntaxOn

\bool_new:N \g_hencetherefore

\NewDocumentCommand{\hence}{ o }{
  \IfValueT{#1}{
    \str_if_eq:noTF {hence} {#1} {
      \bool_gset_true:N \g_hencetherefore
    } {
      \str_if_eq:noTF {Hence} {#1} {
        \bool_gset_true:N \g_hencetherefore
      } {
        \bool_gset_false:N \g_hencetherefore
      }
    }
  }
  \bool_if:nTF { \g_hencetherefore } {
    \bool_gset_false:N \g_hencetherefore
    hence
  } {
    \bool_gset_true:N \g_hencetherefore
    therefore
  }
  \IfValueF{#1}{~}
}

\NewDocumentCommand{\Hence}{ o }{
  \IfValueT{#1}{
    \str_if_eq:noTF {hence} {#1} {
      \bool_gset_true:N \g_hencetherefore
    } {
      \str_if_eq:noTF {Hence} {#1} {
        \bool_gset_true:N \g_hencetherefore
      } {
        \bool_gset_false:N \g_hencetherefore
      }
    }
  }
  \bool_if:nTF { \g_hencetherefore } {
    \bool_gset_false:N \g_hencetherefore
    Hence,~we~obtain
  } {
    \bool_gset_true:N \g_hencetherefore
    Therefore,~we~obtain
  }
  \IfValueF{#1}{~}
}

\ExplSyntaxOff

\ExplSyntaxOn

\seq_const_from_clist:Nn \g_prove_mru {
  establish,
  demonstrate,
  prove,
  show,
  imply,
  ensure
}

\prop_new:N \l__verbs
\prop_put:Nnn \l__verbs {show} {shows}
\prop_put:Nnn \l__verbs {imply} {implies}
\prop_put:Nnn \l__verbs {demonstrate} {demonstrates}
\prop_put:Nnn \l__verbs {prove} {proves}
\prop_put:Nnn \l__verbs {establish} {establishes}
\prop_put:Nnn \l__verbs {ensure} {ensures}
\prop_put:Nnn \l__verbs {assure} {assures}

\tl_new:N \g_wordtmp
\seq_new:N \l_mytmps

\cs_generate_variant:Nn \str_if_in:nnTF { nVTF }
\cs_generate_variant:Nn \str_if_in:nnTF { xVTF }

\NewDocumentCommand{\prove}{ o }{
  \IfValueTF{#1}{
    \seq_clear:N \l_mytmps
    \seq_map_inline:Nn \g_prove_mru {
      \str_if_eq:nnTF {##1} {ensure} {
        \str_set:Nn \l_temps {n}
      } {
        \str_set:Nx \l_temps {\str_head_ignore_spaces:n {##1}}
      }
      \str_if_in:xVTF {#1} \l_temps {
        \seq_put_right:Nn \l_mytmps {##1}
      } { }
    }
    \seq_get_right:NN \l_mytmps \g_wordtmp
  } {
    \seq_get_right:NN \g_prove_mru \g_wordtmp
  }
  \tl_use:N \g_wordtmp
  \IfValueTF{#1}{}{~}
  \seq_gput_left:NV \g_prove_mru \g_wordtmp
  \seq_gremove_duplicates:N \g_prove_mru
}

\NewDocumentCommand{\proves}{ o }{
  \IfValueTF{#1}{
    \seq_clear:N \l_mytmps
    \seq_map_inline:Nn \g_prove_mru {
      \str_if_eq:nnTF {##1} {ensure} {
        \str_set:Nn \l_temps {n}
      } {
        \str_set:Nx \l_temps {\str_head_ignore_spaces:n {##1}}
      }
      \str_if_in:xVTF {#1} \l_temps {
        \seq_put_right:Nn \l_mytmps {##1}
      } { }
    }
    \seq_get_right:NN \l_mytmps \g_wordtmp
  } {
    \seq_get_right:NN \g_prove_mru \g_wordtmp
  }
  \str_set:NV \l_tmpa_str \g_wordtmp
  \prop_get:NVN \l__verbs \l_tmpa_str \l_tmpa_tl
  \tl_use:N \l_tmpa_tl
  \IfValueTF{#1}{}{~}
  \seq_gput_left:NV \g_prove_mru \g_wordtmp
  \seq_gremove_duplicates:N \g_prove_mru
}

\newcommand{\llabel}[1]{\savelabel{#1}\label{\loc.#1}\ignorespaces}

\clist_new:N \l_localreflist
\clist_new:N \l_reflist

\NewDocumentCommand{\lref} { m } {
  \clist_set:No \l_localreflist {#1}
  \clist_clear:N \l_reflist
  \clist_map_inline:Nn \l_localreflist { \clist_put_right:Nn \l_reflist {\loc.##1} }
  \cref{\l_reflist}
}

\NewDocumentCommand{\Lref} { m } {
  \clist_set:No \l_localreflist {#1}
  \clist_clear:N \l_reflist
  \clist_map_inline:Nn \l_localreflist { \clist_put_right:Nn \l_reflist {\loc.##1} }
  \Cref{\l_reflist}
}

\NewDocumentCommand{\itref}{ m m }{
  \clist_set:No \l_localreflist {#2}
  \clist_clear:N \l_reflist
  \clist_map_inline:Nn \l_localreflist { \clist_put_right:Nn \l_reflist {#1.##1} }
  \cref{\l_reflist}~in~\cref{#1}
}

\seq_new:N \l_enum_seq
\int_new:N \l_num_items

\bool_new:N \g_commaused_bool

\providecommand{\comma}{}

\cs_new:Nn \enum_it:nn {
  \int_case:nnF {\l_num_items - #1} {
    {0} {
      \renewcommand{\comma}{}
      #2\space
    }
    {1} {
      \bool_gset_false:N \g_commaused_bool
      \renewcommand{\comma}{,~\bool_gset_true:N \g_commaused_bool}
      #2
      \bool_if:NTF \g_commaused_bool {} {,~}
      and~
    }
  } {
    \bool_gset_false:N \g_commaused_bool
    \renewcommand{\comma}{,~\bool_gset_true:N \g_commaused_bool}
    #2
    \bool_if:NTF \g_commaused_bool {} {,~}
  }
}

\cs_new:Nn \enum_it_U:nn {
  \int_case:nnF {\l_num_items - #1} {
    {0} {
      \renewcommand{\comma}{}
      #2
      \space
    }
    {1} {
      \bool_gset_false:N \g_commaused_bool
      \renewcommand{\comma}{,~\bool_gset_true:N \g_commaused_bool}
      #2
      \bool_if:NTF \g_commaused_bool {} {,~}
      and~
    }
  } {
    \bool_gset_false:N \g_commaused_bool
    \renewcommand{\comma}{,~\bool_gset_true:N \g_commaused_bool}
    \int_compare:nTF {#1=1} {
      \text_titlecase_first:n {#2}
    } {
      #2
    }
    \bool_if:NTF \g_commaused_bool {} {,~}
  }
}

\cs_generate_variant:Nn \tl_if_eq:nnTF {onTF}

\cs_new:Nn \enum:nnnn {
  \seq_set_split:Nnn \l_enum_seq ; {#1}
  \seq_remove_all:Nn \l_enum_seq { }
  \seq_remove_all:Nn \l_enum_seq {#2}
  \seq_log:N \l_enum_seq
  \int_set:Nn \l_num_items {\seq_count:N \l_enum_seq}
  \int_log:N \l_num_items
  \int_case:nnF {\l_num_items} {
    { 0 } { 0 }
    { 1 } {
      \IfBooleanTF{#4} {
        \tl_set:Nn \l_text_case_exclude_arg_tl {\cref}
        \bool_set_false:N \l_text_titlecase_check_letter_bool
        \text_titlecase_first:n {\seq_use:Nn \l_enum_seq {}}
      } {
        \seq_use:Nn \l_enum_seq {}
      }
      \space
      \tl_if_eq:onTF{#3}{-}{}{
        \bool_if:NTF \l_plural_bool {
          \prove[#3]~
        } {
          \proves[#3]~
        }
      }
    }
    { 2 } {
      \IfBooleanTF{#4} {
        \tl_set:Nn \l_text_case_exclude_arg_tl {\cref}
        \bool_set_false:N \l_text_titlecase_check_letter_bool
        \text_titlecase_first:n {\seq_use:Nn \l_enum_seq {~and~}}
      } {
        \seq_use:Nn \l_enum_seq {~and~}
      }
      \space
      \tl_if_eq:onTF{#3}{-}{}{
        \prove[#3]~
      }
    }
  } {
    \IfBooleanTF{#4} {
      \tl_set:Nn \l_text_case_exclude_arg_tl {\cref}
      \seq_indexed_map_function:NN \l_enum_seq \enum_it_U:nn
    } {
      \seq_indexed_map_function:NN \l_enum_seq \enum_it:nn
    }
    \tl_if_eq:onTF{#3}{-}{}{
      \prove[#3]~
    }
  }
}

\cs_generate_variant:Nn \enum:nnnn {nxnn}
\cs_generate_variant:Nn \enum:nnnn {nxxn}

\NewDocumentCommand{\enum}{O{} m O{-} s}{
  \IfBooleanTF{#4}{
    \enum:nxnn {#2} {#1} {sindep} \BooleanFalse
  } {
    \enum:nxxn {#2} {#1} {#3} \BooleanFalse
  }
}

\NewDocumentCommand{\dott}{}{\ifnocf{.}\space}

\bool_new:N \g_arg_start_bool
\bool_gset_true:N \g_arg_start_bool

\NewDocumentCommand{\startnewargseq}{}{\bool_gset_true:N \g_arg_start_bool \tl_set:Nn \g_label_tl {}}

\cs_generate_variant:Nn \seq_if_in:NnTF {NxTF}
\cs_generate_variant:Nn \seq_remove_all:Nn {Nx}

\int_new:N \l_random_int

\cs_generate_variant:Nn \tl_if_head_eq_catcode:nNTF {oNTF}
\cs_generate_variant:Nn \tl_if_head_eq_catcode:nNTF {VNTF}

\cs_generate_variant:Nn \tl_if_head_eq_catcode:nNTF {eNTF}
\cs_generate_variant:Nn \tl_log:n {o}
\cs_generate_variant:Nn \tl_log:n {f}
\cs_generate_variant:Nn \tl_log:n {x}
\cs_generate_variant:Nn \tl_log:n {e}

\cs_generate_variant:Nn \tl_if_in:nnTF {onTF}
\cs_generate_variant:Nn \tl_if_in:NnTF {NeTF}

\cs_generate_variant:Nn \tl_if_head_eq_meaning:nNTF {VNTF}

\seq_const_from_clist:Nn \g_arg_mru_this {
  Ahpr,
  Tapr,
  Ctapr,
  H
}

\seq_const_from_clist:Nn \g_arg_mru_nothis {
  Ia,
  Nwc,
  N,
  Itns,
  Fm,
  Itnswc,
  Mo
}

\seq_new:N \l_arg_seq
\tl_new:N \l_cons_tl
\tl_new:N \l_dummy_tl

\bool_new:N \g_debug_bool
\bool_gset_false:N \g_debug_bool

\bool_new:N \l_insidearg_bool

\bool_new:N \g_firstargletter_bool

\sys_gset_rand_seed:n {0903}

\bool_new:N \l_plural_bool
\tl_new:N \l_arg_verbs_tl

\NewDocumentCommand{\argument}{mom}{
  \bool_set_false:N \l_plural_bool
  \tl_set:Nn \l_arg_verbs_tl {sindep}
  \keys_define:nn { benno/argument } {
    plural .value_forbidden:n = true,
    plural .code:n = {\bool_set_true:N \l_plural_bool},
    verbs .value_required:n = false,
    verbs .tl_set:N = \l_arg_verbs_tl,
  }
  \IfValueT{#2}{
    \keys_set:nn { benno/argument } {#2}
  }
  \bool_log:N \l_plural_bool
  \bool_gset_true:N \l_insidearg_bool
  \seq_set_split:Nnn \l_arg_seq ; {#1}
  \seq_remove_all:Nn \l_arg_seq { }
  \seq_log:N \l_arg_seq
  \tl_set:Nn \l_cons_tl {#3}
  \tl_trim_spaces:N \l_cons_tl
  \seq_if_in:NxTF \l_arg_seq {\lref{\g_label_tl}} {
    \seq_remove_all:Nx \l_arg_seq {\lref{\g_label_tl}}
    \seq_get_left:NNTF \l_arg_seq \l_dummy_tl {
      \tl_trim_spaces:N \l_dummy_tl
      \bool_gset_false:N \g_firstargletter_bool
      \tl_if_head_eq_catcode:VNTF \l_dummy_tl a {
        \bool_gset_true:N \g_firstargletter_bool
      } {
        \tl_if_head_eq_meaning:VNTF \l_dummy_tl {\cref} {
          \tl_set:Nx \l_tmpa_tl {\tl_tail:N \l_dummy_tl}
          \tl_set:Nx \l_tmpb_tl {\tl_head:N \l_tmpa_tl}
          \bool_gset_true:N \g_firstargletter_bool
          \tl_if_in:NeTF \l_tmpb_tl {lem\c_colon_str} {} {
            \tl_if_in:NeTF \l_tmpb_tl {thm\c_colon_str} {} {
              \tl_if_in:NeTF \l_tmpb_tl {prop\c_colon_str} {} {
                \tl_if_in:NeTF \l_tmpb_tl {cor\c_colon_str} {} {
                  \bool_gset_false:N \g_firstargletter_bool
                }
              }
            }
          }
        } {
        }
      }
      \bool_if:NTF \g_firstargletter_bool {
        \seq_set_eq:NN \l_tmpa_seq \g_arg_mru_this
        \seq_remove_all:Nn \l_tmpa_seq {H}
        \seq_get_right:NN \l_tmpa_seq \l_tmpa_tl
        \int_case:nnF {\seq_count:N \l_arg_seq} {
          {1} {
            \str_case:VnF {\l_tmpa_tl} {
              {Ahpr} {
                \bool_if:NT \g_debug_bool {C1.1}
                \seq_gput_left:Nn \g_arg_mru_this {Ahpr}
                \seq_gremove_duplicates:N \g_arg_mru_this
                \enum:nxnn {#1} {\lref{\g_label_tl}} {-} {\BooleanTrue}
                \hence~
                \bool_if:NTF \l_plural_bool {
                  \prove[\l_arg_verbs_tl]~\ignorespaces #3
                } {
                  \proves[\l_arg_verbs_tl]~\ignorespaces #3
                }
              }
              {Tapr} {
                \bool_if:NT \g_debug_bool {C1.2}
                \seq_gput_left:Nn \g_arg_mru_this {Tapr}
                \seq_gremove_duplicates:N \g_arg_mru_this
                \enum[\lref{\g_label_tl}]{
                  This;
                  #1
                }[\l_arg_verbs_tl]\ignorespaces #3
              }
              {Ctapr} {
                \bool_if:NT \g_debug_bool {C1.3}
                \seq_gput_left:Nn \g_arg_mru_this {Ctapr}
                \seq_gremove_duplicates:N \g_arg_mru_this
                Combining~
                \enum[\lref{\g_label_tl}]{
                  this;
                  #1
                } \proves[\l_arg_verbs_tl]~\ignorespaces #3
              }
            } {}
          }
        } {
          \str_case:VnF {\l_tmpa_tl} {
             {Ahpr} {
              \bool_if:NT \g_debug_bool {C2.1}
              \seq_gput_left:Nn \g_arg_mru_this {Ahpr}
              \seq_gremove_duplicates:N \g_arg_mru_this
              \enum:nxnn {#1} {\lref{\g_label_tl}} {-} {\BooleanTrue}
              \hence~
              \prove[\l_arg_verbs_tl]~\ignorespaces #3
            }
            {Tapr} {
              \bool_if:NT \g_debug_bool {C2.2}
              \seq_gput_left:Nn \g_arg_mru_this {Tapr}
              \seq_gremove_duplicates:N \g_arg_mru_this
              \enum[\lref{\g_label_tl}]{
                This;
                #1
              }[\l_arg_verbs_tl]\ignorespaces #3
            }
            {Ctapr} {
              \int_case:nn {\int_rand:nn {0} {1}} {
                {0} {
                  \bool_if:NT \g_debug_bool {C2.3}
                  \seq_gput_left:Nn \g_arg_mru_this {Ctapr}
                  \seq_gremove_duplicates:N \g_arg_mru_this
                  Combining~
                  \enum[\lref{\g_label_tl}]{
                    this;
                    #1
                  } \proves[\l_arg_verbs_tl]~\ignorespaces #3
                }
                {1} {
                  \bool_if:NT \g_debug_bool {C2.4}
                  \seq_gput_left:Nn \g_arg_mru_this {Ctapr}
                  \seq_gremove_duplicates:N \g_arg_mru_this
                  Combining~
                  \enum:nxnn {#1} {\lref{\g_label_tl}} {-} {\BooleanFalse}
                  \hence~
                  \proves[\l_arg_verbs_tl]~\ignorespaces #3
                }
              }
            }
          } {}
        }
      } {
        \seq_set_eq:NN \l_tmpa_seq \g_arg_mru_this
        \seq_remove_all:Nn \l_tmpa_seq {H}
        \seq_remove_all:Nn \l_tmpa_seq {Ahpr}
        \seq_get_right:NN \l_tmpa_seq \l_tmpa_tl
        \int_case:nnF {\seq_count:N \l_arg_seq} {
          {1} {
            \str_case:VnF {\l_tmpa_tl} {
              {Tapr} {
                \bool_if:NT \g_debug_bool {C3.1}
                \seq_gput_left:Nn \g_arg_mru_this {Tapr}
                \seq_gremove_duplicates:N \g_arg_mru_this
                \enum[\lref{\g_label_tl}]{
                  This;
                  #1
                }[\l_arg_verbs_tl]\ignorespaces #3
              }
              {Ctapr} {
                \bool_if:NT \g_debug_bool {C3.2}
                \seq_gput_left:Nn \g_arg_mru_this {Ctapr}
                \seq_gremove_duplicates:N \g_arg_mru_this
                Combining~
                \enum[\lref{\g_label_tl}]{
                  this;
                  #1
                } \proves[\l_arg_verbs_tl]~\ignorespaces #3
              }
            } {}
          }
        } {
          \str_case:VnF {\l_tmpa_tl} {
            {Tapr} {
              \bool_if:NT \g_debug_bool {C4.1}
              \seq_gput_left:Nn \g_arg_mru_this {Tapr}
              \seq_gremove_duplicates:N \g_arg_mru_this
              \enum[\lref{\g_label_tl}]{
                This;
                #1
              }[\l_arg_verbs_tl]\ignorespaces #3		
            }
            {Ctapr} {
              \int_case:nn {\int_rand:nn {0} {1}} {
                {0} {
                  \bool_if:NT \g_debug_bool {C4.2}
                  \seq_gput_left:Nn \g_arg_mru_this {Ctapr}
                  \seq_gremove_duplicates:N \g_arg_mru_this
                  Combining~
                  \enum[\lref{\g_label_tl}]{
                    this;
                    #1
                  } \proves[\l_arg_verbs_tl]~\ignorespaces #3		
                }
                {1} {
                  \bool_if:NT \g_debug_bool {C4.3}
                  \seq_gput_left:Nn \g_arg_mru_this {Ctapr}
                  \seq_gremove_duplicates:N \g_arg_mru_this
                  Combining~
                  \enum:nxnn {#1} {\lref{\g_label_tl}} {-} {\BooleanFalse}
                  \hence~
                  \proves[\l_arg_verbs_tl]~\ignorespaces #3    
                }
              }
            }
          } {}
        }
      }
    } {
      \tl_if_head_eq_catcode:oNTF \l_cons_tl a {
        \seq_set_eq:NN \l_tmpa_seq \g_arg_mru_this
        \seq_remove_all:Nn \l_tmpa_seq {Ctapr}
        \seq_remove_all:Nn \l_tmpa_seq {Ahpr}
        \seq_get_right:NN \l_tmpa_seq \l_tmpa_tl
        \str_case:VnF {\l_tmpa_tl} {
          {H} {
            \bool_if:NT \g_debug_bool {C5.1}
            \seq_gput_left:Nn \g_arg_mru_this {H}
            \seq_gremove_duplicates:N \g_arg_mru_this
            Hence,~we~obtain~\ignorespaces #3
          }
          {Tapr} {
            \bool_if:NT \g_debug_bool {C5.2}
            \seq_gput_left:Nn \g_arg_mru_this {Tapr}
            \seq_gremove_duplicates:N \g_arg_mru_this
            This~\proves[\l_arg_verbs_tl]~\ignorespaces #3
          }
        } {}
      } {
        \bool_if:NT \g_debug_bool {C6.1}
        \seq_gput_left:Nn \g_arg_mru_this {Tapr}
        \seq_gremove_duplicates:N \g_arg_mru_this
        This~\proves[\l_arg_verbs_tl]~\ignorespaces #3
      }
    } 
  } {
    \int_compare:nNnTF {\seq_count:N \l_arg_seq} = {0} {
      \bool_if:NTF \g_arg_start_bool {
        \bool_if:NT \g_debug_bool {C7.1}
        \Nobs\unskip
        #3
      } {
        \bool_if:NT \g_debug_bool {C7.2}
        \Moreover~
        #3
      }
    } {
      \bool_if:NTF \g_arg_start_bool {
        \bool_if:NT \g_debug_bool {C8.1}
        \tl_log:N \l_arg_verbs_tl
        \Nobs~that~
        \enum{
          #1
        }[\l_arg_verbs_tl]\ignorespaces #3
      } {
        \int_compare:nNnTF {\seq_count:N \l_arg_seq} = {1} {
          \seq_set_eq:NN \l_tmpa_seq \g_arg_mru_nothis
          \seq_remove_all:Nn \l_tmpa_seq {Nwc}
          \seq_remove_all:Nn \l_tmpa_seq {Itnswc}
          \seq_get_right:NN \l_tmpa_seq \l_tmpa_tl
        } {
          \seq_get_right:NN \g_arg_mru_nothis \l_tmpa_tl
        }
        \str_case:VnF {\l_tmpa_tl} {
          {Mo} {
            \bool_if:NT \g_debug_bool {C9.1}
            \seq_gput_left:Nn \g_arg_mru_nothis {Mo}
            \seq_gremove_duplicates:N \g_arg_mru_nothis
            Moreover,~\nobs~that~
            \enum{
              #1
            }[\l_arg_verbs_tl]\ignorespaces #3		
          }
          {Fm} {
            \bool_if:NT \g_debug_bool {C9.2}
            \seq_gput_left:Nn \g_arg_mru_nothis {Fm}
            \seq_gremove_duplicates:N \g_arg_mru_nothis
            Furthermore,~\nobs~that~
            \enum{
              #1
            }[\l_arg_verbs_tl]\ignorespaces #3		
          }
          {Ia} {
            \bool_if:NT \g_debug_bool {C9.3}
            \seq_gput_left:Nn \g_arg_mru_nothis {Ia}
            \seq_gremove_duplicates:N \g_arg_mru_nothis
            In~addition,~\nobs~that~
            \enum{
              #1
            }[\l_arg_verbs_tl]\ignorespaces #3		
          }
          {N} {
            \bool_if:NT \g_debug_bool {C9.4}
            \seq_gput_left:Nn \g_arg_mru_nothis {N}
            \seq_gremove_duplicates:N \g_arg_mru_nothis
            Next,~\nobs~that~
            \enum{
              #1
            }[\l_arg_verbs_tl]\ignorespaces #3		
          }
          {Itns} {
            \bool_if:NT \g_debug_bool {C9.5}
            \seq_gput_left:Nn \g_arg_mru_nothis {Itnswc}
            \seq_gput_left:Nn \g_arg_mru_nothis {Itns}
            \seq_gremove_duplicates:N \g_arg_mru_nothis
            In~the~next~step~we~\nobs~that~
            \enum{
              #1
            }[\l_arg_verbs_tl]\ignorespaces #3		
          }
          {Nwc} {
            \bool_if:NT \g_debug_bool {C9.6}
            \seq_gput_left:Nn \g_arg_mru_nothis {Nwc}
            \seq_gremove_duplicates:N \g_arg_mru_nothis
            Next~we~combine~
            \enum{
              #1
            }to~obtain~\ignorespaces #3
          }
          {Itnswc} {
            \bool_if:NT \g_debug_bool {C9.7}
            \seq_gput_left:Nn \g_arg_mru_nothis {Itns}
            \seq_gput_left:Nn \g_arg_mru_nothis {Itnswc}
            \seq_gremove_duplicates:N \g_arg_mru_nothis
            In~the~next~step~we~combine~
            \enum{
              #1
            }to~obtain~\ignorespaces #3
          }
        } {}
      }
    }
  }
  \bool_gset_false:N \g_arg_start_bool
  \bool_gset_false:N \l_insidearg_bool
  \cfload[.]%
  \color{black}
}

\tl_new:N \g_label_tl
\tl_gset:Nn \g_label_tl { }

\NewDocumentCommand{\savelabel}{m}{
  \bool_if:NTF \l_insidearg_bool {
    \tl_gset:Nn \g_label_tl {#1}
  } {
    \tl_gset:Nn \g_label_tl { }
  }
}

\ExplSyntaxOff

\ExplSyntaxOn

\NewDocumentEnvironment {athm} {m m o} {
\str_if_eq:noTF {example} {#1} {
  \bool_gset_true:N \g_example_bool
} {
  \bool_gset_false:N \g_example_bool
}
\cfclear
\IfNoValueTF{#3}{
\begin{#1}\label{#2}\global\def\loc{#2}
}{
\begin{#1}[#3]\label{#2}\global\def\loc{#2}
}
}{
\end{#1}
}

\NewDocumentEnvironment {adef} {m} {
\begin{definition}\label{#1}\global\def\loc{#1}
}{
\end{definition}
}

\NewDocumentEnvironment{aproof} {} {
\bool_if:NTF \g_example_bool {
  \bool_gset_true:N \g_arg_start_bool
  \begin{proof}[Proof~for~\cref{\loc}]
} {
  \bool_gset_true:N \g_arg_start_bool
  \begin{proof}[Proof~of~\cref{\loc}]
}
\bool_gset_false:N \g_finishproof_bool
}{
\bool_if:NTF \g_finishproof_bool {}
{\finishproofthus}
\end{proof}
}

\NewDocumentCommand{\finishproofthus} {} {
  \bool_gset_true:N \g_finishproof_bool 
  \bool_if:NTF \g_example_bool {
    The~proof~for~\cref{\loc}~is~thus~complete.
  } {
    The~proof~of~\cref{\loc}~is~thus~complete.
  }
}
\NewDocumentCommand{\finishproofthis} {} {
  \bool_gset_true:N \g_finishproof_bool 
  \bool_if:NTF \g_example_bool {
    This~completes~the~proof~for~\cref{\loc}.
  } {
    This~completes~the~proof~of~\cref{\loc}.
  }
}

\ExplSyntaxOff

\NewDocumentEnvironment{cproof}{m}
{\begin{proof}[Proof of \cref{#1}]}%
{\noindent The proof of \cref{#1} is thus complete.
\end{proof}}

\NewDocumentEnvironment{cproof2}{m}
{\begin{proof}[Proof of \cref{#1}]}%
{\noindent This completes the proof of \cref{#1}.
\end{proof}}

\hypersetup{
    colorlinks,
    linkcolor={blue!80!black},
    citecolor={green},
    urlcolor={blue!80!black}
}

\makeatletter
\ExplSyntaxOn
\seq_new:N \g_abbrs
\prop_new:N \g_abbr_counts
\tl_new:N \l_abbr_count_tl

\ExplSyntaxOn

\bool_new:N \g_forexample

\NewDocumentCommand{\eg}{ o }{
	\IfValueT{#1}{
		\str_if_eq:noTF {fe} {#1} {
			\bool_gset_true:N \g_forexample
		} {\bool_gset_false:N \g_forexample}
	}
	\bool_if:nTF { \g_forexample } {
		\bool_gset_false:N \g_forexample
		for~example
	}{
		\bool_gset_true:N \g_forexample
		for~instance
	}
}

\NewDocumentCommand{\abbr}{m m O{#1} m m O{#4} m}{
	\expandafter\newcommand\csname#3\endcsname[1][]{
		\seq_if_in:NnTF \g_abbrs {#1} {
			\prop_get:NnN \g_abbr_counts {#1} \l_abbr_count_tl
			\prop_gput:Nnx \g_abbr_counts {#1} {\int_eval:n {\l_abbr_count_tl + 1}}
			\hyperref[#1]{#7}
		} {
			\seq_gput_left:Nn \g_abbrs {#1}
			\prop_gput:Nnn \g_abbr_counts {#1} {1}
			\expandafter\gdef\csname#1@def\endcsname{#2}
			\phantomsection\label{#1}
			\str_if_eq:nnTF{##1}{}{\emph{#2}}{##1}~(\hyperref[#1]{#7})
		}
	}
	\expandafter\newcommand\csname#6\endcsname[1][]{
		\seq_if_in:NnTF \g_abbrs {#1} {
			\prop_get:NnN \g_abbr_counts {#1} \l_abbr_count_tl
			\prop_gput:Nnx \g_abbr_counts {#1} {\int_eval:n {\l_abbr_count_tl + 1}}
			\hyperref[#1]{#4}
		} {
			\expandafter\gdef\csname#1@def\endcsname{#5}
			\seq_gput_left:Nn \g_abbrs {#1}
			\prop_gput:Nnn \g_abbr_counts {#1} {1}
			\phantomsection\label{#1}
			\str_if_eq:nnTF{##1}{}{\emph{#5}}{##1}~(\hyperref[#1]{#4})
		}
	}
}

\ExplSyntaxOff
\makeatother
\abbr{GD}{Gradient descent}{GDs}{gradient descents}{GD}
\abbr{FW}{Frank-Wolfe}{FWs}{Frank-Wolfe}{FW}
\abbr{ANN}{artificial neural network}{ANNs}{artificial neural networks}{ANN}
\abbr{RMSprop}{root mean square propagation}{RMSprops}{the root mean square propagation}{RMSprop}
\abbr{DNN}{deep neural network}{DNNs}{deep neural networks}{DNN}
\abbr{SGD}{Stochastic gradient descent}{SGDs}{Stochastic gradient descent}{SGD}
\abbr{iid}{independent and identically distributed}{i.i.d}{independent and identically distributed}{i.i.d.}
\abbr{SOP}{stochastic optimization problem}{SOPs}{stochastic optimization problems}{SOP}
\abbr{DL}{Deep learning}{DLs}{Deep learning}{DL}
\abbr{Adam}{adaptive moment estimation}{Adams}{adaptive moment estimation \SGD}{Adam}
\abbr{NS}{Newton-Schulz}{NSs}{Newton-Schulz}{NS}
\abbr{AdamW}{\Adam\ with decoupled weight decay}{AdamWs}{\Adam\ with decoupled weight decay}{AdamW}
\abbr{MUON}{momentum orthogonalized by Newton-Schulz}{MUONs}{momentum orthogonalized by Newton-Schulz}{MUON}
\abbr{AI}{artificial intelligence}{AIs}{artificial intelligence}{AI}
\abbr{LLM}{large language model}{LLMs}{large language models}{LLM}
\abbr{PDE}{partial differential equation}{PDEs}{partial differential equations}{PDE}
\abbr{as}{almost surely}{ass}{almost surely}{a.s.}
\abbr{pas}{$\P$-almost surely}{pass}{almost surely}{$\P$-a.s.}
\abbr{ReLUs}{rectified linear unit}{ReLU}{rectified linear unit}
\makeatother

\title{
Optimal higher order convergence\\rates for the MUON optimizer
}
\author{Thomas Beretti$^{1}$, Steffen Dereich$^{2}$, Thang Do$^{3}$, and Arnulf Jentzen$^{4,5}$
 	\bigskip
 	\\
     \small{$^1$ e-mail: \texttt{thomasberetti@outlook.com}}
  \smallskip
 	\\
 	\small{$^2$ Institute for Mathematical Stochastics, Faculty of Mathematics and Computer Science,}\vspace{-0.1cm}\\
 \small{University of M\"unster, Germany, e-mail: \texttt{steffen.dereich@uni-muenster.de}}
 \smallskip
 \\
     \small{$^3$ School of Data Science, The Chinese University of Hong Kong, Shenzhen}
 	\vspace{-0.1cm}\\
 	\small{ (CUHK-Shenzhen), China, e-mail: \texttt{minhthangdo@link.cuhk.edu.cn}}
  \smallskip
 \\
 	\small{$^4$ School of Data Science and School of Artificial Intelligence, The Chinese University}
 	\vspace{-0.1cm}\\
 	\small{of Hong Kong, Shenzhen (CUHK-Shenzhen), China, e-mail: \texttt{ajentzen@cuhk.edu.cn}}
 	\smallskip
 	\\
  \small{$^5$ Applied Mathematics: Institute for Analysis and Numerics, Faculty of Mathematics and}
 	\vspace{-0.1cm}\\
 	\small{Computer Science, University of M{\"u}nster, Germany, e-mail: \texttt{ajentzen@uni-muenster.de}}
 	\smallskip
 	\\
 }

\date{\today}
\allowdisplaybreaks 
\begin{document}
\maketitle
\begin{abstract}
Gradient descent optimization methods are the key tools to train deep neural networks in artificial intelligence systems. In practice, one typically employs suitable accelerated gradient descent optimization methods such as the momentum optimizer, the Adam optimizer, or the MUON optimizer proposed by Jordan et al.\ in 2024. In 1964 Polyak showed that, in a neighborhood of a local minimizer, the momentum optimizer converges with the rate $(\sqrt{\chi}+1)^{-1}(\sqrt{\chi}-1)$, which is strictly faster than the rate $(\chi+1)^{-1}(\chi-1)$ of the standard gradient descent optimizer, where $\chi\in(1,\infty)$ is the condition number (the ratio of the largest to the smallest eigenvalue) of the Hessian of the objective function at the local minimizer. Recently, the same faster rate has been established for the Adam optimizer. In this work we show that the MUON optimizer achieves the same higher speed of convergence and we establish the optimality of the obtained convergence rates. More specifically, we propose a precise concept for measuring the speed of convergence of gradient descent optimization methods, and we establish optimal convergence rates for the standard gradient descent, the momentum, the RMSprop, the Adam, and the MUON optimizers.
\end{abstract}
\tableofcontents
\section{Introduction}
\GD\ optimization methods are the key tools to train \DNNs\ in \AI\ systems (cf., \eg, \cite{Ruder2016AnOO}, \cite[Chapter~8]{GoodfellowBengioCourville2016}, and \cite[Chapters~6 and~7]{ArBePhi2024}).
In practically relevant training  of large-scale \AI\ systems, the standard (stochastic) \GD\ method is typically not the method of choice. Instead, one employs suitable accelerated \GD\ optimization methods such as the momentum optimizer~\cite{Polyak1964SomeMO}, the \Adam\ optimizer~\cite{KingmaBa2024_Adam}, the \AdamW\ optimizer~\cite{Loshchilov2017DecoupledWD}, or the recently introduced \MUON\ optimizer~\cite{jordan2024muon}.
The momentum and the \Adam\ optimizers are known to achieve a faster speed of convergence compared to the standard \GD\ optimizer~\cite{DeArAdsharp2025,Polyak1964SomeMO}. In this work we show that the \MUON\ optimizer achieves the same higher speed of convergence and we establish the optimality of the obtained convergence orders.

In \cref{subsec: rate of convergence} below, we present in \cref{definition: rate of convergence} the precise concept that we propose for measuring the speed of convergence of different \GD\ optimization methods.
In \cref{subsec: main result} below, we present the main result of this work, \cref{main theorem}, which establishes optimal convergence rates for the \GD, the momentum \GD, the \RMSprop, the \Adam, and the \MUON\ optimizers. In particular, we show that the momentum \GD, the \Adam, and the \MUON\ optimizers converge with the same rate, which is strictly faster than the rate of convergence of the standard \GD\ and the \RMSprop\ optimizers.

\subsection{Definition: Rate of convergence}\label{subsec: rate of convergence}

For a general gradient-based optimization method, we specify in \cref{definition: rate of convergence} below what we mean by the rate of convergence of the considered optimization method (cf.\ \cite[Definition~1.1]{DeArAdsharp2025}).
The natural number $\fD\in\N=\{1,2,3,\ldots\}$ in \cref{definition: rate of convergence} denotes the dimension (the number of degrees of freedom) of the optimization problem under consideration.

 \begin{samepage}
    \begin{tcolorbox}[colback=white!95!gray,
                  colframe=black,
                  boxrule=0.5pt,
                  sharp corners,
                  enhanced,
                  unbreakable,
                 ]
		\begin{definition}[\textcolor{red}{Rate of convergence}]\label{definition: rate of convergence}
			Let $\fD\in \N$, $\chi\in (1,\infty)$, $\varrho\in [0,1)$
			and let
			$ \Phi_n \colon(\R^\fD) ^{ n } \to \R^\fD $,
			$ n \in \N $,
			be functions.
			Then we say that
			$ ( \Phi_n )_{ n \in \N } $ converges
			with $ \cond $-rate $ \varrho $
			 if and only if it holds for all $\rho \in (0,\infty) \backslash \{ \varrho \}$ that the following two statements are equivalent:
             \begin{enumerate}[label=(\roman*)]
                 \item  It holds that $\rho > 
                 \varrho$.
                \item For all $\convexcoe \in (0,\infty)$, $\cL \in C^2( \R^\fD, \R )$, $\vartheta\in(\nabla\cL)^{-1}(\{0\})$ with $\forall\,\theta\in\R^\fD \colon \convexcoe\|\theta\|^2\leq  \cL'' (\vartheta)(\theta,\theta)\leq \convexcoe\chi\|\theta\|^2
			$ there exist $\fC,\Gamma, \delta \in (0,\infty) $ such that for every $\gamma=(\gamma_n)_{n\in\N}\colon\N\to (\Gamma-\delta,\Gamma+\delta)$ and every $ \Theta = ( \Theta_n )_{ n \in \N _0} \colon \N_0 \to \R^\fD $
			with $\| \Theta_0 - \vartheta \|< \delta$ and
			\begin{equation}
				 \textcolor{magenta}{\forall \, n \in \N \colon
				\Theta_n
				=
				\Theta_{ n - 1 }
				-
				\gamma_n \Phi_n\bigl(
				( \nabla \cL )( \Theta_0 ) ,
				( \nabla \cL )( \Theta_1 ) ,
				\dots ,
				( \nabla \cL )( \Theta_{ n - 1 } )
				\bigr)}
			\end{equation}
			we have for all $n\in\N_0$ that
			$
			\textcolor{magenta}{\| \Theta_n - \vartheta\| 
			\leq \fC\rho^{  n } \|\Theta_0-\vartheta\|}
			$.
             \end{enumerate}
		\end{definition}
        \end{tcolorbox}
        \end{samepage}

The general functions $\Phi_n\colon(\R^\fD)^n\to\R^\fD$, $n\in\N$, from \cref{definition: rate of convergence} allow us to apply \cref{definition: rate of convergence} to a large class of different gradient-based optimization methods (cf.\ \cite[Subsections~6.4 and~6.5]{MR4688424}).
For example, for the standard \GD\ optimizer (see \cref{definition: GD}) it holds for all $n\in\N$, $g_1,g_2,\dots,g_n\in\R^\fD$ that $\Phi_n(g_1,g_2,\dots,g_n)=g_n$,
for every $\alpha\in[0,1)$ and the $\alpha$-momentum \GD\ optimizer (see \cref{definition: momentum}) it holds for all $n\in\N$, $g_1,g_2,\dots,g_n\in\R^\fD$ that $\Phi_n(g_1,g_2,\dots,g_n)=(1-\alpha)\sum_{k=1}^n\alpha^{n-k}g_k$,
and for every $\beta\in(0,1)$, $\varepsilon\in(0,\infty)$ and the $\beta$-$\varepsilon$-\RMSprop\ optimizer (see \cref{definition: RMSprop}) it holds for all $n\in\N$, $g=((g_{k,j})_{j\in\{1,2,\dots,\fD\}})_{k\in\{1,2,\dots,n\}}\in(\R^\fD)^n$, $j\in\{1,2,\dots,\fD\}$ that $\Phi_n^j(g)=[\varepsilon+((1-\beta)\sum_{k=1}^n\beta^{n-k}\abs{g_{k,j}}^2)^{1/2}]^{-1}g_{n,j}$.
We refer to \cref{definition: Adam} for the description of the \Adam\ optimizer using the functions $\Phi_n\colon(\R^\fD)^n\to\R^\fD$, $n\in\N$, from \cref{definition: rate of convergence} and we refer to \cref{definition: muon} for the description of the \MUON\ optimizer using the functions $\Phi_n\colon(\R^\fD)^n\to\R^\fD$, $n\in\N$, from \cref{definition: rate of convergence}.

        \subsection{Main result: Optimal convergence rates of the MUON optimizer}\label{subsec: main result}

        Using the notion of the rate of convergence in \cref{definition: rate of convergence}, we now specify in \cref{main theorem} below the optimal rates of convergence for the \GD, the momentum \GD, the \RMSprop, the \Adam, and the \MUON\ optimizers.

        \begin{samepage}
         \begin{tcolorbox}[colback=white!95!gray,
                  colframe=black,
                  boxrule=0.5pt,
                  sharp corners,
                  enhanced,
                  unbreakable,
                 ]
                 \begin{athm}{theorem}{main theorem}
                 Let $\nb\in\N$, $\fD\in\N\backslash\{1\}$, $d=(d_1,\dots,d_\nb)$, $\fd=(\fd_1,\dots,\fd_\nb)\in\N^\nb$, $\chi\in (1,\infty)$, $\varepsilon\in (0,\infty)$, $\alpha,\beta\in (0,1)$ satisfy $\fD=\sum_{i=1}^\nb d_i\fd_i$ and $\alpha= (\sqrt{ \chi } + 1) ^{-2}(\sqrt{ \chi } - 1 )^2$. Then 
                 \begin{enumerate}[label=(\roman*)]
                     \item \label{item 1: main theorem} for every \GD\ optimizer $\Phi_n\colon(\R^\fD)^n\to\R^\fD$, $n\in\N$, it holds that $(\Phi_n)_{n\in\N}$ converges with $\chi$-rate $\textcolor{magenta}{\frac{\chi-1}{\chi+1}}$,
                     \item \label{item 2: main theorem} for every $\alpha$-momentum \GD\ optimizer $\Phi_n\colon(\R^\fD)^n\to \R^\fD$, $n\in\N$, it holds that $(\Phi_n)_{n\in\N}$ converges with $\chi$-rate $\textcolor{magenta}{\frac{\sqrt{\chi}-1}{\sqrt{\chi}+1}}$,
        \item \label{item 3: main theorem} for every $\beta$-$\varepsilon$-\RMSprop\ optimizer $\Phi_n\colon(\R^\fD)^n\to\R^\fD$, $n\in\N$, it holds that $(\Phi_n)_{n\in\N}$ converges with $\chi$-rate $\textcolor{magenta}{\frac{\chi-1}{\chi+1}}$,
                      \item \label{item 4: main theorem} for every $\alpha$-$\beta$-$\varepsilon$-\Adam\ optimizer $\Phi_n\colon(\R^\fD)^n\to\R^\fD$, $n\in\N$, it holds that $(\Phi_n)_{n\in\N}$ converges with $\chi$-rate $\textcolor{magenta}{\frac{\sqrt{\chi}-1}{\sqrt{\chi}+1}}$, and
                       \item \label{item 5: main theorem} for every  $\alpha$-$d$-$\fd$-\MUON\ type optimizer $\Phi_n\colon(\R^\fD)^n\to\R^\fD$, $n\in\N$, it holds that $(\Phi_n)_{n\in\N}$ converges with $\chi$-rate $\textcolor{magenta}{\frac{\sqrt{\chi}-1}{\sqrt{\chi}+1}}$.
                 \end{enumerate}
                 (cf.\ Definitions \ref{definition: rate of convergence}, \ref{definition: momentum}, \ref{definition: muon type}, \ref{definition: GD}, \ref{definition: RMSprop}, and \ref{definition: Adam}).
                 \end{athm}
                 \end{tcolorbox}
                 \end{samepage}
\Cref{item 1: main theorem} in \cref{main theorem} is a direct consequence of \cref{cor: sharp rate rmsprop} in \cref{subsec: sharp rmsprop} below, \cref{item 2: main theorem} in \cref{main theorem} is a direct consequence of \cref{cor: sharp rate} in \cref{subsec: sharp muon} below, \cref{item 3: main theorem} in \cref{main theorem} is a direct consequence of \cref{cor: sharp rate rmsprop} in \cref{subsec: sharp rmsprop} below, \cref{item 4: main theorem} in \cref{main theorem} is a direct consequence of \cref{cor: sharp rate adam} in \cref{subsec: sharp adam} below, and \cref{item 5: main theorem} in \cref{main theorem} is a direct consequence of \cref{cor: sharp rate muon} in \cref{subsec: sharp muon} below.

\subsection{Literature review}

In this subsection we review selected articles in the scientific literature that are related to the findings of this work. Our work builds on the analysis of convergence rates for adaptive \GD\ optimization methods in \cite{DeArAdsharp2025}, and we refer to \cite[Subsection~1.3]{DeArAdsharp2025} for a review of selected articles on gradient-based optimization methods.

The upper bounds in \cref{item 1: main theorem,item 2: main theorem} in \cref{main theorem} are essentially due to Polyak \cite[Theorem~9]{Polyak1964SomeMO}, who proved that, in a neighborhood of a local minimizer at which the eigenvalues of the Hessian of the objective function lie in $[\convexcoe,\convexcoe\chi]$, the standard \GD\ optimizer (with a suitable constant learning rate) converges with rate $(\chi+1)^{-1}(\chi-1)$, and that the momentum \GD\ optimizer (with a suitable constant learning rate and a suitable momentum parameter) converges with the strictly faster rate $(\sqrt{\chi}+1)^{-1}(\sqrt{\chi}-1)$. The upper bounds in \cref{item 3: main theorem,item 4: main theorem} in \cref{main theorem} have essentially been established in \cite[Corollary~2.9 and Corollary~4.11]{DeArAdsharp2025}. We remark that, by \cite[Theorem~2.1.13]{Nesterov2004}, no first-order method whose iterates differ from the initial value by linear combinations of the previous gradients achieves a better rate than $(\sqrt{\chi}+1)^{-1}(\sqrt{\chi}-1)$ (uniformly over all smooth and strongly convex objective functions with condition number $\chi$ on an infinite-dimensional space) and, in this sense, this rate can not be improved within a large class of optimization methods. The main contribution of this work is the convergence rate for the \MUON\ optimizer in \cref{item 5: main theorem} in \cref{main theorem}. We note that, by \cref{definition: rate of convergence}, every optimizer converges with at most one $\chi$-rate and the rates in \cref{main theorem} can not be improved.

The \MUON\ optimizer has been recently introduced in \cite{jordan2024muon} and can be regarded as an accelerated (due to the involved momentum) and adaptive (due to the approximate matrix normalization, which adapts the learning rate separately for every singular direction of the momentum) \GD\ optimization method which, in addition, accounts for the matrix structure of the parameters to optimize. More specifically, in every step \MUON\ replaces the momentum matrix $M=U\Sigma V^\top$ of each parameter block by an approximation of its polar factor $UV^\top$ which, up to a scalar factor, is the direction of steepest descent with respect to the spectral norm, so that \MUON\ can be viewed as a computationally efficient realization of spectral descent (cf.\ \cite{BernLaker2024}). Moreover, the \MUON\ optimizer has been successfully employed in the training of large language models (cf., \eg, \cite{ar2505.02222,ar2507.20534,ar2502.16982}).

In practice, the polar factor is not computed exactly, but is approximated by finitely many steps of the \NS\ method (cf., \eg, \cite{MR295547} and \cite[Section~8.3]{MR2396439}) or of variants of it, such as the quintic iteration proposed in \cite{jordan2024muon} and the Polar Express method \cite{ar2505.16932} (cf.\ \cref{lem: abc verify 0}, \cref{lem: abc verify}, and \cref{lem: abc verify 2}). We emphasize that the analysis in our work does not rely on the approximation of the polar factor by the \NS\ iterations, but rather on the behavior of the generalized \NS\ functions (cf.\ \cref{definition: NS}) close to zero; indeed, due to the positive constant $\nscoe_{0,1}$ in the normalization step, the function $\NewtonSchulzAlgorithm{\nscoe}{K}$ is, to first order, the multiplication by $\bfb=(\nscoe_{0,1})^{-1}\prod_{k=1}^K\nscoe_{k,0}$ close to zero (cf.\ \cref{lem: NS represent}) and therefore, close to a local minimizer, \MUON\ behaves like the momentum \GD\ optimizer with a rescaled learning rate. In particular, our analysis also covers the \MUON$^p$ optimizer proposed in \cite{ar2606.13867}, which, for a rational $p\in(0,1)$, approximates a positive multiple of $U\Sigma^pV^\top$ instead of $UV^\top$; indeed, with a positive momentum parameter, a fixed number of iterations, and the normalization step in \cref{definition: NS} below, it is a \MUON\ type optimizer with one generalized \NS\ step (cf.\ \cref{definition: muon type} and \cite[Algorithm~1, Theorem~2.3, and Section~4.1]{ar2606.13867}).

There is also a growing literature on the convergence analysis of \MUON-type optimization methods. Many of these works establish, for stochastic or non-convex optimization problems, upper bounds for averages or minima of norms of the gradients for idealized \MUON, in which the \NS\ iterations are replaced by the exact polar factor (cf., \eg, \cite{ar2503.12645,muon1,Weimuonconvergence2025}). With constant learning rates, idealized \MUON\ does in general not converge to a minimizer, since its update does not become small when the gradient becomes small (cf., \eg, \cite{ar2602.11948,ar2505.21799}). In \cite{DoDeAr2026muon} an error analysis for \MUON\ (with the \NS\ polynomial from \cite{jordan2024muon} and with the Polar Express method as implemented in \cite{Jordanmuoncode}) has been established for stochastic optimization problems, and it has been shown that, for skewed data, \MUON\ does in general not converge to the minimizer. We also mention \cite{ar2610.02158}, where a regularized variant of idealized \MUON\ is obtained as the zero-temperature limit of a Langevin-based sampling algorithm. We refer to \cite[Subsection~1.5]{DoDeAr2026muon} for a more detailed review of the literature on the \MUON\ optimizer. 

The observation that the \NS\ functions act near zero as the multiplication by a positive constant, so that, close to a minimizer, \MUON\ behaves like the momentum \GD\ optimizer with a rescaled learning rate, has also been made in the recent work \cite{ar2610.34158} and, for a smooth surrogate of the \NS\ functions, in \cite{ar2609.30546}; these works, however, do not establish an explicit linear rate of convergence. For \MUON-type methods without momentum, local linear convergence with the rate of the standard \GD\ optimizer has been established for strongly convex quadratic problems in \cite{ar2610.27747}, where the case with momentum is left open. To the best of our knowledge, \cref{main theorem} is the first result in the scientific literature that establishes local linear convergence with a sharp rate for \MUON\ with momentum and with the \NS\ iterations that are employed in practice.

\subsection{Structure of this article}

The remainder of this article is organized as follows. In \cref{sec: muon} we establish optimal convergence rates for the momentum \GD\ and the \MUON\ optimizers and, thereby, prove \cref{item 2: main theorem,item 5: main theorem} in \cref{main theorem} above. In \cref{sec: rmsprop} we establish optimal convergence rates for the \GD\ and the \RMSprop\ optimizers and, thereby, prove \cref{item 1: main theorem,item 3: main theorem} in \cref{main theorem} above. In \cref{sec: adam} we establish optimal convergence rates for the \Adam\ optimizer and, thereby, prove \cref{item 4: main theorem} in \cref{main theorem} above.

\section{Optimal convergence rates of the momentum and the MUON optimizers}
\label{sec: muon}
In this section we establish in \cref{cor: sharp rate} and \cref{cor: sharp rate muon} below that the momentum \GD\ optimizer and every \MUON\ type optimizer with momentum parameter $\alpha=(\sqrt{\chi}+1)^{-2}(\sqrt{\chi}-1)^2$ converge with $\chi$-rate $(\sqrt{\chi}+1)^{-1}(\sqrt{\chi}-1)$. This is based on an upper bound (see \cref{theo: abstract optimal rate} below) and a lower bound (see \cref{main theorem lower bound} below) for the rate of convergence of a general class of optimization methods (which includes the momentum \GD\ optimizer and all \MUON\ type optimizers as special cases) in which the \NS\ functions are replaced by continuously differentiable functions which vanish at zero and whose derivatives at zero are all equal to the same positive multiple of the identity (cf.\ \cref{lem: NS represent} and \cref{definition: muon type}).

\subsection{Momentum optimizer}
In the next \cref{definition: momentum} we recall the concept of the momentum \GD\ optimizer (see \cite{Polyak1964SomeMO}) using the functions $\Phi_n\colon(\R^d)^n\to\R^d$, $n\in\N$, from \cref{definition: rate of convergence} above; see, \eg, \cite[Definition~3.1]{DeArAdsharp2025}.

\begin{tcolorbox}[colback=white!95!gray,
                  colframe=black,
                  boxrule=0.5pt,
                  sharp corners,
                  enhanced,
                  breakable,
                 ]
 \begin{definition}[\textcolor{red}{Momentum \GD}]
    \label{definition: momentum}
    Let $ d \in \N $, $ \alpha \in [0,1) $
    and let
    $
      \Phi_n = ( \Phi_n^1, \dots, \Phi_n^{ d } )
      \colon ( \R^{ d } )^n
      \allowbreak
      \to \R^{ d }
    $,
    $ n \in \N $, be functions. Then we say that
    $ ( \Phi_n )_{ n \in \N } $ is the
    $ \alpha $-momentum \GD\ optimizer
    on $ \R^{ d } $ (we say that
    $ ( \Phi_n )_{ n \in \N } $ is
    the $ \alpha $-momentum \GD\ optimizer)
    if and only if it holds for all
    $ n \in \N $,
    $
      g =
      (
        ( g_{ i, j } )_{ j \in \{ 1, 2, \dots, d \} }
      )_{ i \in \{ 1, 2, \dots, n \} } \in ( \R^{ d } )^n
    $,
    $ j \in \{ 1, 2, \dots, d \} $
    that
    \begin{equation}
    \textstyle
	  \Phi_n^j( g )
	  =
	  ( 1 - \alpha )
	  \bigl[
	  \sum\nolimits_{ i = 1 }^n
	  \alpha^{ n - i }
	  g_{ i, j }
	  \bigr]
	  .
    \end{equation}
	\end{definition}
    \end{tcolorbox}

\subsection{Generalized Newton-Schulz (NS) methods and MUON optimizer}
In the next \cref{definition: NS} we recall the concept of the generalized \NS\ method (cf.\ \cite[Definition~1.1]{DoDeAr2026muon}), which is a direct generalization of the classical \NS\ method (cf., \eg, \cite{MR295547}, \cite[Subsection~1.1]{DoDeAr2026muon}, \cite[Section~8.3]{MR2396439}, and \cite[Subsection~1.2]{ar2505.16932}) applied to a regularized normalization of the input matrix.

\begin{tcolorbox}[colback=white!95!gray,
                  colframe=black,
                  boxrule=0.5pt,
                  sharp corners,
                  enhanced,
                  breakable,
                 ]
\begin{definition}[\textcolor{red}{Generalized \NS\ method}]\label{definition: NS}
   Let $\nscoe=(\nscoe_{i,j})_{(i,j)\in(\N_0)^2}\colon (\N_0)^2\to\R$ satisfy $\min\{\nscoe_{0,0},\nscoe_{0,1}\}>0$ and $\#(\nscoe^{ - 1 }( \R\backslash\{0\} )) <\infty$ and let $n,m\in \N$, $A\in \R^{n\times m}$. Then we denote by
$
	\NewtonSchulzAlgorithm{\nscoe}{k}\colon \R^{n \times m}\to\R^{n\times m}
$, $k\in \N_0$, the sequence of functions
 which satisfy for all
	$k \in \N$, $ A,M \in \R^{n\times m}$ with $\textcolor{magenta}{M=\TNewtonSchulzAlgorithm{\nscoe}{k-1}{A}} $
that\footnotemark\
$\textcolor{magenta}{\TNewtonSchulzAlgorithm{\nscoe}{0}{A} = (\nscoe_{0,0}\|A\| + \nscoe_{0,1})^{-1} A}$
and
\begin{equation}\label{T_B_D}
\begin{split}
\textcolor{magenta}{\TNewtonSchulzAlgorithm{\nscoe}{k}{A}= \textstyle\sum_{ i = 0 }^{ \infty } \nscoe_{ k, i } ( M M^{\top} )^{ i} M}
.
\end{split}
\end{equation}
\end{definition}
\end{tcolorbox}
\footnotetext{Note that for all $N\in \N$, $v_1, v_2, \dots, v_{N }, w_1, w_2, \dots, w_{ N } \in \N$, $\theta=( ( \theta^n_{ i, j } )_{ (i, j) \in \{1,2,\dots,v_n\}\times\{1,2,\dots,w_n\} } )_{ n \in \{1,2,\dots,N\} } \in \times_{ n = 1 }^{N} \R^{ v_n \times w_n }$ it holds that $\| \theta \| = \bigl[\sum_{n=1}^N\sum_{i=1}^{v_n}\sum_{j=1}^{w_n}|\theta_{i,j}^n|^2\bigr]^{1/2}$ (Frobenius norm).}

In the next \cref{definition: muon} we reformulate the concept of the \MUON\ optimizer with generalized \NS\ iterations (cf.\ \cite{jordan2024muon}, \cite{ar2505.16932}, and \cite[Definition~3.14]{DoDeAr2026muon}) using the functions $\Phi_n\colon(\R^\fD)^n\to\R^\fD$, $n\in\N$, from \cref{definition: rate of convergence} above.

\begin{tcolorbox}[colback=white!95!gray,
                  colframe=black,
                  boxrule=0.5pt,
                  sharp corners,
                  enhanced,
                  breakable,
                 ]
\begin{definition}[\textcolor{red}{\MUON\ optimizer}]\label{definition: muon}
Let $\nb\in\N$, let $d=(d_1,\dots,d_\nb)$, $\fd=(\fd_1,\dots,\fd_\nb)\in \N^\nb$, let $\cd=\times_{i=1}^{\nb}\R^{d_i\times\fd_i}$, let $\fD,D_1,D_2,\dots,D_\nb\in\N$ satisfy for all $k\in\{1,2,\dots,\nb\}$ that $D_k=d_k\fd_k$ and $\fD=\sum_{i=1}^\nb D_i$, let $\operatorname{vec}\colon\cd\to\R^\fD$ satisfy for all $A=((a_{l,i,j})_{(i,j)\in \{1,2,\dots,d_l\}\times\{1,2,\dots,\fd_l\}})_{l\in\{1,2,\dots,\nb\}}\in\cd$ and all $v_1=(v_{1,1},\dots,v_{1,D_1})\in\R^{D_1}$, $v_2=(v_{2,1},\dots,v_{2,D_2})\in\R^{D_2}$, $\dots$, $v_\nb=(v_{\nb,1},\dots,v_{\nb,D_\nb})\in\R^{D_\nb}$ with
\begin{equation}
    \forall\,l\in\{1,2,\dots,\nb\}\colon v_l=(a_{l,1,1},\dots,a_{l,d_l,1}, a_{l,1,2},\dots,a_{l,d_l,2}, \dots,a_{l,1,\fd_l},\dots,a_{l,d_l,\fd_l})
\end{equation}
that
\begin{equation}
    \operatorname{vec}(A)=(v_{1,1},\dots,v_{1,D_1},v_{2,1},\dots,v_{2,D_2},\dots,v_{\nb,1},\dots,v_{\nb,D_\nb}),
\end{equation}
let $\operatorname{mat}\colon\R^\fD\to\cd$ satisfy for all $v\in \R^{\sum_{i=1}^\nb D_i}$ that $\operatorname{vec}(\operatorname{mat}(v))=v$, let $\alpha\in [0,1)$, $K\in\N_0$, let $\nscoe=(\nscoe_{i,j})_{(i,j)\in(\N_0)^2}\colon (\N_0)^2\to\R$ satisfy $\min\{\nscoe_{0,0},\nscoe_{0,1}\}>0$ and $\#(\nscoe^{ - 1 }( \R\backslash\{0\} )) <\infty$, assume for all $i\in \mathbb{N}\cap[0,K]$ that $\nscoe_{i,0}>0$, and let $\Phi_n\colon(\R^\fD)^n\to\R^\fD$, $n\in\N$, be functions. Then we say that $(\Phi_n)_{n\in\N}$ is an $\alpha$-$\nscoe$-$K$-$d$-$\fd$-\MUON\ optimizer on $\R^{\fD}$ (we say that $(\Phi_n)_{n\in\N}$ is the $\alpha$-$\nscoe$-$K$-$d$-$\fd$-\MUON\ optimizer) if and only if it holds for all $n\in\N$, $g_1,g_2,\dots,g_n\in \R^\fD$, $(M_1,M_2,\dots,M_\nb)\in\cd$ with
\begin{equation}
    (M_1,M_2,\dots,M_\nb)=\operatorname{mat}\bigl(\textstyle (1-\alpha)\sum_{k=1}^n\alpha^{n-k}g_k\bigr)
\end{equation}
that
\begin{equation}
    \Phi_n(g_1,g_2,\dots,g_n)=\operatorname{vec}(\NewtonSchulzAlgorithm{\nscoe}{K}(M_1),\dots,\NewtonSchulzAlgorithm{\nscoe}{K}(M_\nb))
\end{equation}
(cf.\ \cref{definition: NS}).
\end{definition}
\end{tcolorbox}
In the next \cref{definition: muon type} we introduce the class of \MUON\ type optimizers, which appears in \cref{item 5: main theorem} in \cref{main theorem} above.
\begin{tcolorbox}[colback=white!95!gray,
                  colframe=black,
                  boxrule=0.5pt,
                  sharp corners,
                  enhanced,
                  breakable,
                 ]
                 \begin{definition}[\textcolor{red}{\MUON\ type optimizer}]\label{definition: muon type}
                 Let $\fD,\nb\in\N$, $\alpha\in (0,1)$, $d=(d_1,\dots,d_\nb)$, $\fd=(\fd_1,\dots,\fd_\nb)\in \N^\nb$ satisfy $\fD=\sum_{i=1}^\nb d_i\fd_i$ and let $\Phi_n\colon(\R^\fD)^n\to\R^\fD$, $n\in\N$, be functions. Then we say that $(\Phi_n)_{n\in\N}$ is an $\alpha$-$d$-$\fd$-\MUON\ type optimizer on $\R^{\fD}$ if and only if there exist $K\in\N_0$ and $\nscoe=(\nscoe_{i,j})_{(i,j)\in(\N_0)^2}\colon (\N_0)^2\to\R$ such that
                 \begin{enumerate}[label=(\roman*)]
                     \item it holds that $\min\{\nscoe_{0,0},\nscoe_{0,1}\}>0$ and $\#(\nscoe^{ - 1 }( \R\backslash\{0\} )) <\infty$,
                     \item it holds for all $i\in \mathbb{N}\cap[0,K]$ that $\nscoe_{i,0}>0$, and
                     \item it holds that $(\Phi_n)_{n\in\N}$ is an $\alpha$-$\nscoe$-$K$-$d$-$\fd$-\MUON\ optimizer on $\R^{\fD}$ 
                 \end{enumerate}
                 (cf.\ \cref{definition: muon}).
                 \end{definition}
                 \end{tcolorbox}
\subsection{Upper bounds on convergence rates of momentum and MUON}
\begin{samepage}
 \begin{tcolorbox}[colback=white!95!gray,
                  colframe=black,
                  boxrule=0.5pt,
                  sharp corners,
                  enhanced,
                  breakable,
                 ]
\begin{athm}{theorem}{theo: muon higher bound}
      Let $\nb\in \mathbb{N}$,  $\chi\in (1,\infty)$, $\convexcoe,\bfb\in (0,\infty)$, $\varrho\in (\frac{\sqrt{ \chi } - 1 }  { \sqrt{ \chi } + 1 },1)$, $\Gamma=(\bfb\convexcoe\sqrt{\chi})^{-1}$, $\alpha=(\frac{\sqrt{ \chi } - 1 }  { \sqrt{ \chi } + 1 })^2$, for every $i\in\{1,2,\dots,\nb\}$ let $d_i,\fd_i\in\N$, $\Phi_i\in C^1(\R^{d_i\times\fd_i},\R^{d_i\times\fd_i})$  satisfy for all $\theta\in \R^{d_i\times\fd_i}$ that $(\Phi_i)'(0)\theta-\bfb\theta=0=\Phi_i(0)$, let $\cd=\times_{i=1}^{\nb}\R^{d_i\times\fd_i}$, $\vartheta\in\cd$, and let $\smalll\in \allowbreak C^{2}(\cd,\mathbb{R})$ satisfy for all $\theta \in \cd$ that \begin{equation}\llabel{def: bbL}
     \textcolor{magenta}{(\nabla\smalll)(\vartheta)=0}\qqandqq\textcolor{magenta}{\convexcoe\|\theta\|^2\leq \smalll'' (\vartheta)(\theta,\theta)\leq \convexcoe\chi \| \theta \|^2}.
 \end{equation}  
 Then 
there exist $\fC,\delta \in (0,\infty)$ such that for every $\gamma\colon\N\to(\Gamma-\delta,\Gamma+\delta)$, every $\bfm=(\bfm^1,\dots,\bfm^{\nb})\colon \mathbb{N}_0\to\cd$, and every $\Theta=(\Theta^1,\dots,\Theta^{\nb})\colon \mathbb{N}_0\to\cd$ with the property that for all $n \in \mathbb{N}$, $i\in\{1,2,\dots,\nb\}$ it holds that 
         \begin{equation}
            \textcolor{magenta}{\bfm_0=0},\qquad \textcolor{magenta}{\bfm_n= \alpha \bfm_{n-1}+(1-\alpha)(\nabla\smalll)(\Theta_{n-1})},
            \end{equation}
            \begin{equation}
            \textcolor{magenta}{\|\Theta_0-\vartheta\|< \delta},\qquad\text{and}\qquad  \textcolor{magenta}{\Theta_n^i=\Theta_{n-1}^i-\gamma_n
          \Phi_i(\bfm_n^i)}
         \end{equation}
          we have for all $n\in\N_0$ that
$    \textcolor{magenta}{ \textstyle
    \|\Theta_n-\vartheta\|\leq \fC\|\Theta_0-\vartheta\|\varrho^{n}}$.
\end{athm}
\end{tcolorbox}
\end{samepage}
\begin{aproof}
    Throughout this proof let $\cL=(\cL(\theta))_{\theta=(\theta_1,\dots,\theta_\nb)\in \cd}\colon \cd\to\R$ satisfy for all $\theta\in \cd$ that $\cL(\theta)=\smalll(\theta)$, let $\bfv=(\bfv_1,\dots,\bfv_\nb)\in\cd$ satisfy $\bfv=\vartheta$,  for every $\xi\in\cd$ and every $\gamma\colon\N\to\R$ let $\Theta^{\xi,\gamma}=(\Theta^{\xi,\gamma,1},\dots,\Theta^{\xi,\gamma,\nb})\colon\N_0\to\cd$ and $\bfm^{\xi,\gamma}=(\bfm^{\xi,\gamma,1},\dots,\bfm^{\xi,\gamma,\nb})\colon\N_0\to\cd$ satisfy for all $n\in\N$, $i\in\{1,2,\dots,\nb\}$ that
     \begin{equation}\llabel{def: bfm}
      \bfm^{\xi,\gamma}_0=0,\qquad \bfm^{\xi,\gamma}_n= \alpha \bfm^{\xi,\gamma}_{n-1}+(1-\alpha)(\nabla\cL)(\Theta^{\xi,\gamma}_{n-1}),
            \end{equation}
            \begin{equation}\llabel{def: Theta}
        \Theta^{\xi,\gamma}_0=\xi,\qquad\text{and}\qquad \Theta^{\xi,\gamma,i}_n=\Theta^{\xi,\gamma,i}_{n-1}-\gamma_n
           \Phi_i(\bfm^{\xi,\gamma,i}_n),
         \end{equation}
         let $\cA= (\cA^{1},\cA^{2})=((\cA^{1,1},\dots,\cA^{1,\nb}),\cA^2)\colon \cd\times\cd\to\cd\times\cd$ satisfy for all $\theta,m\in\cd$ that
         \begin{equation}\llabel{def: cA}
             \cA^{1}(\theta,m)=\theta-\Gamma \bfb(\alpha m+(1-\alpha)(\nabla\cL)(\theta+\bfv))\qqandqq \cA^{2}(\theta,m)=\alpha m+(1-\alpha)(\nabla\cL)(\theta+\bfv), 
         \end{equation} 
         and for every $n\in\N$ and every $\gamma\colon\N\to\R$ let $\psi_{\gamma,n}= (\psi_{\gamma,n}^{1},\psi_{\gamma,n}^{2})=((\psi_{\gamma,n,1}^{1},\dots,\psi_{\gamma,n,\nb}^{1}),\psi_{\gamma,n}^{2})\colon \cd\times\cd\to\cd\times\cd$ satisfy for all $\theta=(\theta_1,\dots,\theta_\nb)$, $m=(m_1,\dots,m_\nb)\in\cd$, $i\in\{1,2,\dots,\nb\}$ that
         \begin{equation}\llabel{def: Phi}
             \psi_{\gamma,n,i}^{1}(\theta,m)=\theta_i-\gamma_n\Phi_i(\alpha m_i+(1-\alpha)(\nabla_{\theta_i}\cL)(\theta+\bfv))\qandq \psi_{\gamma,n}^{2}(\theta,m)=\alpha m+(1-\alpha)(\nabla\cL)(\theta+\bfv).
         \end{equation} 
         \argument{\lref{def: bfm};\lref{def: Theta};\lref{def: Phi}}{that for all $n\in \N$, $\xi\in \cd$, $i\in\{1,2,\dots,\nb\}$ and all $\gamma\colon\N\to\R$ it holds that
         \begin{equation}\llabel{eq1}
         \begin{split}
            & \Theta_n^{\xi,\gamma,i}-\bfv_i=\Theta_{n-1}^{\xi,\gamma,i}-\bfv_i-\gamma_n\Phi_i(\alpha \bfm^{\xi,\gamma,i}_{n-1}+(1-\alpha)(\nabla_{\theta_i}\cL)(\Theta^{\xi,\gamma}_{n-1}))\\
             &=\Theta_{n-1}^{\xi,\gamma,i}-\bfv_i-\gamma_n\Phi_i(\alpha \bfm^{\xi,\gamma,i}_{n-1}+(1-\alpha)(\nabla_{\theta_i}\cL)(\Theta^{\xi,\gamma}_{n-1}-\bfv+\bfv))=\psi_{\gamma,n,i}^{1}(\Theta_{n-1}^{\xi,\gamma}-\bfv,\bfm^{\xi,\gamma}_{n-1})
             \end{split}
         \end{equation}
         and
         \begin{equation}\llabel{eq2}
             \bfm_n^{\xi,\gamma}=\alpha \bfm^{\xi,\gamma}_{n-1}+(1-\alpha)(\nabla\cL)(\Theta^{\xi,\gamma}_{n-1}-\bfv+\bfv)=\psi_{\gamma,n}^{2}(\Theta^{\xi,\gamma}_{n-1}-\bfv,\bfm^{\xi,\gamma}_{n-1}).
         \end{equation}}
         \argument{\lref{eq2};}{for all $n\in\N$, $\xi\in \cd$ and all $\gamma\colon\N\to\R$ that
         \begin{equation}\llabel{eq3}
         \begin{split}
             (\Theta_n^{\xi,\gamma}-\bfv,\bfm_n^{\xi,\gamma})&=\psi_{\gamma,n}(\Theta^{\xi,\gamma}_{n-1}-\bfv,\bfm^{\xi,\gamma}_{n-1})\\
             &=\cA(\Theta^{\xi,\gamma}_{n-1}-\bfv,\bfm^{\xi,\gamma}_{n-1})+\psi_{\gamma,n}(\Theta^{\xi,\gamma}_{n-1}-\bfv,\bfm^{\xi,\gamma}_{n-1})-\cA(\Theta^{\xi,\gamma}_{n-1}-\bfv,\bfm^{\xi,\gamma}_{n-1}).
             \end{split}
         \end{equation}}
         \argument{\lref{def: cA};}{for all $\theta_1,\theta_2\in \cd$ that 
        \begin{equation}\llabel{eq5}
            (\cA^1)'(0)(\theta_1,\theta_2)=\theta_1-\Gamma \bfb (\alpha\theta_2+(1-\alpha)(\operatorname{Hess} \cL)(\bfv)\theta_1)
        \end{equation}
        and 
        \begin{equation}\llabel{eq6}
            (\cA^2)'(0)(\theta_1,\theta_2)=\alpha\theta_2+(1-\alpha)(\operatorname{Hess} \cL)(\bfv)\theta_1.
        \end{equation}}
         \argument{the fact that $\alpha=(\frac{\sqrt{\chi}-1}{\sqrt{\chi}+1})^2$}{that
    \begin{equation}\llabel{eqq0}
        \textstyle\chi=(\frac{1+\sqrt{\alpha}}{1-\sqrt{\alpha}})^2.
    \end{equation}}
    \argument{\lref{eqq0};the assumption that $\Gamma=(\bfb \convexcoe\sqrt{\chi})^{-1}$;}{that
    \begin{equation}\llabel{eqq1}
    \begin{split}
          &(1+\alpha-\Gamma\bfb \convexcoe(1-\alpha))^2= (1+\alpha-\chi^{-1/2}(1-\alpha))^2=\biggl(1+\alpha-\frac{(1-\alpha)(1-\sqrt{\alpha})}{1+\sqrt{\alpha}}\biggr)^2=4\alpha\\
          &\biggl(1+\alpha-\frac{(1-\alpha)(1+\sqrt{\alpha})}{1-\sqrt{\alpha}}\biggr)^2=(1+\alpha-\chi^{1/2}(1-\alpha))^2=(1+\alpha-\Gamma\bfb \convexcoe\chi(1-\alpha))^2.
          \end{split}
    \end{equation}}
    \argument{\lref{eqq1};the fact that $\Gamma\bfb>0$}{that for all $\lambda\in [\convexcoe,\convexcoe\chi]$ it holds that
    \begin{equation}\llabel{eqq1'}
        (1+\alpha-\Gamma\bfb \lambda(1-\alpha))^2\leq 4\alpha.
    \end{equation}}
       \argument{\lref{def: bbL};}{that for all $\lambda\in \R$, $\theta\in (\times_{i=1}^{\nb}\bbC^{d_i\times\fd_i})\backslash\{0\}$ with $(\operatorname{Hess} \cL)(\bfv)\theta=\lambda\theta$ it holds that
        \begin{equation}\llabel{eq7}
\convexcoe\leq \lambda\leq \convexcoe\chi.
        \end{equation}}
         \argument{\lref{eq7}; \lref{eqq1'}}
      {that for all $\lambda\in \R$, $\theta\in (\times_{i=1}^{\nb}\bbC^{d_i\times\fd_i})\backslash\{0\}$ with $(\operatorname{Hess} \cL)(\bfv)\theta=\lambda\theta$ it holds that
        \begin{equation}\llabel{eqq1''}
          \textstyle\bigl(  \frac{1+\alpha-\Gamma\bfb\lambda(1-\alpha)}{2}\bigr)^2+\bigl(  \frac{1+\alpha-\Gamma\bfb\lambda(1-\alpha)}{2}\bigr)^2-\alpha=\frac{(1+\alpha-\Gamma\bfb\lambda(1-\alpha))^2}{2}-\alpha\leq \alpha.
        \end{equation}}
\argument{\lref{eq5};\lref{eq6};\lref{eqq1''};\cite[item (i) in Lemma 3.2]{DeArAdsharp2025} (applied with $d\curvearrowleft\sum_{i=1}^\nb d_i\fd_i$, $\alpha\curvearrowleft\alpha$, $\gamma\curvearrowleft \Gamma\bfb$, $H\curvearrowleft(\operatorname{Hess}\cL)(\bfv)$, $A\curvearrowleft \cA'(0)$ in the notation of \cite[Lemma 3.2]{DeArAdsharp2025})}{that for all $\lambda\in \bbC$, $\theta\in ((\times_{i=1}^{\nb}\bbC^{d_i\times\fd_i})\times(\times_{i=1}^{\nb}\bbC^{d_i\times\fd_i}))\backslash\{0\}$ with $\cA'(0)\theta=\lambda\theta$ it holds that
\begin{equation}\llabel{eq8}
    |\lambda|\leq \sqrt{\alpha}.
\end{equation}}
\argument{\lref{eq8};the fact that $\varrho>\sqrt{\alpha}$;\unskip, \eg, \cite[Lemma 2.3]{DeArAdsharp2025}}{that there exists a norm $\nnorm{\cdot}\colon \cd\times\cd\to[0,\infty)$ which satisfies for all $\theta\in\cd\times\cd$ that
\begin{equation}\llabel{eq9}
    \nnorm{\cA'(0)\theta}\leq \frac{\varrho+\sqrt{\alpha}}{2}\nnorm{\theta}.
\end{equation}}
      \startnewargseq
       \argument{the fact that $\|\cdot\|$ is a norm; the fact that $\nnorm{\cdot}$ is a norm; the fact that $\cd\times\cd$ is finite dimensional vector space;}{that there exist $\mu_1,\mu_2\in (0,\infty)$ which satisfy for all $\theta\in \cd\times\cd$ that
         \begin{equation}\llabel{eq14}
             \mu_1\nnorm{\theta}\leq \|\theta\|\leq \mu_2\nnorm{\theta}.
         \end{equation}}
         \startnewargseq
         \argument{\lref{def: cA};the fact that $\cL\in C^2(\cd,\R)$}{that \llabel{arg1} $\cA\in C^1(\cd\times\cd,\cd\times\cd)$\dott}
         \argument{\lref{arg1};\lref{eq14}; the fact that $\varrho>\sqrt{\alpha}$}{that there exists $\delta\in (0,\infty)$ such that for all $\theta\in \cd\times\cd$ with $\nnorm{\theta}\leq \delta$ it holds that
         \begin{equation}\llabel{eq4}
         \begin{split}
             \nnorm{\cA(\theta)-\cA(0)-\cA'(0)\theta}\leq  (\mu_1)^{-1}\|\cA(\theta)-\cA(0)-\cA'(0)\theta\|\leq \frac{\varrho-\sqrt{\alpha}}{4\mu_2}\|\theta\|\leq \frac{\varrho-\sqrt{\alpha}}{4}\nnorm{\theta}.
             \end{split}
         \end{equation}}
         \argument{\lref{def: bbL};\lref{def: cA}; the fact that $\vartheta=\bfv$}{that
         \begin{equation}\llabel{eq10}
             \cA(0)=(0-\Gamma\bfb(0\alpha+(1-\alpha)(\nabla\cL)(\bfv)),0\alpha+(1-\alpha)(\nabla\cL)(\bfv))=0.
         \end{equation}}
         \argument{\lref{eq10};\lref{eq4};}{that there exists $\delta\in (0,\infty)$ such that for all $\theta\in \cd\times\cd$ with $\nnorm{\theta}\leq \delta$ it holds that
         \begin{equation}\llabel{evd1}
             \nnorm{\cA(\theta)-\cA'(0)\theta}\leq \frac{\varrho-\sqrt{\alpha}}{4}\nnorm{\theta}.
         \end{equation}}
               \argument{\lref{def: bbL};the fact that $\cL\in C^2(\cd,\R)$}{that there exist $\delta,\fC_1,\fC_2\in (0,\infty)$ such that for all $\theta\in \cd$ with $\|\theta-\bfv\|\leq \delta$ it holds that
         \begin{equation}\llabel{eq15}
            \fC_1\|\theta-\bfv\|\leq \|(\nabla\cL)(\theta)\|=\|(\nabla\cL)(\theta)-(\nabla\cL)(\bfv)\|\leq \fC_2\|\theta-\bfv\|.
         \end{equation}}
         \argument{\lref{eq15};\lref{eq14}}{that there exist $\delta,\fC_1,\fC_2,\fC_3\in (0,\infty)$ such that for all $\theta,m\in \cd$ with $\|\theta-\bfv\|\leq \delta$ it holds that
         \begin{equation}\llabel{def: nu1'}
         \begin{split}
             \|\alpha m+(1-\alpha)(\nabla\cL)(\theta)\|&\leq \alpha \|m\|+(1-\alpha)\|(\nabla\cL)(\theta)\|\leq \|m\|+\fC_1\|\theta-\bfv\|
             \leq \fC_2\|(\theta-\bfv,m)\|\\
             &\leq \fC_3\nnorm{(\theta-\bfv,m)}.
             \end{split}
         \end{equation}}
\argument{\lref{eq14};}{that for all $\theta,m\in \cd$ it holds that
\begin{equation}\llabel{eqqq1}
    \|\theta-\bfv\|\leq \|(\theta-\bfv,m)\|\leq \mu_2\nnorm{(\theta-\bfv,m)}.
\end{equation}}
\argument{\lref{def: nu1'};\lref{eqqq1};}{that there exist $\nu,\fC\in (0,\infty)$ which satisfy for all $\theta,m\in \cd$ with $\nnorm{(\theta-\bfv,m)}\leq \nu$ it holds that
         \begin{equation}\llabel{def: nu}
         \begin{split}
             \|\alpha m+(1-\alpha)(\nabla\cL)(\theta)\|\leq \fC\nnorm{(\theta-\bfv,m)}.
             \end{split}
         \end{equation}}
        \startnewargseq
         \argument{the fact that for all $i\in\{1,2,\dots,\nb\}$ it holds that $\Phi_i\in C^1(\R^{d_i\times\fd_i},\R^{d_i\times\fd_i})$; the fundamental theorem of calculus}{that for all $i\in\{1,2,\dots,\nb\}$, $\theta\in \R^{d_i\times\fd_i}$  it holds that
         \begin{equation}\llabel{eq12.1}
             \begin{split}
                 \Phi_i(\theta)-\Phi_i(0)=\int_{0}^1(\Phi_i)'(t\theta)\theta\,\d t=(\Phi_i)'(0)\theta+\int_{0}^1((\Phi_i)'(t\theta)-(\Phi_i)'(0))\theta\,\d t.
             \end{split}
         \end{equation}}
         \startnewargseq
         \argument{\lref{eq12.1};the fact that for all $i\in\{1,2,\dots,\nb\}$ it holds that $\Phi_i(0)=0$; the fact that for all $i\in\{1,2,\dots,\nb\}$, $\theta\in\R^{d_i\times\fd_i}$ with it holds that $(\Phi_i)'$ is continuous;the fact that $\varrho>\sqrt{\alpha}$}{that there exists $\delta\in (0,\infty)$ such that for all $i\in\{1,2,\dots,\nb\}$, $\theta\in\R^{d_i\times\fd_i}$ with $\|\theta\|\leq \delta$ it holds that
         \begin{equation}\llabel{eq12.2}
             \|(\Phi_i)(\theta)-(\Phi_i)'(0)\theta\|= \|\Phi_i(\theta)-\Phi_i(0)-(\Phi_i)'(0)\theta\|\leq \int_{0}^1\frac{\mu_1(\varrho-\sqrt{\alpha})}{16\fC\Gamma\nb}\|\theta\|\,\d t=\frac{\mu_1(\varrho-\sqrt{\alpha})}{16\fC\Gamma\nb}\|\theta\|.
         \end{equation}}
         \argument{\lref{eq12.2};the fact that for all $i\in\{1,2,\dots,\nb\}$,  $\theta\in \R^{d_i\times\fd_i}$ it holds that $(\Phi_i)'(0)\theta=\bfb\theta$}{that there exists $\delta\in (0,\infty)$ such that for all $i\in\{1,2,\dots,\nb\}$, $\theta\in \R^{d_i\times\fd_i}$ with $\|\theta\|\leq \delta$ it holds that
         \begin{equation}\llabel{eq12}
             \|\Phi_i(\theta)-\bfb\theta\|\leq \frac{\mu_1(\varrho-\sqrt{\alpha})}{16\fC\Gamma\nb}\|\theta\|.
         \end{equation}}
         \argument{\lref{eq12};\lref{def: cA};\lref{def: Phi};}{that there exists $\delta\in (0,1)$ such that for every $n\in\N$, $\theta=(\theta_{1},\dots,\theta_{\nb})$, $m=(m_{1},\allowbreak\dots,m_{\nb})\in \cd$ and every $\gamma\colon\N\to\R$ with $\|\alpha m+(1-\alpha)(\nabla\cL)(\theta+\bfv)\|\leq \delta$ and $\sup_{k\in\N}|\gamma_k|\leq 2\Gamma$ it holds that
         \begin{equation}\llabel{eq13'}
         \begin{split}
             &\|\psi_{\gamma,n}(\theta,m)-\cA(\theta,m)\|\\
             &\textstyle\leq\bigl[ \sum_{i=1}^{\nb}\|\psi_{\gamma,n,i}^1(\theta,m)-\cA^{1,i}(\theta,m)\|\bigr]+\|\psi_{\gamma,n}^2(\theta,m)-\cA^2(\theta,m)\|\\
             &\textstyle=\bigl[ \sum_{i=1}^{\nb}\|\psi_{\gamma,n,i}^1(\theta,m)-\cA^{1,i}(\theta,m)\|\bigr]\\
&=\textstyle\sum_{i=1}^{\nb}\|\gamma_n\Phi_i(\alpha m_i+(1-\alpha)(\nabla_{\theta_i}\cL)(\theta+\bfv))-\Gamma\bfb(\alpha m_i+(1-\alpha)(\nabla_{\theta_i}\cL)(\theta+\bfv))\|\\
             &\leq |\gamma_n|\textstyle\sum_{i=1}^{\nb}\|\Phi_i(\alpha m_i+(1-\alpha)(\nabla_{\theta_i}\cL)(\theta+\bfv))-\bfb(\alpha m_i+(1-\alpha)(\nabla_{\theta_i}\cL)(\theta+\bfv))\|\\
             &+|\gamma_n-\Gamma|\textstyle\sum_{i=1}^{\nb}\|\bfb(\alpha m_i+(1-\alpha)(\nabla_{\theta_i}\cL)(\theta+\bfv))\|\\
             &\leq \frac{\mu_1(\varrho-\sqrt{\alpha})}{16\fC\Gamma\nb}|\gamma_n|\textstyle\sum_{i=1}^{\nb}\|\alpha m_i+(1-\alpha)(\nabla_{\theta_i}\cL)(\theta+\bfv)\|\\
             &+|\gamma_n-\Gamma|\bfb\textstyle\sum_{i=1}^{\nb}\|\alpha m_i+(1-\alpha)(\nabla_{\theta_i}\cL)(\theta+\bfv)\|\\
             &\leq \frac{\mu_1(\varrho-\sqrt{\alpha})}{8\fC}\|\alpha m+(1-\alpha)(\nabla\cL)(\theta+\bfv)\|+\bfb|\gamma_n-\Gamma|\nb\|\alpha m+(1-\alpha)(\nabla\cL)(\theta+\bfv)\|.
             \end{split}
         \end{equation}}
         \argument{\lref{eq13'};\lref{eq14};\lref{def: nu}}{that there exists $\delta\in (0,1)$ such that for every $n\in\N$, $\theta,m\in \cd$ and every $\gamma\colon\N\to\R$ with $\|\alpha m+(1-\alpha)(\nabla\cL)(\theta+\bfv)\|\leq \delta$, $\nnorm{(\theta,m)}\leq \delta$, and $\sup_{k\in\N}|\gamma_k-\Gamma|\leq \delta$ it holds that
         \begin{equation}
             \llabel{eq13''}
             \begin{split}
             &\nnorm{\psi_{\gamma,n}(\theta,m)-\cA(\theta,m)}\leq (\mu_1)^{-1}\|\psi_{\gamma,n}(\theta,m)-\cA(\theta,m)\|\\
             &\leq \biggl( \frac{\varrho-\sqrt{\alpha}}{8}+(\mu_1)^{-1}\bfb|\gamma_n-\Gamma| \nb\fC\biggr)\nnorm{(\theta,m)}\leq\frac{\varrho-\sqrt{\alpha}}{4}\nnorm{(\theta,m)}.
             \end{split}
         \end{equation}}
         \argument{\lref{eq13''};\lref{def: nu}}{that there exists $\delta\in (0,1)$ such that for every $n\in\N$, $\theta,m\in \cd$ and every $\gamma\colon\N\to\R$ with $\nnorm{(\theta,m)}\leq \delta$ and $\sup_{k\in\N}|\gamma_k-\Gamma|\leq \delta$ it holds that
         \begin{equation}
             \llabel{eq13}
             \begin{split}
             \nnorm{\psi_{\gamma,n}(\theta,m)-\cA(\theta,m)}\leq \frac{\varrho-\sqrt{\alpha}}{4}\nnorm{(\theta,m)}.
             \end{split}
         \end{equation}}
         \argument{\lref{eq13};}{that there exists $\delta\in (0,1)$ such that for every $n\in\N$, $\xi\in\cd$ and every $\gamma\colon\N\to\R$ with $\nnorm{(\Theta_{n-1}^{\xi,\gamma}-\bfv,\bfm_{n-1}^{\xi,\gamma})}\leq \delta$ and $\sup_{k\in\N}|\gamma_k-\Gamma|\leq \delta$ it holds that
         \begin{equation}\llabel{evd2}
         \begin{split}
             &\nnorm{\psi_{\gamma,n}(\Theta^{\xi,\gamma}_{n-1}-\bfv,\bfm^{\xi,\gamma}_{n-1})-\cA(\Theta^{\xi,\gamma}_{n-1}-\bfv,\bfm^{\xi,\gamma}_{n-1})}
            \leq \frac{\varrho-\sqrt{\alpha}}{4}\nnorm{(\Theta^{\xi,\gamma}_{n-1}-\bfv,\bfm^{\xi,\gamma}_{n-1})}.
             \end{split}
         \end{equation}}
\argument{\lref{evd2};\lref{eq3};\lref{eq9};\lref{evd1};}{that there exists $\delta\in (0,1)$ which satisfies for every $n\in\N$, $\xi\in \cd$ and every $\gamma\colon\N\to\R$ with  $\nnorm{(\Theta_{n-1}^{\xi,\gamma}-\bfv,\bfm_{n-1}^{\xi,\gamma})}\leq \delta$ and $\sup_{k\in\N}|\gamma_k-\Gamma|\leq \delta$ it holds that
         \begin{equation}\llabel{eqmain}
         \begin{split}
            &\nnorm{(\Theta_n^{\xi,\gamma}-\bfv,\bfm_n^{\xi,\gamma})}
             \\
             &\leq \nnorm{\cA(\Theta^{\xi,\gamma}_{n-1}-\bfv,\bfm^{\xi,\gamma}_{n-1})}+\nnorm{\psi_{\gamma,n}(\Theta^{\xi,\gamma}_{n-1}-\bfv,\bfm^{\xi,\gamma}_{n-1})-\cA(\Theta^{\xi,\gamma}_{n-1}-\bfv,\bfm^{\xi,\gamma}_{n-1})}\\
             &\leq \nnorm{\cA'(0)(\Theta^{\xi,\gamma}_{n-1}-\bfv,\bfm^{\xi,\gamma}_{n-1})}+\nnorm{\cA(\Theta^{\xi,\gamma}_{n-1}-\bfv,\bfm^{\xi,\gamma}_{n-1})-\cA'(0)(\Theta^{\xi,\gamma}_{n-1}-\bfv,\bfm^{\xi,\gamma}_{n-1})}\\
             &+\nnorm{\psi_{\gamma,n}(\Theta^{\xi,\gamma}_{n-1}-\bfv,\bfm^{\xi,\gamma}_{n-1})-\cA(\Theta^{\xi,\gamma}_{n-1}-\bfv,\bfm^{\xi,\gamma}_{n-1})}\\
             &\leq \frac{\varrho+\sqrt{\alpha}}{2}\nnorm{(\Theta^{\xi,\gamma}_{n-1}-\bfv,\bfm^{\xi,\gamma}_{n-1})}+\frac{\varrho-\sqrt{\alpha}}{4}\nnorm{(\Theta^{\xi,\gamma}_{n-1}-\bfv,\bfm^{\xi,\gamma}_{n-1})}+\frac{\varrho-\sqrt{\alpha}}{4}\nnorm{(\Theta^{\xi,\gamma}_{n-1}-\bfv,\bfm^{\xi,\gamma}_{n-1})}\\
             &=\varrho\nnorm{(\Theta^{\xi,\gamma}_{n-1}-\bfv,\bfm^{\xi,\gamma}_{n-1})}.
             \end{split}  
         \end{equation}}
         \startnewargseq
         In the following we prove that for every $n\in\N_0$, $\xi\in \cd$ and every $\gamma\colon\N\to\R$ with $\sup_{k\in\N}|\gamma_k-\Gamma|< \delta$ and $\|\xi-\bfv\|\leq \delta\mu_1$ it holds that
         \begin{equation}\llabel{need to prove}
             \nnorm{(\Theta_n^{\xi,\gamma}-\bfv,\bfm_n^{\xi,\gamma})}\leq(\mu_1)^{-1}\|\xi-\bfv\| \varrho^n.
         \end{equation}
         We prove \lref{need to prove} by induction on $n\in\N_0$. For the base case $n=0$ note that \lref{def: Theta} and \lref{eq14} show that for every $\xi\in \cd$ and every $\gamma\colon\N\to\R$ with $\|\xi-\bfv\|\leq \delta\mu_1$ it holds that
         \begin{equation}\llabel{ntp: base 1}
             \nnorm{(\Theta_0^{\xi,\gamma}-\bfv,\bfm_0^{\xi,\gamma})}=\nnorm{(\xi-\bfv,0)}\leq (\mu_1)^{-1}\|(\xi-\bfv,0)\|= (\mu_1)^{-1}\|\xi-\bfv\|=(\mu_1)^{-1}\|\xi-\bfv\|\varrho^{0}.
         \end{equation}
         This shows \lref{need to prove} in the base case $n=0$. For the induction step we assume that there exists $n\in\N$ which satisfies for every  $\xi\in \cd$ and every $\gamma\colon\N\to\R$ with $\sup_{k\in\N}|\gamma_k-\Gamma|< \delta$ and $\|\xi-\bfv\|\leq \delta\mu_1$ that
         \begin{equation}\llabel{ntp: assume}
              \nnorm{(\Theta_{n-1}^{\xi,\gamma}-\bfv,\bfm_{n-1}^{\xi,\gamma})}\leq(\mu_1)^{-1}\|\xi-\bfv\|\varrho^{n-1}.
         \end{equation}
         \startnewargseq
         \argument{\lref{eq14};\lref{ntp: assume};the fact that $0<\varrho<1$}{that for every $\xi\in \cd$ and every $\gamma\colon\N\to\R$ with $\sup_{k\in\N}|\gamma_k-\Gamma|< \delta$ and $\|\xi-\bfv\|\leq \delta\mu_1$ it holds that
         \begin{equation}\llabel{ntp: eq1}
         \begin{split}
           \textstyle  \nnorm{(\Theta_{n-1}^{\xi,\gamma}-\bfv,\bfm_{n-1}^{\xi,\gamma})}\leq (\mu_1)^{-1}\varrho^{n-1}\|\xi-\bfv\|\leq (\mu_1)^{-1}\|\xi-\bfv\|\leq \delta.
             \end{split}
         \end{equation}}
         \argument{\lref{ntp: eq1};\lref{eqmain};\lref{ntp: assume}}{for every $\xi\in \cd$ and every $\gamma\colon\N\to\R$ with $\sup_{k\in\N}|\gamma_k-\Gamma|< \delta$  and $\|\xi-\bfv\|\leq \delta\mu_1$ that
         \begin{equation}\llabel{ntp: eq2}
             \begin{split}
             &\nnorm{(\Theta_n^{\xi,\gamma}-\bfv,\bfm_n^{\xi,\gamma})}\leq  \varrho\nnorm{(\Theta_{n-1}^{\xi,\gamma}-\bfv,\bfm_{n-1}^{\xi,\gamma})}\leq \varrho(\mu_1)^{-1} \|\xi-\bfv\|\varrho^{n-1}=(\mu_1)^{-1}\|\xi-\bfv\|\varrho^{n}.
             \end{split}
         \end{equation}}
         \argument{\lref{ntp: eq2};\lref{ntp: assume};induction}{\lref{need to prove}\dott}
         \startnewargseq
         \argument{\lref{eq14};\lref{need to prove};}{that for every $n\in\N_0$, $\xi\in \cd$ and every $\gamma\colon\N\to\R$ with $\sup_{k\in\N}|\gamma_k-\Gamma|<\delta$ and $\|\xi-\bfv\|\leq \delta\mu_1$ it holds that
         \begin{equation}\llabel{eq19}
             \|\Theta_{n}^{\xi,\gamma}-\bfv\|\leq \|(\Theta_n^{\xi,\gamma}-\bfv,\bfm_n^{\xi,\gamma})\|\leq \mu_2\nnorm{(\Theta_n^{\xi,\gamma}-\bfv,\bfm_n^{\xi,\gamma})}\leq \mu_2(\mu_1)^{-1}\|\xi-\bfv\|\varrho^{n}.
         \end{equation}}
         \argument{\lref{eq19};\lref{def: bfm};\lref{def: Theta};}{that there exist $\scrC,\eta \in (0,\infty)$ such that for every $\gamma\colon\N\to(\Gamma-\eta,\Gamma+\eta)$, every $\bfm=(\bfm^1,\dots,\bfm^\nb)\colon \mathbb{N}_0\to\cd$, and every $\Theta=(\Theta^1,\dots,\Theta^\nb)\colon \mathbb{N}_0\to\cd$ with the property that for all $n \in \mathbb{N}$, $i\in\{1,2,\dots,\nb\}$ it holds that
         \begin{equation}
         \bfm_0=0,\qquad \bfm_n= \alpha \bfm_{n-1}+(1-\alpha)(\nabla\cL)(\Theta_{n-1}),
            \end{equation}
            \begin{equation}
           \|\Theta_0-\bfv\|< \eta,\qquad\text{and}\qquad  \Theta_n^i=\Theta_{n-1}^i-\gamma_n
           \Phi_i(\bfm_n^i)
         \end{equation}
          we have for all $n\in\N_0$ that
$    \textstyle
    \|\Theta_n-\bfv\|\leq \scrC \|\Theta_0-\bfv\|\varrho^n$\dott}
\end{aproof}
\begin{samepage}
\begin{tcolorbox}[colback=white!95!gray,
                  colframe=black,
                  boxrule=0.5pt,
                  sharp corners,
                  enhanced,
                  breakable,
                 ]
\begin{athm}{cor}{theo: abstract optimal rate}
      Let $\nb\in \mathbb{N}$, $\bfb\in (0,\infty)$,  $\chi\in (1,\infty)$, $\alpha=(\frac{\sqrt{ \chi } - 1 }  { \sqrt{ \chi } + 1 })^2$, for every $i\in\{1,2,\dots,\nb\}$ let $d_i,\fd_i\in\N$, $\Phi_i\in C^1(\R^{d_i\times\fd_i},\R^{d_i\times\fd_i})$ satisfy for all $\theta\in \R^{d_i\times\fd_i}$ that $(\Phi_i)'(0)\theta-\bfb\theta=0=\Phi_i(0)$, and let $\cd=\times_{i=1}^{\nb}\R^{d_i\times\fd_i}$. Then for every $\convexcoe \in (0,\infty)$, $\vartheta \in \cd$, $\smalll\in \allowbreak C^{2}(\cd,\mathbb{R})$, $\varrho\in (\frac{\sqrt{ \chi } - 1 }  { \sqrt{ \chi } + 1 },1)$ 
there exist $\fC,\Gamma, \delta \in (0,\infty)$ such that for every $\gamma\colon\N\to(\Gamma-\delta,\Gamma+\delta)$, every $\bfm=(\bfm^1,\dots,\bfm^\nb)\colon \mathbb{N}_0\to\cd$, and every $\Theta=(\Theta^1,\dots,\Theta^\nb)\colon \mathbb{N}_0\to\cd$ with the property that for all $n \in \mathbb{N}$, $i\in\{1,2,\dots,\nb\}$, $\theta\in\cd$ it holds that
 \begin{equation}
     \textcolor{magenta}{(\nabla\smalll)(\vartheta)=0},\qquad\textcolor{magenta}{\convexcoe\|\theta\|^2\leq \smalll'' (\vartheta)(\theta,\theta)\leq \convexcoe\chi \| \theta \|^2},
 \end{equation} 
         \begin{equation}
            \textcolor{magenta}{\bfm_0=0},\qquad \textcolor{magenta}{\bfm_n= \alpha \bfm_{n-1}+(1-\alpha)(\nabla\smalll)(\Theta_{n-1})},
            \end{equation}
            \begin{equation}
            \textcolor{magenta}{\|\Theta_0-\vartheta\|< \delta},\qquad\text{and}\qquad  \textcolor{magenta}{\Theta_n^i=\Theta_{n-1}^i-\gamma_n
           \Phi_i(\bfm_n^i})
         \end{equation}
          we have for all $n\in\N_0$ that
$    \textcolor{magenta}{ \textstyle
     \|\Theta_n-\vartheta\|\leq \fC\|\Theta_0-\vartheta\|\varrho^{n}}$.
\end{athm}
\end{tcolorbox}
\end{samepage}
\begin{aproof}
    \argument{\cref{theo: muon higher bound};}{that for every $\convexcoe \in (0,\infty)$, $\vartheta \in \cd$, $\smalll\in \allowbreak C^{2}(\cd,\mathbb{R})$, $\varrho\in (\frac{\sqrt{ \chi } - 1 }  { \sqrt{ \chi } + 1 },1)$ 
there exist $\fC,\Gamma, \delta \in (0,\infty)$ such that for every $\gamma\colon\N\to(\Gamma-\delta,\Gamma+\delta)$, every $\bfm=(\bfm^1,\dots,\bfm^\nb)\colon \mathbb{N}_0\to\cd$, and every $\Theta=(\Theta^1,\dots,\Theta^\nb)\colon \mathbb{N}_0\to\cd$ with the property that for all $n \in \mathbb{N}$, $i\in\{1,2,\dots,\nb\}$, $\theta\in\cd$ it holds that
 \begin{equation}
     (\nabla\smalll)(\vartheta)=0,\qquad\convexcoe\|\theta\|^2\leq \smalll'' (\vartheta)(\theta,\theta)\leq \convexcoe\chi \| \theta \|^2,
 \end{equation} 
         \begin{equation}
            \bfm_0=0,\qquad \bfm_n= \alpha \bfm_{n-1}+(1-\alpha)(\nabla\smalll)(\Theta_{n-1}),
            \end{equation}
            \begin{equation}
            \|\Theta_0-\vartheta\|< \delta,\qquad\text{and}\qquad  \Theta_n^i=\Theta_{n-1}^i-\gamma_n
          \Phi_i(\bfm_n^i)
         \end{equation}
          we have for all $n\in\N_0$ that
$     \textstyle
     \|\Theta_n-\vartheta\|\leq \fC\|\Theta_0-\vartheta\|\varrho^{n}$\dott}
\end{aproof}
In the next elementary result, \cref{lem: NS represent}, we show that the generalized \NS\ functions are continuously differentiable, vanish at zero, and have the derivative $\theta\mapsto(\nscoe_{0,1})^{-1}\bigl[\prod_{k=1}^K\nscoe_{k,0}\bigr]\theta$ at zero.

\begin{tcolorbox}[colback=white!95!gray,
                  colframe=black,
                  boxrule=0.5pt,
                  sharp corners,
                  enhanced,
                  breakable,
                 ]
                 \begin{athm}{lemma}{lem: NS represent}
                     Let $d,\fd\in\N$, $K\in\N_0$ and let $\nscoe=(\nscoe_{i,j})_{(i,j)\in(\N_0)^2}\colon (\N_0)^2\to\R$ satisfy $\min\{\nscoe_{0,0},\nscoe_{0,1}\}>0$ and $\#(\nscoe^{ - 1 }( \R\backslash\{0\} )) <\infty$. Then 
                     \begin{enumerate}[label=(\roman*)]
                         \item \label{item 1: NS represent} it holds that $ \NewtonSchulzAlgorithm{\nscoe}{K}\in C^1(\R^{d\times\fd},\R^{d\times\fd})$,
                         \item \label{item 2: NS represent} it holds that $\NewtonSchulzAlgorithm{\nscoe}{K}(0)=0$, and 
                         \item \label{item 3: NS represent} it holds for all $\theta\in \R^{d\times\fd}$ that $(\NewtonSchulzAlgorithm{\nscoe}{K})'(0)\theta=(\nscoe_{0,1})^{-1}\bigl[\prod_{k=1}^K\nscoe_{k,0}\bigr] \theta$.
                     \end{enumerate}
                 \end{athm}
                 \end{tcolorbox}
    \begin{aproof}
     Throughout this proof let $\psi\colon \R^{d\times\fd}\to\R^{d\times\fd}$ satisfy for all $\theta\in \R^{d\times\fd}$ that
                 \begin{equation}\llabel{def: T}
\psi(\theta)=\theta(\nscoe_{0,0}\|\theta\|+\nscoe_{0,1})^{-1}.
                 \end{equation}
   \argument{\cref{T_B_D};}{that for all $\theta\in\R^{d\times\fd}$ it holds that
   \begin{equation}\llabel{eqq1}
       \TNewtonSchulzAlgorithm{\nscoe}{0}{\theta}=(\nscoe_{0,0}\|\theta\| + \nscoe_{0,1})^{-1} \theta.
   \end{equation}}
   \argument{\lref{eqq1};the fact that $\nscoe_{0,0}>0$; the fact that $\nscoe_{0,1}>0$}{that \llabel{eqq2}
       $\NewtonSchulzAlgorithm{\nscoe}{0}\in C^1(\R^{d\times\fd},\R^{d\times\fd})$\dott}
   \argument{\lref{eqq2};\cref{T_B_D};induction}{that for all $n\in\N_0$ it holds that
   \begin{equation}\llabel{eqq2.1}
        \NewtonSchulzAlgorithm{\nscoe}{n}\in C^1(\R^{d\times\fd},\R^{d\times\fd}).
   \end{equation}}
   \argument{\lref{eqq2.1};}{\cref{item 1: NS represent}\dott}
   \startnewargseq
   \argument{\lref{eqq1};}{that
   \llabel{eqq3} $\TNewtonSchulzAlgorithm{\nscoe}{0}{0}=0$\dott}
   \argument{\lref{eqq3};\cref{T_B_D};induction}{that for all $n\in\N_0$ it holds that
   \begin{equation}\llabel{eqq4}
       \TNewtonSchulzAlgorithm{\nscoe}{n}{0}=0.
   \end{equation}}
   \argument{\lref{eqq4};}{\cref{item 2: NS represent}\dott}
   \startnewargseq
   In the following we prove that for all $k\in \N_0$ there exists a polynomial $P\colon \R^{d\times\fd}\to\R^{\fd\times\fd}$ such that for all $\theta\in\R^{d\times\fd}$ it holds that
\begin{equation}\llabel{need to prove}
   \textstyle \TNewtonSchulzAlgorithm{\nscoe}{k}{\theta}=\bigl[\prod_{n=1}^k\nscoe_{n,0}\bigr]\psi(\theta)+\psi(\theta) P(\psi(\theta))(\psi(\theta))^\top\psi(\theta).
\end{equation}
We prove \lref{need to prove} by induction on $k\in\N_0$. For the base case $k=0$ note that \cref{T_B_D} and \lref{def: T} show that for all $\theta\in \R^{d\times\fd}$ it holds that
\begin{equation}
\begin{split}
    \textstyle \TNewtonSchulzAlgorithm{\nscoe}{0}{\theta}=\theta(\nscoe_{0,0}\|\theta\|+\nscoe_{0,1})^{-1}\textstyle=\bigl[\prod_{n=1}^0\nscoe_{n,0}\bigr]\psi(\theta)+\psi(\theta)0(\psi(\theta))^\top\psi(\theta) .
    \end{split}
\end{equation}
This proves \lref{need to prove} in the case $k=0$. For the induction step we assume that there exist $k\in\N$ and a polynomial $P\colon \R^{d\times\fd}\to\R^{\fd\times\fd}$ which satisfy for all $\theta\in\R^{d\times\fd}$ that
\begin{equation}\llabel{assume}
   \textstyle \TNewtonSchulzAlgorithm{\nscoe}{k-1}{\theta}=\bigl[\prod_{n=1}^{k-1}\nscoe_{n,0}\bigr]\psi(\theta)+\psi(\theta)P(\psi(\theta))(\psi(\theta))^\top\psi(\theta) .
\end{equation}
\argument{\lref{assume};}{that there exist polynomials $R\colon \R^{d\times\fd}\to\R^{\fd\times\fd}$ and $S\colon \R^{d\times\fd}\to\R^{d\times\fd}$ which satisfy for all $\theta\in \R^{d\times\fd}$ that
\begin{equation}\llabel{eq2.1}
    \TNewtonSchulzAlgorithm{\nscoe}{k-1}{\theta}=\psi(\theta)R(\psi(\theta))=S(\psi(\theta)).
\end{equation}}
\startnewargseq
\argument{\cref{T_B_D};\lref{assume};\lref{eq2.1}}{that for all $\theta\in \R^{d\times\fd}$ it holds that
\begin{equation}\llabel{eq2}
\begin{split}
\TNewtonSchulzAlgorithm{\nscoe}{k}{\theta}
&=\textstyle\sum_{i=0}^\infty\nscoe_{k,i} (\TNewtonSchulzAlgorithm{\nscoe}{k-1}{\theta}(\TNewtonSchulzAlgorithm{\nscoe}{k-1}{\theta})^\top)^{i}\TNewtonSchulzAlgorithm{\nscoe}{k-1}{\theta}\\
&=\textstyle\sum_{i=0}^\infty\nscoe_{k,i} \big(\TNewtonSchulzAlgorithm{\nscoe}{k-1}{\theta}(\TNewtonSchulzAlgorithm{\nscoe}{k-1}{\theta})^\top\bigr)^i\bigl(\bigl[\prod_{n=1}^{k-1}\nscoe_{n,0}\bigr]\psi(\theta)+\psi(\theta)P(\psi(\theta))(\psi(\theta))^\top\psi(\theta) \bigr)\\
&=\textstyle\nscoe_{k,0}\bigl(\bigl[\prod_{n=1}^{k-1}\nscoe_{n,0}\bigr]\psi(\theta)+ \psi(\theta)P(\psi(\theta))(\psi(\theta))^\top\psi(\theta)\bigr)\\
&+\textstyle\sum_{i=1}^\infty\nscoe_{k,i} \big(\TNewtonSchulzAlgorithm{\nscoe}{k-1}{\theta}(\TNewtonSchulzAlgorithm{\nscoe}{k-1}{\theta})^\top\bigr)^{i}\bigl[\prod_{n=1}^{k-1}\nscoe_{n,0}\bigr]\psi(\theta)\\
&+\textstyle\sum_{i=1}^\infty\nscoe_{k,i} \big(\TNewtonSchulzAlgorithm{\nscoe}{k-1}{\theta}(\TNewtonSchulzAlgorithm{\nscoe}{k-1}{\theta})^\top\bigr)^{i}\psi(\theta)P(\psi(\theta))(\psi(\theta))^\top\psi(\theta)\\
&=\textstyle\bigl[\prod_{n=1}^{k}\nscoe_{n,0}\bigr]\psi(\theta)+\nscoe_{k,0}\psi(\theta)P(\psi(\theta))(\psi(\theta))^\top\psi(\theta)\\
&+\textstyle\sum_{i=1}^\infty\nscoe_{k,i} \bigl[\prod_{n=1}^{k-1}\nscoe_{n,0}\bigr]\psi(\theta)R(\psi(\theta))\big((S(\psi(\theta)))^\top S(\psi(\theta))\bigr)^{i-1}(R(\psi(\theta)))^\top(\psi(\theta))^\top\psi(\theta)\\
&+\textstyle\sum_{i=1}^\infty\nscoe_{k,i} \psi(\theta)R(\psi(\theta))(S(\psi(\theta)))^\top\big(S(\psi(\theta))(S(\psi(\theta)))^\top\bigr)^{i-1} \psi(\theta)P(\psi(\theta))(\psi(\theta))^\top\psi(\theta).
\end{split}
\end{equation}}
\argument{the assumption that $\#(\nscoe^{ - 1 }( \R\backslash\{0\} )) <\infty$}{that \llabel{argg1} $\#\{i\in\N\colon \nscoe_{k,i}\neq 0\}<\infty$\dott}
\argument{\lref{argg1};the fact that $P$ is a polynomial;the fact that $R$ is a polynomial; the fact that $S$ is a polynomial}{that
\begin{equation}
    \textstyle\R^{d\times\fd}\ni x\mapsto \nscoe_{k,0}P(x)\in \R^{\fd\times\fd},
\end{equation}
\begin{equation}
     \textstyle\R^{d\times\fd}\ni x\mapsto \sum_{i=1}^\infty\nscoe_{k,i} \bigl[\prod_{n=1}^{k-1}\nscoe_{n,0}\bigr]R(x)\big((S(x))^\top S(x)\bigr)^{i-1}(R(x))^\top\in \R^{\fd\times\fd},
\end{equation}
and
\begin{equation}
    \textstyle\R^{d\times\fd}\ni x\mapsto \textstyle\sum_{i=1}^\infty\nscoe_{k,i} R(x)(S(x))^\top\big(S(x)(S(x))^\top\bigr)^{i-1} xP(x) \in \R^{\fd\times\fd}
\end{equation}
are \llabel{arg1} polynomials\dott}
\argument{\lref{arg1};\lref{eq2}}{that there exists a polynomial $Q\colon \R^{d\times\fd}\to\R^{\fd\times\fd}$ such that for all $\theta\in \R^{d\times\fd}$ it holds that
\begin{equation}\llabel{eq3}
    \TNewtonSchulzAlgorithm{\nscoe}{k}{\theta}=\textstyle\big[\prod_{n=1}^{k}\nscoe_{n,0}\big]\psi(\theta)+\psi(\theta)Q(\psi(\theta))(\psi(\theta))^\top\psi(\theta) .
\end{equation}}
\argument{\lref{eq3}; \lref{assume}; induction}{\lref{need to prove}\dott}
\startnewargseq
\argument{\lref{def: T};}{that for all $\theta\in \R^{d\times\fd}$ it holds that
\begin{equation}\llabel{eq9}
    \psi'(0)\theta=(\nscoe_{0,1})^{-1}\theta.
\end{equation}}
\argument{\lref{def: T};}{that
\begin{equation}\llabel{eq10}
    \psi(0)=0.
\end{equation}}
\argument{\lref{need to prove};\lref{eq9};\lref{eq10}; the chain rule;}{\cref{item 3: NS represent}\dott}
    \end{aproof}
\begin{tcolorbox}[colback=white!95!gray,
                  colframe=black,
                  boxrule=0.5pt,
                  sharp corners,
                  enhanced,
                  breakable,
                 ]
\begin{athm}{cor}{main theorem higher bound}
Let $\nb\in \mathbb{N}$, $K\in\N_0$, $d_1,d_2,\dots,d_\nb,\fd_1,\fd_2,\dots,\fd_\nb\in\N$, $\cd=\times_{i=1}^{\nb}\R^{d_i\times\fd_i}$,   $\chi\in (1,\infty)$, $\alpha=(\frac{\sqrt{ \chi } - 1 }  { \sqrt{ \chi } + 1 })^2$, let $
\nscoe=(\nscoe_{i,j})_{(i,j)\in (\mathbb{N}_0)^2}\colon(\mathbb{N}_0)^2\allowbreak\to \mathbb{R}$ satisfy $\textcolor{magenta}{\min\{\nscoe_{0,0},\nscoe_{0,1}\}>0}$ and $\textcolor{magenta}{\#(\nscoe^{ - 1 }( \mathbb{R}\backslash\{0\} )) <\infty}$, and assume for all $i\in \mathbb{N}\cap[0,K]$ that $\textcolor{magenta}{\nscoe_{i,0}>0}$. Then for every $\convexcoe \in (0,\infty)$, $\vartheta \in \cd$, $\smalll\in \allowbreak C^{2}(\cd,\mathbb{R})$, $\varrho\in (\frac{\sqrt{ \chi } - 1 }  { \sqrt{ \chi } + 1 },1)$ 
there exist $\fC,\Gamma, \delta \in (0,\infty)$ such that for every $\gamma\colon\N\to(\Gamma-\delta,\Gamma+\delta)$, every $\bfm=(\bfm^1,\dots,\bfm^\nb)\colon \mathbb{N}_0\to\cd$, and every $\Theta=(\Theta^1,\dots,\Theta^\nb)\colon \mathbb{N}_0\to\cd$ with the property that for all $n \in \mathbb{N}$, $i\in\{1,2,\dots,\nb\}$, $\theta\in \cd$ it holds that
 \begin{equation}
     \textcolor{magenta}{(\nabla\smalll)(\vartheta)=0},\qquad\textcolor{magenta}{\convexcoe\|\theta\|^2\leq \smalll'' (\vartheta)(\theta,\theta)\leq \convexcoe\chi \| \theta \|^2},
 \end{equation} 
         \begin{equation}
            \textcolor{magenta}{\bfm_0=0},\qquad \textcolor{magenta}{\bfm_n= \alpha \bfm_{n-1}+(1-\alpha)(\nabla\smalll)(\Theta_{n-1})},
            \end{equation}
            \begin{equation}
            \textcolor{magenta}{\|\Theta_0-\vartheta\|< \delta},\qquad\text{and}\qquad  \textcolor{magenta}{\Theta_n^i=\Theta_{n-1}^i-\gamma_n
           \TNewtonSchulzAlgorithm{\nscoe}{K}{\bfm_n^i}}
         \end{equation}
          we have for all $n\in\N_0$ that
$    \textcolor{magenta}{ \|\Theta_n-\vartheta\|\leq \fC\|\Theta_0-\vartheta\|\varrho^n}$
(cf.\ Definition \ref{definition: NS}).
\end{athm}
\end{tcolorbox}
\begin{aproof}
    \argument{\cref{theo: abstract optimal rate};\cref{lem: NS represent};the fact that $(\nscoe_{0,1})^{-1}\prod_{k=1}^K\nscoe_{k,0}>0$}{that for every $\convexcoe \in (0,\infty)$, $\vartheta \in\cd$, $\smalll\in \allowbreak C^{2}(\cd,\mathbb{R})$, $\varrho\in (\frac{\sqrt{ \chi } - 1 }  { \sqrt{ \chi } + 1 },1)$ 
there exist $\fC,\Gamma, \delta \in (0,\infty)$ such that for every $\gamma\colon\N\to(\Gamma-\delta,\Gamma+\delta)$, every $\bfm=(\bfm^1,\dots,\bfm^\nb)\colon \mathbb{N}_0\to\cd$, and every $\Theta=(\Theta^1,\dots,\Theta^\nb)\colon \mathbb{N}_0\to\cd$ with the property that for all $n \in \mathbb{N}$, $i\in\{1,2,\dots,\nb\}$, $\theta\in\cd$ it holds that
 \begin{equation}
     (\nabla\smalll)(\vartheta)=0,\qquad\convexcoe\|\theta\|^2\leq \smalll'' (\vartheta)(\theta,\theta)\leq \convexcoe\chi \| \theta \|^2,
 \end{equation} 
         \begin{equation}
           \bfm_0=0,\qquad \bfm_n= \alpha \bfm_{n-1}+(1-\alpha)(\nabla\smalll)(\Theta_{n-1}),
            \end{equation}
            \begin{equation}
            \|\Theta_0-\vartheta\|< \delta,\qquad\text{and}\qquad  \Theta_n^i=\Theta_{n-1}^i-\gamma_n
           \TNewtonSchulzAlgorithm{\nscoe}{K}{\bfm_n^i}
         \end{equation}
          we have for all $n\in\N_0$ that
$    \|\Theta_n-\vartheta\|\leq \fC\|\Theta_0-\vartheta\|\varrho^n$\dott}
\end{aproof}
\subsection{Lower bounds on convergence rates of momentum and MUON}
\begin{tcolorbox}[colback=white!95!gray,
                  colframe=black,
                  boxrule=0.5pt,
                  sharp corners,
                  enhanced,
                  breakable,
                 ]
\begin{athm}{lemma}{lem: lower bound}
      Let $\nb\in\N$, $\bfb,\convexcoe,\delta,\gamma,\fC\in (0,\infty)$, $\alpha\in (0,1)$, $\varrho\in (0,\sqrt{\alpha} )$, for every $i\in\{1,2,\dots,\nb\}$ let $d_i,\fd_i\in\N$, $\Phi_i\in C^1(\R^{d_i\times\fd_i},\R^{d_i\times\fd_i})$ satisfy for all $\theta\in \R^{d_i\times\fd_i}$ that $(\Phi_i)'(0)\theta-\bfb\theta=0=\Phi_i(0)$, let $\cd=\times_{i=1}^{\nb}\R^{d_i\times\fd_i}$, and let $\smalll\colon \cd\to\R$ satisfy for all $\theta\in \cd$ that $\smalll(\theta)=\convexcoe\|\theta\|^2$. Then there exist $\bfm=(\bfm^1,\dots,\bfm^\nb)\colon \mathbb{N}_0\to\cd$ and $\Theta=(\Theta^1,\dots,\Theta^\nb)\colon \mathbb{N}_0\to\cd$  such that
\begin{enumerate}[label=(\roman*)]
    \item \llabel{item 1} it holds for all $n \in \mathbb{N}$, $i\in\{1,2,\dots,\nb\}$ that
         \begin{equation}
            \textcolor{magenta}{\bfm_0=0},\qquad \textcolor{magenta}{\bfm_n= \alpha \bfm_{n-1}+(1-\alpha)(\nabla\smalll)(\Theta_{n-1})},
            \end{equation}
            \begin{equation}
            \textcolor{magenta}{\|\Theta_0\|< \delta},\qquad\text{and}\qquad  \textcolor{magenta}{\Theta_n^i=\Theta_{n-1}^i-\gamma
           \Phi_i(\bfm_n^i)}
         \end{equation}
         and
         \item \llabel{item 2} there exists $n\in\N$ such that $\textcolor{magenta}{ \textstyle\|\Theta_n\|>\fC\|\Theta_0\|\varrho^n}$.
\end{enumerate} 
\end{athm}
\end{tcolorbox}
\begin{aproof}
    Throughout this proof let $\lambda\in\R$ satisfy $\lambda=2r(1-\alpha)$, let $\bfe=(\bfe^k)_{k\in\{1,2,\dots,\nb\}}=((\bfe_{i,j}^k)_{(i,j)\in\{1,2,\dots,d_k\}\times\{1,2,\dots,\fd_k\}})_{k\in\{1,2,\dots,\nb\}}\allowbreak\in \cd$ satisfy for all $k\in\{1,2,\dots,\nb\}$, $i\in\{1,2,\dots,d_k\}$, $j\in\{1,2,\dots,\fd_k\}$ that
    \begin{equation}\llabel{def: bfe}
        \bfe_{i,j}^k=\mathbbm 1_{(1,1,1)}(i,j,k),
    \end{equation}
     for every $t\in\R$ let $\bfm^{t}=(\bfm^{t,1},\dots,\bfm^{t,\nb})\colon\N_0\to\cd$ and $\Theta^{t}=(\Theta^{t,1},\dots,\Theta^{t,\nb})\colon\N_0\to\cd$ satisfy for all $n\in\N$, $i\in\{1,2,\dots,\nb\}$ that
     \begin{equation}\llabel{def: bfm}
      \bfm^{t}_0=0,\qquad \bfm^{t}_n= \alpha \bfm^{t}_{n-1}+\lambda\Theta^{t}_{n-1},
            \end{equation}
            \begin{equation}\llabel{def: Theta}
        \Theta^{t}_0=t\bfe,\qquad\text{and}\qquad \Theta^{t,i}_n=\Theta^{t,i}_{n-1}-\gamma
           \Phi_i(\bfm^{t,i}_n),
         \end{equation}
         and let $\bfu\colon \N_0\to\R$, $\bfv\colon\N_0\to\R$ satisfy for all $n\in\N$ that
         \begin{equation}\llabel{def: bfu}
             \bfu_0=1,\qquad \bfu_n=\bfu_{n-1}-\gamma\bfb \bfv_{n},\qquad \bfv_0=0,\qqandqq\bfv_n=\alpha\bfv_{n-1}+\lambda \bfu_{n-1}.
         \end{equation}
         \argument{\lref{def: bfu}}{that for all $n\in \N$ it holds that
         \begin{equation}\llabel{eqq1}
         \begin{split}
            & \bfu_n=\bfu_{n-1}-\gamma\bfb(\alpha \bfv_{n-1}+\lambda\bfu_{n-1})=\bfu_{n-1}(1-\gamma\lambda\bfb)-\gamma\alpha\bfb \bfv_{n-1}\qandq\bfv_n=\alpha \bfv_{n-1}+\lambda\bfu_{n-1}.
             \end{split}
         \end{equation}}
          \argument{\lref{eqq1};\lref{def: bfu}}{for all $n\in\N\backslash\{1\}$ that
         \begin{equation}\llabel{eqq2}
         \begin{split}
\bfu_n=\bfu_{n-1}(1-\gamma\lambda\bfb)-\alpha(\bfu_{n-2}-\bfu_{n-1})=\bfu_{n-1}(1+\alpha-\gamma\lambda\bfb)-\alpha\bfu_{n-2}.
             \end{split}
         \end{equation}}
    \argument{\lref{eqq2};}{for all $n\in\N\backslash\{1\}$ that
    \begin{equation}\llabel{eqq3}
    \begin{split}
        \bfu_{n+1}\bfu_{n-1}-(\bfu_{n})^2&=[\bfu_{n}(1+\alpha-\gamma\lambda\bfb)-\alpha\bfu_{n-1}]\bfu_{n-1}-[\bfu_{n-1}(1+\alpha-\gamma\lambda\bfb)-\alpha\bfu_{n-2}]\bfu_n\\
        &=-\alpha(\bfu_{n-1})^2+\alpha\bfu_{n-2}\bfu_n=\alpha(\bfu_{n}\bfu_{n-2}-(\bfu_{n-1})^2).
        \end{split}
    \end{equation}}
    \argument{\lref{def: bfu};\lref{eqq2}}{that
    \begin{equation}\llabel{eqq4}
        \bfv_1=\lambda,\qquad \bfu_1=1-\gamma\bfb\lambda,\qqandqq\bfu_2=(1+\alpha-\gamma\bfb\lambda)(1-\gamma\bfb\lambda)-\alpha.
    \end{equation}}
    \argument{\lref{eqq4};\lref{def: bfu}}{that
    \begin{equation}\llabel{eqq5}
    \begin{split}
        \bfu_2\bfu_0-(\bfu_1)^2&=(1+\alpha-\gamma\bfb\lambda)(1-\gamma\bfb\lambda)-\alpha-(1-\gamma\bfb\lambda)^2=-\alpha\gamma\bfb\lambda.
        \end{split}
    \end{equation}}
    \argument{\lref{eqq5};\lref{eqq3};induction}{that for all $n\in\N$ it holds that
    \begin{equation}\llabel{eqq6}
        \bfu_{n+1}\bfu_{n-1}-(\bfu_n)^2=-\gamma\bfb\lambda\alpha^n.
    \end{equation}}
    \argument{\lref{eqq6};the fact that $\sqrt{\alpha}>\varrho$}{that there exists $\ell\in\N$ which satisfies that
    \begin{equation}\llabel{def: N}
        |\bfu_\ell|\geq 2\fC\varrho^{\ell}.
    \end{equation}}
    \startnewargseq
       In the following we prove that for all $n\in\{0,1,\dots,\ell\}$ it holds that
       \begin{equation}\llabel{need to prove}
           \limsup_{t\to 0}\biggl\|\frac{\Theta_n^{t}-t \bfu_n\bfe}{t}\biggr\|=0\qqandqq  \limsup_{t\to 0}\biggl\|\frac{\bfm_n^{t}-t \bfv_n\bfe}{t}\biggr\|=0.
       \end{equation}
       We prove \lref{need to prove} by induction on $n\in\{0,1,\dots,\ell\}$. For the base case $n=0$ note that \lref{def: bfm}, \lref{def: Theta}, and \lref{def: bfu} show that for all $t\in\R\backslash\{0\}$ it holds that
       \begin{equation}\llabel{base: eq1}
           \frac{\Theta_0^{t}-t \bfu_0\bfe}{t}=\frac{t\bfe-1t\bfe}{t}=0\qqandqq \frac{\bfm_0^{t}-t \bfv_0\bfe}{t}=\frac{0-0t\bfe}{t}=0.
       \end{equation}
       This proves \lref{need to prove} in the base case $n=0$. For the induction step we assume that there exists $n\in\N\cap[0,\ell]$ which satisfies that
       \begin{equation}\llabel{ntp: assume}
             \limsup_{t\to 0}\biggl\|\frac{\Theta_{n-1}^{t}-t \bfu_{n-1}\bfe}{t}\biggr\|=0\qqandqq  \limsup_{t\to 0}\biggl\|\frac{\bfm_{n-1}^{t}-t \bfv_{n-1}\bfe}{t}\biggr\|=0.
       \end{equation}
       \argument{\lref{def: bfm};\lref{def: bfu};\lref{ntp: assume};}{that
       \begin{equation}\llabel{ntp: eq1}
       \begin{split}
           &\limsup_{t\to 0}\biggl\|\frac{\bfm_{n}^{t}-t \bfv_{n}\bfe}{t}\biggr\|= \limsup_{t\to 0}\biggl\|\frac{\alpha \bfm^{t}_{n-1}+\lambda\Theta^{t}_{n-1}-t (\alpha\bfv_{n-1}+\lambda \bfu_{n-1})\bfe}{t}\biggr\|\\
           &\leq \limsup_{t\to 0}\biggl\|\frac{\alpha (\bfm^{t}_{n-1}-t \bfv_{n-1}\bfe)}{t}\biggr\|+\limsup_{t\to 0}\biggl\|\frac{\lambda(\Theta^{t}_{n-1}-t\bfu_{n-1}\bfe)}{t}\biggr\|=0.
           \end{split}
       \end{equation}}
       \argument{\lref{ntp: eq1};}{that
       \begin{equation}\llabel{ntp: eq2}
           \limsup_{t\to0}\biggl\|\frac{\bfm_n^{t}}{t}\biggr\|=\|\bfv_n\bfe\|<\infty.
       \end{equation}}
       \argument{\lref{ntp: eq2};}{that
       \begin{equation}\llabel{ntp: eq3}
         \limsup_{t\to0}\|\bfm_n^{t}\|=0.
       \end{equation}}
       \argument{the fact that for all $i\in\{1,2,\dots,\nb\}$ it holds that $\Phi_i\in C^1(\R^{d_i\times\fd_i},\R^{d_i\times\fd_i})$;the fact that for all $i\in\{1,2,\dots,\nb\}$, $\theta\in\R^{d_i\times\fd_i}$ it holds that $(\Phi_i)'(0)\theta=\bfb\theta$;the fact that for all $i\in\{1,2,\dots,\nb\}$ it holds that $\Phi_i(0)=0$}{that for all $i\in\{1,2,\dots,\nb\}$ it holds that 
       \begin{equation}\llabel{ntp: eq5}
       \begin{split}
&\limsup_{h\searrow 0}\sup_{\theta\in\R^{d_i\times\fd_i}\backslash\{0\},\,\|\theta\|\leq h}\frac{\|\Phi_i(\theta)-\bfb\theta\|}{\|\theta\|}\\
&=  \limsup_{h\searrow 0}\sup_{\theta\in\R^{d_i\times\fd_i}\backslash\{0\},\,\|\theta\|\leq h}\frac{\|\Phi_i(\theta)-\Phi_i(0)-(\Phi_i)'(0)\theta\|}{\|\theta\|}=0.
\end{split}
       \end{equation}}
       \argument{\lref{ntp: eq5};\lref{ntp: eq3}}{that for all $i\in\{1,2,\dots,\nb\}$ it holds that 
       \begin{equation}\llabel{ntp: eq5'}
           \limsup_{t\to 0}\frac{\|\Phi_i(\bfm_n^{t,i})-\bfb\bfm_n^{t,i}\|}{\|\bfm_n^{t,i}\|+t\mathbbm 1_{\{0\}}(\bfm_n^{t,i})}=0.
       \end{equation}}
       \argument{\lref{ntp: eq5'};\lref{ntp: eq2};}{that for all $i\in\{1,2,\dots,\nb\}$ it holds that 
       \begin{equation}\llabel{ntp: eq6}
       \begin{split}
         & \limsup_{t\to 0}\frac{\|\Phi_i(\bfm_n^{t,i})-\bfb\bfm_n^{t,i}\|}{|t|} \\
          &\leq\biggl( \limsup_{t\to 0}\frac{\|\Phi_i(\bfm_n^{t,i})-\bfb\bfm_n^{t,i}\|}{\|\bfm_n^{t,i}\|+t\mathbbm 1_{\{0\}}(\bfm_n^{t,i})}\biggr)\biggl( \limsup_{t\to 0}\frac{\|\bfm_n^{t,i}\|+t\mathbbm 1_{\{0\}}(\bfm_n^{t,i})}{|t|}\biggr)=0.
          \end{split}
       \end{equation}}
       \argument{\lref{def: Theta};\lref{def: bfu};\lref{ntp: assume};\lref{ntp: eq1};\lref{ntp: eq6}}{that for all $i\in\{1,2,\dots,\nb\}$ it holds that 
       \begin{equation}\llabel{ntp: eq8'}
       \begin{split}
           \limsup_{t\to 0}\biggl\|\frac{\Theta_{n}^{t,i}-t \bfu_{n}\bfe^i}{t}\biggr\|&=\limsup_{t\to 0}\biggl\|\frac{\Theta^{t,i}_{n-1}-\gamma
           \Phi_i(\bfm^{t,i}_n)-t \bfe^i(\bfu_{n-1}-\gamma\bfb \bfv_{n})}{t}\biggr\|\\
           &\leq \biggl(\limsup_{t\to 0}\biggl\|\frac{\Theta_{n-1}^{t,i}-t \bfu_{n-1}\bfe^i}{t}\biggr\|\biggr)+\gamma\biggl(\limsup_{t\to 0}\biggl\|\frac{\Phi_i(\bfm_n^{t,i})-\bfb\bfm_n^{t,i}}{t}\biggr\|\biggr)\\
           &+\gamma\bfb\biggl(\limsup_{t\to 0}\biggl\|\frac{\bfm_n^{t,i}-t\bfv_n\bfe^i}{t}\biggr\|\biggr)\\
           &=0.
           \end{split}
       \end{equation}}
       \argument{\lref{ntp: eq8'};}{that
       \begin{equation}\llabel{ntp: eq8}
           \limsup_{t\to 0}\biggl\|\frac{\Theta_{n}^{t}-t \bfu_{n}\bfe}{t}\biggr\|=0.
       \end{equation}}
       \argument{\lref{ntp: eq8};\lref{ntp: assume};\lref{ntp: eq1};induction}{\lref{need to prove}\dott}
       \startnewargseq
       \argument{\lref{def: N};\lref{need to prove};}{that
       \begin{equation}\llabel{eq5}
           \lim_{t\to 0}\biggl\|\frac{\Theta_\ell^{t}}{t}\biggr\|=|\bfu_\ell|\geq 2\fC\varrho^\ell.
       \end{equation}}
       \argument{\lref{eq5};}{that there exists $t\in(-\delta,\delta)$ such that
       \begin{equation}\llabel{eq6}
           \frac{\|\Theta_\ell^{t}\|}{|t|}>\fC\varrho^\ell.
       \end{equation}}
       \argument{\lref{eq6};\lref{def: Theta}}{that there exists $t\in(-\delta,\delta)$ such that
       \begin{equation}\llabel{eq7}
           \|\Theta_\ell^{t}\|>\fC |t| \varrho^\ell=\fC \|\Theta_0^t\| \varrho^\ell.
       \end{equation}}
       \argument{\lref{eq7};\lref{def: bfm};\lref{def: Theta}; the fact that for all $\theta\in \cd$ it holds that $(\nabla\smalll)(\theta)=2r\theta$; the fact that $\lambda=2r(1-\alpha)$}{\lref{item 1,item 2}\dott}
\end{aproof}
\begin{tcolorbox}[colback=white!95!gray,
                  colframe=black,
                  boxrule=0.5pt,
                  sharp corners,
                  enhanced,
                  breakable,
                 ]
\begin{athm}{cor}{main theorem lower bound}
     Let $\nb\in\N$, $\bfb\in (0,\infty)$, $\chi\in (1,\infty)$, $\varrho\in (0,\frac{\sqrt{ \chi } - 1 }  { \sqrt{ \chi } + 1 })$, $\alpha=(\frac{\sqrt{ \chi } - 1 }  { \sqrt{ \chi } + 1 })^2$, for every $i\in\{1,2,\dots,\nb\}$ let $d_i,\fd_i\in\N$, $\Phi_i\in C^1(\R^{d_i\times\fd_i},\R^{d_i\times\fd_i})$ satisfy for all $\theta\in \R^{d_i\times\fd_i}$ that $(\Phi_i)'(0)\theta-\bfb\theta=0=\Phi_i(0)$, and let $\cd=\times_{i=1}^{\nb}\R^{d_i\times\fd_i}$. Then there exist $\convexcoe\in (0,\infty)$, $\vartheta\in\cd$ and $\smalll\in C^2( \cd,\R)$ with $(\nabla \smalll )( \vartheta ) = 0$ and $
			\forall\, \theta\in\cd\colon \convexcoe\|\theta\|^2\leq \smalll'' (\vartheta)(\theta,\theta)\leq \convexcoe\chi \| \theta \|^2
			$ such that for all $\delta,\Gamma,\fC \in (0,\infty)$ there exist $\gamma\colon\N\to (\Gamma-\delta,\Gamma+\delta)$, $\bfm=(\bfm^1,\dots,\bfm^\nb)\colon \mathbb{N}_0\to\cd$, and $\Theta=(\Theta^1,\dots,\Theta^\nb)\colon \mathbb{N}_0\to\cd$  such that
\begin{enumerate}[label=(\roman*)]
    \item \llabel{item 1} it holds for all $n \in \mathbb{N}$, $i\in\{1,2,\dots,\nb\}$  that
         \begin{equation}
            \textcolor{magenta}{\bfm_0=0},\qquad \textcolor{magenta}{\bfm_n= \alpha \bfm_{n-1}+(1-\alpha)(\nabla\smalll)(\Theta_{n-1})},
            \end{equation}
            \begin{equation}
            \textcolor{magenta}{\|\Theta_0-\vartheta\|<\delta},\qquad\text{and}\qquad  \textcolor{magenta}{\Theta_n^i=\Theta_{n-1}^i-\gamma_n
          \Phi_i(\bfm_n^i)}
         \end{equation}
         and
         \item \llabel{item 2} there exists $n\in\N$ such that $\textcolor{magenta}{ \textstyle
   \textstyle\|\Theta_n-\vartheta\|>\fC\|\Theta_0-\vartheta\|\varrho^n}$.
\end{enumerate} 
\end{athm}
\end{tcolorbox}
    \begin{aproof}
    Throughout this proof let $\convexcoe\in \R$, $\vartheta\in\cd$ satisfy $\convexcoe=2$ and $\vartheta=0$ and let $\smalll\colon\cd\to\R$ satisfy for all $\theta\in\cd$ that $\smalll(\theta)=\|\theta\|^2$.
\argument{\cref{lem: lower bound};}{\lref{item 1,item 2}\dott}
    \end{aproof}
\subsection{Sharp convergence rates of momentum and MUON}\label{subsec: sharp muon}
 \begin{tcolorbox}[colback=white!95!gray,
                  colframe=black,
                  boxrule=0.5pt,
                  sharp corners,
                  enhanced,
                  breakable,
                 ]
                 \begin{athm}{theorem}{theo: sharp rate}
                 Let $\nb\in \mathbb{N}$, $\bfb\in (0,\infty)$,  $\chi\in (1,\infty)$, $\alpha=(\frac{\sqrt{ \chi } - 1 }  { \sqrt{ \chi } + 1 })^2$, $\varrho=\frac{\sqrt{ \chi } - 1 }  { \sqrt{ \chi } + 1 }$, for every $i\in\{1,2,\dots,\nb\}$ let $d_i,\fd_i\in\N$, $\Phi_i\in C^1(\R^{d_i\times\fd_i},\R^{d_i\times\fd_i})$ satisfy for all $\theta\in \R^{d_i\times\fd_i}$ that $(\Phi_i)'(0)\theta-\bfb\theta=0=\Phi_i(0)$, let $\cd=\times_{i=1}^{\nb}\R^{d_i\times\fd_i}$, and for every $n\in\N$ let $\varphi_n=(\varphi_n^1,\dots,\varphi_n^\nb)\colon \cd^n\to\cd$ satisfy for all $\theta=((\theta_k^i)_{i\in\{1,2,\dots,\nb\}})_{k\in\{1,2,\dots,n\}}\in \cd^n$, $i\in\{1,2,\dots,\nb\}$ that
\begin{equation}\llabel{def: Phi}
\textcolor{magenta}{\textstyle\varphi_n^i(\theta)=\Phi_i\bigl((1-\alpha)\sum_{k=1}^{n}\alpha^{n-k}\theta_k^i\bigr)}.
\end{equation}
Then it holds for all $\rho \in (0,\infty) \backslash \{ \varrho \}$ that the following three statements are equivalent:
             \begin{enumerate}[label=(\roman*)]
                 \item \label{item 1: sharp rate} It holds that $\rho > 
                 \varrho$.
                 \item \label{item 2: sharp rate} For all $\convexcoe \in (0,\infty)$, $\cL \in C^2( \cd, \R )$, $\vartheta\in(\nabla\cL)^{-1}(\{0\})$ with $
			\forall\, \theta\in\cd\colon \convexcoe\|\theta\|^2\leq \cL'' (\vartheta)(\theta,\theta)\leq \convexcoe\chi \| \theta \|^2
			$ there exist $\delta,\Gamma,\fC \in (0,\infty)$ such that for every $\gamma=(\gamma_n)_{n\in\N}\colon\N\to (\Gamma-\delta,\Gamma+\delta)$, every $\bfm=(\bfm^1,\dots,\bfm^\nb)\colon \mathbb{N}_0\to\cd$, and every $\Theta=(\Theta^1,\dots,\Theta^\nb)\colon \mathbb{N}_0\to\cd$ with the property that 
    for all $n \in \mathbb{N}$, $i\in\{1,2,\dots,\nb\}$ it holds that
         \begin{equation}
            \textcolor{magenta}{\bfm_0=0},\qquad \textcolor{magenta}{\bfm_n= \alpha \bfm_{n-1}+(1-\alpha)(\nabla\cL)(\Theta_{n-1})},
            \end{equation}
            \begin{equation}
            \textcolor{magenta}{\|\Theta_0-\vartheta\|< \delta},\qquad\text{and}\qquad  \textcolor{magenta}{\Theta_n^i=\Theta_{n-1}^i-\gamma_n
           \Phi_i(\bfm_n^i)}
         \end{equation}
         we have for all $n\in\N_0$ that $\textcolor{magenta}{ \textstyle
   \textstyle\|\Theta_n-\vartheta\|\leq\fC\rho^n\|\Theta_0-\vartheta\|}$.
                 \item \label{item 3: sharp rate} For all $\convexcoe \in (0,\infty)$, $\cL \in C^2(\cd, \R )$, $\vartheta\in(\nabla\cL)^{-1}(\{0\})$ with $
			\forall\, \theta\in\cd\colon \convexcoe\|\theta\|^2\leq \cL'' (\vartheta)(\theta,\theta)\leq \convexcoe\chi \| \theta \|^2
			$ there exist $\delta,\Gamma, \fC \in (0,\infty) $ such that 
                  for every $\gamma=(\gamma_n)_{n\in\N}\colon\N\to (\Gamma-\delta,\Gamma+\delta)$ and every $ \Theta \colon \N_0\to\cd $
			with $\| \Theta_0 - \vartheta \| < \delta$ and
			\begin{equation}
				 \textcolor{magenta}{\forall \, n \in \N \colon
				\Theta_n
				=
				\Theta_{ n - 1 }
				-
				\gamma_n \varphi_n\bigl(
				( \nabla \cL )( \Theta_0 ) ,
				( \nabla \cL )( \Theta_1 ) ,
				\dots ,
				( \nabla \cL )( \Theta_{ n - 1 } )
				\bigr)}
			\end{equation}
			we have for all $n\in\N_0$ that
			$
			\textcolor{magenta}{\| \Theta_n - \vartheta \|
			\leq \fC\rho^{  n } \|\Theta_0-\vartheta\|}
			$.
             \end{enumerate}
                 \end{athm}
                 \end{tcolorbox}
    \begin{aproof}
        Throughout this proof for every $\Theta\colon\N_0\to\cd$ and $\cL\in C^1(\cd,\allowbreak\R)$ let $\bfm(\Theta,\cL)=(\bfm_n^i(\Theta,\cL))_{(i,n)\in\{1,2,\dots,\nb\}\times\N_0}\colon\N_0\to\cd$ satisfy for all $n\in\N_0$ that
        \begin{equation}\llabel{def: bfm}
            \bfm_n(\Theta,\cL)=(1-\alpha)\sum_{k=1}^n\alpha^{n-k}(\nabla\cL)(\Theta_{k-1}).
        \end{equation}
        \argument{\cref{theo: abstract optimal rate}}{that (\cref{item 1: sharp rate}$\rightarrow$\cref{item 2: sharp rate})\dott}
        \startnewargseq
         \argument{\cref{main theorem lower bound}}{that (\cref{item 2: sharp rate}$\rightarrow$\cref{item 1: sharp rate})\dott}
        \startnewargseq
        \argument{\lref{def: bfm};}{that for all $n\in\N$, $\cL\in C^1(\cd,\R)$ and all $\Theta\colon\N_0\to\cd$ it holds that
        \begin{equation}\llabel{eq1}
            \bfm_{n}(\Theta,\cL)=\alpha\bfm_{n-1}(\Theta,\cL)+(1-\alpha)(\nabla\cL)(\Theta_{n-1}).
        \end{equation}}
       \argument{\lref{def: Phi};\lref{def: bfm}}{that for all $n\in\N$, $\cL=(\cL(\theta))_{\theta=(\theta_1,\dots,\theta_\nb)}\in C^1(\cd,\R)$, $i\in\{1,2,\dots,\nb\}$ and all $\Theta\colon\N_0\to\cd$ it holds that
       \begin{equation}\llabel{eq2}
       \begin{split}
        \varphi_n^i\bigl(
				( \nabla \cL )( \Theta_0 ) ,
				( \nabla \cL )( \Theta_1 ) ,
				\dots ,
				( \nabla \cL )( \Theta_{ n - 1 } )
				\bigr)
                &\textstyle = \Phi_i\bigl((1-\alpha)\sum_{k=1}^{n}\alpha^{n-k}(\nabla_{\theta_i}\cL)(\Theta_{k-1})\bigr)\\
                &=\Phi_i(\bfm_n^i(\Theta,\cL)).
                \end{split}
       \end{equation}}
       \argument{\lref{eq2};\lref{eq1}}{that (\cref{item 2: sharp rate}$\rightarrow$\cref{item 3: sharp rate})\dott}
       \startnewargseq
       \argument{\unskip, \eg, \cite[Lemma 2.7]{DeRoArAd2025}}{that for every $\cL\in C^2(\cd,\R)$, every $\bfm\colon\N_0\to\cd$, and every $\Theta\colon\N_0\to\cd$, with the property that for all $n\in\N$ it holds that $ \bfm_0=0$ and $\bfm_n= \alpha \bfm_{n-1}+(1-\alpha)(\nabla\cL)(\Theta_{n-1})$ we have for all $n\in\N_0$ that
       \begin{equation}\llabel{eq4}
           \textstyle \bfm_n=(1-\alpha)\sum_{k=1}^n\alpha^{n-k}(\nabla\cL)(\Theta_{k-1}).
       \end{equation}}
       \argument{\lref{eq4};\lref{def: Phi}}{that for every $\cL=(\cL(\theta))_{\theta=(\theta_1,\dots,\theta_\nb)}\in C^2(\cd,\R)$, every $\bfm=(\bfm^1,\dots,\bfm^\nb)\colon\allowbreak\N_0\to\cd$, and every $\Theta\colon\N_0\to\cd$ with the property that for all $n\in\N$ it holds that $ \bfm_0=0$ and $\bfm_n= \alpha \bfm_{n-1}+(1-\alpha)(\nabla\cL)(\Theta_{n-1})$ we have for all $n\in\N$, $i\in\{1,2,\dots,\nb\}$ that
       \begin{equation}\llabel{eq5}
       \begin{split}
           \varphi^i_n\bigl(
				( \nabla \cL )( \Theta_0 ) ,
				( \nabla \cL )( \Theta_1 ) ,
				\dots ,
				( \nabla \cL )( \Theta_{ n - 1 } )
				\bigr)&=\textstyle\Phi_i\bigl((1-\alpha)\sum_{k=1}^n\alpha^{n-k}
				( \nabla_{\theta_i} \cL )( \Theta_{k-1} )
				\bigr)\\
                &=\Phi_i(\bfm_n^i).
                \end{split}
       \end{equation}}
       \argument{\lref{eq5};}{that (\cref{item 3: sharp rate}$\rightarrow$\cref{item 2: sharp rate})\dott}
    \end{aproof}
            \begin{samepage}
 \begin{tcolorbox}[colback=white!95!gray,
                  colframe=black,
                  boxrule=0.5pt,
                  sharp corners,
                  enhanced,
                  breakable,
                 ]
                 \begin{athm}{cor}{cor: sharp rate}
                 Let $\nb\in \mathbb{N}$, $\bfb\in (0,\infty)$, $\chi\in (1,\infty)$, $\alpha=(\frac{\sqrt{ \chi } - 1 }  { \sqrt{ \chi } + 1 })^2$, $\varrho=\frac{\sqrt{ \chi } - 1 }  { \sqrt{ \chi } + 1 }$, for every $i\in\{1,2,\dots,\nb\}$ let $d_i,\fd_i\in\N$, $\Phi_i\in C^1(\R^{d_i\times\fd_i},\R^{d_i\times\fd_i})$ satisfy for all $\theta\in \R^{d_i\times\fd_i}$ that $(\Phi_i)'(0)\theta-\bfb\theta=0=\Phi_i(0)$, let $\cd=\times_{i=1}^{\nb}\R^{d_i\times\fd_i}$, and for every $n\in\N$ let $\varphi_n=(\varphi_n^1,\dots,\varphi_n^\nb)\colon \cd^n\to\cd$ satisfy for all $\theta=((\theta_k^i)_{i\in\{1,2,\dots,\nb\}})_{k\in\{1,2,\dots,n\}}\in \cd^n$, $i\in\{1,2,\dots,\nb\}$ that
\begin{equation}\llabel{def: Phi}
\textcolor{magenta}{\textstyle\varphi_n^i(\theta)=\Phi_i\bigl((1-\alpha)\sum_{k=1}^{n}\alpha^{n-k}\theta_k^i\bigr)}.
\end{equation}
Then $(\varphi_n)_{n\in\N}$ converges with $\chi$-rate $\varrho$ (cf.\ \cref{definition: rate of convergence}).
                 \end{athm}
                 \end{tcolorbox}
                 \end{samepage}
                 \begin{aproof}
                     \argument{\cref{theo: sharp rate};}{that $(\varphi_n)_{n\in\N}$ converges with $\chi$-rate $\varrho$\dott}
                 \end{aproof}
  \begin{samepage}
 \begin{tcolorbox}[colback=white!95!gray,
                  colframe=black,
                  boxrule=0.5pt,
                  sharp corners,
                  enhanced,
                  breakable,
                 ]
                 \begin{athm}{cor}{cor: sharp rate muon}
                 Let $K\in\N_0$, $\nb,d_1,d_2,\dots,d_\nb,\fd_1,\fd_2,\dots,\fd_\nb\in\N$, $\cd=\times_{i=1}^{\nb}\R^{d_i\times\fd_i}$, $\chi\in (1,\infty)$, $\alpha=(\frac{\sqrt{ \chi } - 1 }  { \sqrt{ \chi } + 1 })^2$, $\varrho=\frac{\sqrt{ \chi } - 1 }  { \sqrt{ \chi } + 1 }$, let $
\nscoe=(\nscoe_{i,j})_{(i,j)\in (\mathbb{N}_0)^2}\colon(\mathbb{N}_0)^2\allowbreak\to \mathbb{R}$ satisfy $\textcolor{magenta}{\min\{\nscoe_{0,0},\nscoe_{0,1}\}>0}$ and $\textcolor{magenta}{\#(\nscoe^{ - 1 }( \mathbb{R}\backslash\{0\} )) <\infty}$, assume for all $i\in \mathbb{N}\cap[0,K]$ that $\textcolor{magenta}{\nscoe_{i,0}>0}$, and for every $n\in\N$ let $\varphi_n=(\varphi_n^1,\dots,\varphi_n^\nb)\colon \cd^n\to\cd$ satisfy for all $\theta=((\theta_k^i)_{i\in\{1,2,\dots,\nb\}})_{k\in\{1,2,\dots,n\}}\in \cd^n$, $i\in\{1,2,\dots,\nb\}$ that
\begin{equation}\llabel{def: Phi}
\textcolor{magenta}{\textstyle\varphi_n^i(\theta)=\TbigNewtonSchulzAlgorithm{\nscoe}{K}{(1-\alpha)\sum_{k=1}^{n}\alpha^{n-k}\theta_k^i}}
\end{equation}
(cf.\ \cref{definition: NS}). Then $(\varphi_n)_{n\in\N}$ converges with $\chi$-rate $\varrho$ (cf.\ \cref{definition: rate of convergence}).
                 \end{athm}
                 \end{tcolorbox}
                 \end{samepage}
                 \begin{aproof}
     \argument{\cref{lem: NS represent};\cref{cor: sharp rate};the fact that $(\nscoe_{0,1})^{-1}\prod_{k=1}^K\nscoe_{k,0}>0$}{that $(\varphi_n)_{n\in\N}$ converges with $\chi$-rate $\varrho$\dott}
                 \end{aproof}
 \subsection{Approximative orthogonalization methods within MUON}\label{sec: practical}
In this subsection we show that \cref{cor: sharp rate muon} covers \MUON\ with $K\in\N$ steps of the classical \NS\ method (see \cref{lem: abc verify 0} below; cf.\ \cite{MR295547} and \cite[Section~8.3]{MR2396439}), with five \NS\ steps with the polynomial proposed in \cite{jordan2024muon} (see \cref{lem: abc verify} below; cf.\ \cite[Lemma~5.1]{DoDeAr2026muon}), and with the Polar Express method \cite{ar2505.16932} as implemented in \cite{Jordanmuoncode} (see \cref{lem: abc verify 2} below; cf.\ \cite[Lemma~5.2]{DoDeAr2026muon}).

 \begin{tcolorbox}[colback=white!95!gray,
                  colframe=black,
                  boxrule=0.5pt,
                  sharp corners,
                  enhanced,
                  breakable,
                 ]
                 \begin{athm}{lemma}{lem: abc verify 0}[\textcolor{red}{classical \MUON}]
                    Let $K\in\N$, $\varepsilon\in (0,\infty)$, $\nscoe=(\nscoe_{i,j})_{(i,j)\in (\N_0)^2}\colon(\N_0)^2\to\R$ satisfy for all $i, j \in \N_0$, $k \in \{0,1,\dots,K\}$ with $ j \geq 2$ or $i>K$ that $\nscoe_{ i, j } = 0$ and
                   \begin{equation}\llabel{def: a}
                   \begin{split}
                      (\nscoe_{k,0},\nscoe_{k,1})=
                       \begin{cases}
                       (1, \varepsilon) &\colon k=0,\\
(\frac 32,-\frac 12) &\colon k>0.
                       \end{cases}
                       \end{split}
                   \end{equation}
                    Then 
                    \begin{enumerate}[label=(\roman*)]
                        \item \llabel{item 1} it holds that $\min\{\nscoe_{0,0},\nscoe_{0,1}\}>0$ and $\#(\nscoe^{ - 1 }( \R\backslash\{0\} )) <\infty$ and
                     \item \llabel{item 2} it holds for all $n\in \{1,2,\dots,K\}$ that $\nscoe_{n,0}>0$.
                    \end{enumerate}
                 \end{athm}
                 \end{tcolorbox}
\begin{aproof}
     \argument{\lref{def: a}}{\lref{item 1,item 2}\dott} 
\end{aproof}
\begin{tcolorbox}[colback=white!95!gray,
                  colframe=black,
                  boxrule=0.5pt,
                  sharp corners,
                  enhanced,
                  breakable,
                 ]
                 \begin{athm}{lemma}{lem: abc verify}[\textcolor{red}{5 \NS\ steps with the original \MUON\ polynomial \cite{jordan2024muon}}]
                    Let $\nscoe=(\nscoe_{i,j})_{(i,j)\in (\N_0)^2}\colon(\N_0)^2\to\R$ satisfy for all $i, j \in \N_0$, $k \in \{0,1,\dots,5\}$ with $\max\{ i, 2 j \} \geq 6$ that $\nscoe_{ i, j } = 0$ and
                   \begin{equation}\llabel{def: a}
                   \begin{split}
                      (\nscoe_{k,0},\nscoe_{k,1},\nscoe_{k,2})=
                       \begin{cases}
                       (1, 10^{-7},0) &\colon k=0,\\
(3.4445,-4.7750,2.0315) &\colon k>0.
                       \end{cases}
                       \end{split}
                   \end{equation}
                    Then 
                    \begin{enumerate}[label=(\roman*)]
                        \item \llabel{item 1} it holds that $\min\{\nscoe_{0,0},\nscoe_{0,1}\}>0$ and $\#(\nscoe^{ - 1 }( \R\backslash\{0\} )) <\infty$ and
                     \item \llabel{item 2} it holds for all $n\in \{1,2,\dots,5\}$ that $\nscoe_{n,0}>0$.
                    \end{enumerate}
                 \end{athm}
                 \end{tcolorbox}
\begin{aproof}
     \argument{\lref{def: a}}{\lref{item 1,item 2}\dott} 
\end{aproof}

  \begin{tcolorbox}[colback=white!95!gray,
                  colframe=black,
                  boxrule=0.5pt,
                  sharp corners,
                  enhanced,
                  breakable,
                 ]
                 \begin{athm}{lemma}{lem: abc verify 2}[\textcolor{red}{Polar Express method \cite{ar2505.16932} as implemented in \cite{Jordanmuoncode}}]
                   Let $\nscoe=(\nscoe_{i,j})_{(i,j)\in (\N_0)^2}\colon(\N_0)^2\to\R$ satisfy for all $i, j \in \N_0$, $k \in \{0,1,\dots,5\}$ with $\max\{ i, 2 j \} \geq 6$ that $\nscoe_{ i, j } = 0$ and
                   \begin{equation}\llabel{def: a}
                   \begin{split}
                     & (\nscoe_{k,0},\nscoe_{k,1},\nscoe_{k,2})\\
                     &=
                       \begin{cases}
                       (1.02, 10^{-6},0) &\colon k=0,\\
(8.156554524902461,-22.48329292557795,15.878769915207462) &\colon k=1,\\
(4.042929935166739,-2.808917465908714,0.5000178451051316) &\colon k=2,\\
(3.8916678022926607,-2.772484153217685,0.5060648178503393) &\colon k=3,\\
(3.285753657755655,-2.3681294933425376,0.46449024233003106) &\colon k=4,\\
(2.3465413258596377,-1.7097828382687081,0.42323551169305323) &\colon k=5.
                       \end{cases}
                       \end{split}
                   \end{equation}
                   
                   Then  \begin{enumerate}[label=(\roman*)]
                        \item \llabel{item 1} it holds that $\min\{\nscoe_{0,0},\nscoe_{0,1}\}>0$ and $\#(\nscoe^{ - 1 }( \R\backslash\{0\} )) <\infty$ and
                     \item \llabel{item 2} it holds for all $n\in \{1,2,\dots,5\}$ that $ \nscoe_{n,0}>0$.
                    \end{enumerate}
                 \end{athm}
                 \end{tcolorbox}
                 \begin{aproof}
                     \argument{\lref{def: a};}{\lref{item 1,item 2}\dott} 
                 \end{aproof}
\section{Optimal convergence rates of the GD and the RMSprop optimizers}
\label{sec: rmsprop}
In this section we establish in \cref{cor: sharp rate rmsprop} below that for every $d\in\N\backslash\{1\}$ the \RMSprop\ optimizer on $\R^d$ converges with $\chi$-rate $(\chi+1)^{-1}(\chi-1)$. \Cref{cor: sharp rate rmsprop} also includes the standard \GD\ optimizer as a special case by choosing $\beta=\varepsilon=1$.

\subsection{GD and RMSprop optimizers}

In the next \cref{definition: GD} we recall the concept of the standard \GD\ optimizer (see \cite{Cauchy1847}) using the functions $\Phi_n\colon(\R^d)^n\to\R^d$, $n\in\N$, from \cref{definition: rate of convergence} above; see, \eg, \cite[Definition~2.1]{DeArAdsharp2025}.
\begin{tcolorbox}[colback=white!95!gray,
                  colframe=black,
                  boxrule=0.5pt,
                  sharp corners,
                  enhanced,
                  breakable,
                 ]
\begin{definition}[\GD]
	\label{definition: GD}
	Let $ d \in \N $ and let
	$ \Phi_n \colon ( \R^{ d } )^n \to \R^{ d } $,
	$ n \in \N $, be functions. Then we say that
    $ ( \Phi_n )_{ n \in \N } $ is the \GD\ optimizer
    on $ \R^{ d } $ (we say that
    $ ( \Phi_n )_{ n \in \N } $ is the \GD\ optimizer)
    if and only if it holds for all
    $ n \in \N $, $ g_1, g_2, \dots, g_n \in \R^{ d } $
    that
    \begin{equation}
	  \Phi_n( g_1, g_2, \dots, g_n )
	  =
	  g_n
	  .
    \end{equation}
	\end{definition}
\end{tcolorbox}
In the next \cref{definition: RMSprop} we recall the concept of the \RMSprop\ optimizer (see \cite{Hinton24_RMSprop}) using the functions $\Phi_n\colon(\R^d)^n\to\R^d$, $n\in\N$, from \cref{definition: rate of convergence} above; cf., \eg, \cite[Definition~2.2]{DeArAdsharp2025}.
\begin{tcolorbox}[colback=white!95!gray,
                  colframe=black,
                  boxrule=0.5pt,
                  sharp corners,
                  enhanced,
                  breakable,
                 ]
    \begin{definition}[\RMSprop]
	\label{definition: RMSprop}
	Let $ d \in \N $,
	$ \beta \in (0,1) $,
	$ \varepsilon \in (0,\infty) $
	and let
	$
      \Phi_n = ( \Phi_n^1, \dots, \Phi_n^{ d } )
	  \colon
	  \allowbreak
	  ( \R^{ d } )^n \to \R^{ d }
	$,
	$ n \in \N $, be functions. Then we say that
    $ ( \Phi_n )_{ n \in \N } $ is the
    $ \beta $-$ \varepsilon $-\RMSprop\ optimizer
    on $ \R^{ d } $ (we say that
    $ ( \Phi_n )_{ n \in \N } $ is
    the $ \beta $-$ \varepsilon $-\RMSprop\ optimizer)
    if and only if it holds for all
    $ n \in \N $,
    $
      g =
      (
        ( g_{ i, j } )_{ j \in \{ 1, 2, \dots, d \} }
      )_{ i \in \{ 1, 2, \dots, n \} } \in ( \R^{ d } )^n
    $,
    $ j \in \{ 1, 2, \dots, d \} $
    that
    \begin{equation}
    \textstyle
	  \Phi_n^j( g )
	  =
	  \bigl[
	    \varepsilon
	    +
	    (
	      ( 1 - \beta )
	      \sum_{ i = 1 }^n
	      \beta^{ n - i }
	      | g_{ i, j } |^2
	    )^{ \nicefrac{ 1 }{ 2 } }
	  \bigr]^{ - 1 }
	  g_{ n, j }
	  .
    \end{equation}
	\end{definition}
    \end{tcolorbox}
\subsection{Upper bounds on convergence rates of RMSprop}
\begin{samepage}
 \begin{tcolorbox}[colback=white!95!gray,
                  colframe=black,
                  boxrule=0.5pt,
                  sharp corners,
                  enhanced,
                  breakable,
                 ]
\begin{athm}{theorem}{theo: higher bound rmsprop}
      Let $d\in \mathbb{N}$,  $\chi\in (1,\infty)$, $\convexcoe,\varepsilon\in (0,\infty)$, $\varrho\in (\frac{ \chi  - 1 }  { \chi  + 1 },1)$, $\Gamma=\frac{2\varepsilon}{\convexcoe(1+\chi)}$, $\beta\in (0,1]$, $\vartheta\in\R^d$ and let $\smalll=(\smalll(\theta))_{\theta=(\theta_1,\dots,\theta_d)\in\R^d}\in \allowbreak C^{2}(\R^d,\mathbb{R})$ satisfy for all $\theta \in \R^d$ that \begin{equation}\llabel{def: bbL}
     \textcolor{magenta}{(\nabla\smalll)(\vartheta)=0}\qqandqq\textcolor{magenta}{\convexcoe\|\theta\|^2\leq \theta^\top(\operatorname{Hess}\smalll) (\vartheta)\theta\leq \convexcoe\chi \| \theta \|^2}.
 \end{equation}  
 Then 
there exist $\fC,\delta \in (0,\infty)$ such that for every $\gamma\colon\N\to(\Gamma-\delta,\Gamma+\delta)$, every $\bbV=(\bbV^1,\dots,\bbV^{d})\colon \mathbb{N}_0\to\R^d$, and every $\Theta=(\Theta^1,\dots,\Theta^{d})\colon \mathbb{N}_0\to\R^d$ with the property that for all $n \in \mathbb{N}$, $i\in\{1,2,\dots,d\}$ it holds that 
         \begin{equation}
            \textcolor{magenta}{\bbV_0=0},\qquad \textcolor{magenta}{\bbV_n^i= \beta \bbV_{n-1}^i+(1-\beta)|(\nabla_{\theta_i}\smalll)(\Theta_{n-1})|^2},
            \end{equation}
            \begin{equation}
            \textcolor{magenta}{\|\Theta_0-\vartheta\|<\delta},\qquad\text{and}\qquad  \textcolor{magenta}{\Theta_n^i=\Theta_{n-1}^i-\gamma_n
          [\varepsilon+(\bbV_n^i)^{1/2}]^{-1}(\nabla_{\theta_i}\smalll)(\Theta_{n-1})}
         \end{equation}
          we have for all $n\in\N_0$ that
$    \textcolor{magenta}{ \textstyle
    \|\Theta_n-\vartheta\|\leq \fC\|\Theta_0-\vartheta\|\varrho^{n}}$.
\end{athm}
\end{tcolorbox}
\end{samepage}
\begin{aproof}
     Throughout this proof let $\rho=\frac{\chi-1}{\chi+1}$, for every $\xi\in\R^d$ and every $\gamma\colon\N\to\R$ let $\bbV^{\xi,\gamma}=(\bbV^{\xi,\gamma,1},\dots,\bbV^{\xi,\gamma,d})\colon\N_0\to\R^d$ and $\Theta^{\xi,\gamma}=(\Theta^{\xi,\gamma,1},\dots,\Theta^{\xi,\gamma,d})\colon\N_0\to\R^d$ satisfy for all $n\in\N$, $i\in\{1,2,\dots,d\}$ that
     \begin{equation}\llabel{def: bfm}
      \bbV^{\xi,\gamma}_0=0,\qquad \bbV^{\xi,\gamma,i}_n= \beta \bbV^{\xi,\gamma,i}_{n-1}+(1-\beta)|(\nabla_{\theta_i}\smalll)(\Theta^{\xi,\gamma}_{n-1})|^2,
            \end{equation}
            \begin{equation}\llabel{def: Theta}
        \Theta^{\xi,\gamma}_0=\xi,\qquad\text{and}\qquad \Theta_n^{\xi,\gamma,i}=\Theta_{n-1}^{\xi,\gamma,i}-\gamma_n
          [\varepsilon+(\bbV_n^{\xi,\gamma,i})^{1/2}]^{-1}(\nabla_{\theta_i}\smalll)(\Theta_{n-1}^{\xi,\gamma}),
         \end{equation}
         let $I=(I_{i,j})_{(i,j)\in\{1,2,\dots,d\}\times\{1,2,\dots,d\}}\in\R^{d\times d}$ satisfy for all $i,j\in\{1,2,\dots,d\}$ that $I_{i,j}=\mathbbm 1_{\{j\}}(i)$,
      and for every $n\in \N$, $\xi\in\R^d$ and every $\gamma\colon\N\to\R$ let $\Lambda^{\xi,\gamma}_n=(\Lambda^{\xi,\gamma,i,j}_n)_{(i,j)\in\{1,2,\dots,d\}\times\{1,2,\dots,d\}}\in\R^{d\times d}$ satisfy for all $i,j\in\{1,2,\dots,d\}$ that
      \begin{equation}\llabel{def: Lambda}
         \Lambda^{\xi,\gamma,i,j}_n=\gamma_n[\varepsilon+(\bbV_n^{\xi,\gamma,i})^{1/2}]^{-1}\mathbbm 1_{\{j\}}(i). 
      \end{equation}
      \argument{\lref{def: Theta};\lref{def: Lambda}}{that for every $n\in\N$, $\xi\in\R^d$ and every $\gamma\colon\N\to\R$ it holds that
      \begin{equation}\llabel{eq1}
      \begin{split}
         & \Theta_n^{\xi,\gamma}-\vartheta=\Theta_{n-1}^{\xi,\gamma}-\vartheta-\Lambda_n^{\xi,\gamma}(\nabla\smalll)(\Theta_{n-1}^{\xi,\gamma})\\
          &=\Theta_{n-1}^{\xi,\gamma}-\vartheta-\Gamma\varepsilon^{-1} I(\operatorname{Hess}\smalll)(\vartheta)(\Theta_{n-1}^{\xi,\gamma}-\vartheta)-\Gamma\varepsilon^{-1} I((\nabla\smalll)(\Theta_{n-1}^{\xi,\gamma})-(\operatorname{Hess}\smalll)(\vartheta)(\Theta_{n-1}^{\xi,\gamma}-\vartheta))\\
          &-(\gamma_n-\Gamma)\varepsilon^{-1}I(\nabla\smalll)(\Theta_{n-1}^{\xi,\gamma})+(\gamma_n\varepsilon^{-1}I-\Lambda_n^{\xi,\gamma})(\nabla\smalll)(\Theta_{n-1}^{\xi,\gamma})\\
          &=(I-\Gamma\varepsilon^{-1}(\operatorname{Hess}\smalll)(\vartheta))(\Theta_{n-1}^{\xi,\gamma}-\vartheta)-\Gamma\varepsilon^{-1} I((\nabla\smalll)(\Theta_{n-1}^{\xi,\gamma})-(\operatorname{Hess}\smalll)(\vartheta)(\Theta_{n-1}^{\xi,\gamma}-\vartheta))\\
          &-(\gamma_n-\Gamma)\varepsilon^{-1}I(\nabla\smalll)(\Theta_{n-1}^{\xi,\gamma})+(\gamma_n\varepsilon^{-1}I-\Lambda_n^{\xi,\gamma})(\nabla\smalll)(\Theta_{n-1}^{\xi,\gamma}).
          \end{split}
      \end{equation}}
      \argument{the assumption that $\Gamma=\frac{2\varepsilon}{\convexcoe(1+\chi)}$}{that for all $\lambda\in \bbC$, $\theta\in\bbC^d\backslash\{0\}$ with $(I-\Gamma\varepsilon^{-1}(\operatorname{Hess}\smalll)(\vartheta))\theta=\lambda\theta$ it holds that
      \begin{equation}\llabel{eq2}
          (\operatorname{Hess}\smalll)(\vartheta)\theta=\frac{(1-\lambda)\convexcoe(1+\chi)}{2}\theta.
      \end{equation}}
      \argument{the fact that $\smalll\in C^2(\R^d,\R)$}{that $(\operatorname{Hess}\smalll)(\vartheta)$ \llabel{arg1} is a symmetric matrix\dott}
      \argument{\lref{arg1};}{that $I-\Gamma\varepsilon^{-1}(\operatorname{Hess}\smalll)(\vartheta)$ is \llabel{arg2} symmetric\dott}
      \argument{\lref{arg2};}{for all $\lambda\in \bbC$, $\theta\in\bbC^d\backslash\{0\}$ with \llabel{arg3} $(I-\Gamma\varepsilon^{-1}(\operatorname{Hess}\smalll)(\vartheta))\theta=\lambda\theta$ that $\lambda\in \R$\dott}
      \argument{\lref{def: bbL};\lref{eq2};\lref{arg3}}{for all $\lambda\in \bbC$, $\theta\in\bbC^d\backslash\{0\}$ with $(I-\Gamma\varepsilon^{-1}(\operatorname{Hess}\smalll)(\vartheta))\theta=\lambda\theta$ that $\lambda\in \R$ and
      \begin{equation}\llabel{eq3'}
       \textstyle   \convexcoe\leq \frac{(1-\lambda)\convexcoe(1+\chi)}{2}\leq \convexcoe\chi.
      \end{equation}}
      \argument{\lref{eq3'};}{for all $\lambda\in \bbC$, $\theta\in\bbC^d\backslash\{0\}$ with $(I-\Gamma\varepsilon^{-1}(\operatorname{Hess}\smalll)(\vartheta))\theta=\lambda\theta$ that
      \begin{equation}\llabel{eq4}
      |\lambda|\leq \frac{\chi-1}{\chi+1}.
      \end{equation}}
      \argument{\lref{eq4};the fact that $\varrho>\rho=\frac{\chi-1}{\chi+1}$;\unskip, \eg, \cite[Lemma 2.3]{DeArAdsharp2025}}{that there exists a norm $\nnorm{\cdot}\colon \R^d\to[0,\infty)$ which satisfies for all $\theta\in\R^d$ that
\begin{equation}\llabel{eq9}
    \nnorm{(I-\Gamma\varepsilon^{-1}(\operatorname{Hess}\smalll)(\vartheta))\theta}\leq \frac{\varrho+\rho}{2}\nnorm{\theta}.
\end{equation}}
\startnewargseq
    \argument{the fact that $\|\cdot\|$ is a norm; the fact that $\nnorm{\cdot}$ is a norm; the fact that $\R^d$ is finite dimensional vector space;}{that there exist $\mu_1,\mu_2\in (0,\infty)$ which satisfy for all $\theta\in \R^d$ that
         \begin{equation}\llabel{def: mu}
             \mu_1\nnorm{\theta}\leq \|\theta\|\leq \mu_2\nnorm{\theta}.
         \end{equation}}
         \startnewargseq
\argument{\lref{eq1};\lref{eq9};\lref{def: mu}}{that for every $n\in\N$, $\xi\in\R^d$ and every $\gamma\colon\N\to\R$ it holds that
\begin{equation}\llabel{eq10}
\begin{split}
    \nnorm{\Theta_n^{\xi,\gamma}-\vartheta}&\leq \frac{\rho+\varrho}{2}\nnorm{\Theta_{n-1}^{\xi,\gamma}-\vartheta}+\Gamma\varepsilon^{-1}(\mu_1)^{-1}\|(\nabla\smalll)(\Theta_{n-1}^{\xi,\gamma})-(\operatorname{Hess}\smalll)(\vartheta)(\Theta_{n-1}^{\xi,\gamma}-\vartheta)\|\\
    &+|\gamma_n-\Gamma|\varepsilon^{-1}(\mu_1)^{-1}\|(\nabla\smalll)(\Theta_{n-1}^{\xi,\gamma})\|+(\mu_1)^{-1}\|(\gamma_n\varepsilon^{-1}I-\Lambda_n^{\xi,\gamma})(\nabla\smalll)(\Theta_{n-1}^{\xi,\gamma})\|.
    \end{split}
\end{equation}}
 \argument{\lref{def: bbL};the fact that $\smalll\in C^2(\R^d,\R)$}{that there exist $\delta,\fC_1,\fC_2\in (0,\infty)$ such that for all $\theta\in \R^d$ with $\|\theta-\vartheta\|\leq \delta$ it holds that
         \begin{equation}\llabel{eq11}
            \fC_1\|\theta-\vartheta\|\leq \|(\nabla\smalll)(\theta)\|=\|(\nabla\smalll)(\theta)-(\nabla\smalll)(\vartheta)\|\leq \fC_2\|\theta-\vartheta\|.
         \end{equation}}
         \argument{\lref{eq11};\lref{def: mu};}{that there exist $\delta,\fC_1,\fC_2\in (0,\infty)$ such that for every $n\in\N$, $\xi\in\R^d$ and every $\gamma\colon\N\to\R$ with $\|\Theta_{n-1}^{\xi,\gamma}-\vartheta\|\leq \delta$ it holds that
         \begin{equation}\llabel{eq12}
         |\gamma_n-\Gamma|\varepsilon^{-1}(\mu_1)^{-1}\|(\nabla\smalll)(\Theta_{n-1}^{\xi,\gamma})\|\leq \fC _1|\gamma_n-\Gamma|\|\Theta_{n-1}^{\xi,\gamma}-\vartheta\|\leq\fC _2|\gamma_n-\Gamma|\nnorm{\Theta_{n-1}^{\xi,\gamma}-\vartheta} .
         \end{equation}}
    \argument{\lref{eq12};the fact that $\varrho>\rho$}{that there exists $\delta\in (0,\infty)$ such that for every $n\in\N$, $\xi\in\R^d$ and every $\gamma\colon\N\to\R$ with $\|\Theta_{n-1}^{\xi,\gamma}-\vartheta\|\leq \delta$ and $\sup_{k\in\N}|\gamma_k-\Gamma|\leq \delta$ it holds that
    \begin{equation}\llabel{evidence 1}
             |\gamma_n-\Gamma|\varepsilon^{-1}(\mu_1)^{-1}\|(\nabla\smalll)(\Theta_{n-1}^{\xi,\gamma})\|\leq \frac{\varrho-\rho}{8}\nnorm{\Theta_{n-1}^{\xi,\gamma}-\vartheta} .
    \end{equation}}
    \argument{\lref{def: bfm};\lref{eq11};\unskip, \eg, \cite[Lemma 2.7]{DeRoArAd2025}}{that there exist $\delta,\fC\in (0,\infty)$ such that for every $n\in\N$, $\xi\in\R^d$, $i\in\{1,2,\dots,d\}$, $\eta\in (0,\delta)$ and every $\gamma\colon\N\to\R$ with $\sup_{k\in\{0,1,\dots,n-1\}}\|\Theta_{k}^{\xi,\gamma}-\vartheta\|\leq \eta$ it holds that
    \begin{equation}\llabel{eq14}
        \bbV_n^{\xi,\gamma, i}=(1-\beta)\bigl[\textstyle\sum_{k=1}^n\beta^{n-k}|(\nabla_{\theta_i}\smalll)(\Theta_{k-1}^{\xi,\gamma})|^2 \textstyle \bigr]\leq (1-\beta)\bigl[\textstyle\sum_{k=1}^n\beta^{n-k}\fC\|\Theta_{k-1}^{\xi,\gamma}-\vartheta\|^2 \textstyle \bigr]\leq \fC\delta^2.
    \end{equation}}
    \argument{\lref{def: Lambda};\lref{eq14};}{that there exist $\delta,\fC\in (0,\infty)$ such that for every $n\in\N$, $\xi\in\R^d$, $\eta\in (0,\delta)$ and every $\gamma\colon\N\to\R$ with $\sup_{k\in\{0,1,\dots,n-1\}}\|\Theta_{k}^{\xi,\gamma}-\vartheta\|\leq \eta$ and $\sup_{k\in\N}|\gamma_k|\leq 2\Gamma$ it holds that  
    \begin{equation}\llabel{eq15}
    \begin{split}
        \|\gamma_n\varepsilon^{-1} I-\Lambda_n^{\xi,\gamma}\|\leq \sum_{i=1}^d|\gamma_n|\bigl|[\varepsilon+(\bbV_n^{\xi,\gamma,i})^{1/2}]^{-1}-\varepsilon^{-1}\bigr|\leq \sum_{i=1}^d|\gamma_n||\bbV_n^{\xi,\gamma,i}|^{1/2}\varepsilon^{-2}\leq 2\Gamma d\varepsilon^{-2}\fC\delta.
        \end{split}
    \end{equation}}
    \argument{\lref{eq15};\lref{def: mu};\lref{eq11}}{that there exist $\fC_1,\fC_2,\fC_3,\fC_4\in (0,\infty)$ such that for every $n\in\N$, $\xi\in\R^d$, $\delta\in (0,1)$ and every $\gamma\colon\N\to\R$ with $\sup_{k\in\{0,1,\dots,n-1\}}\|\Theta_{k}^{\xi,\gamma}-\vartheta\|\leq \delta$ and $\sup_{k\in\N}|\gamma_k|\leq 2\Gamma$ it holds that 
    \begin{equation}\llabel{eq16}
        \begin{split}
            \|(\gamma_n\varepsilon^{-1}I-\Lambda_n^{\xi,\gamma})(\nabla\smalll)(\Theta_{n-1}^{\xi,\gamma})\|&\leq \|\gamma_n\varepsilon^{-1}I-\Lambda_n^{\xi,\gamma}\|\|(\nabla\smalll)(\Theta_{n-1}^{\xi,\gamma})\|\leq \fC_1\delta \fC_2\|\Theta_{n-1}^{\xi,\gamma}-\vartheta\|\\
            &=\fC_3\delta\|\Theta_{n-1}^{\xi,\gamma}-\vartheta\|\leq \fC_4\delta\nnorm{\Theta_{n-1}^{\xi,\gamma}-\vartheta}.
        \end{split}
    \end{equation}}
    \argument{\lref{eq16};the fact that $\varrho>\rho$}{that there exists $\delta\in (0,\infty)$ such that for every $n\in\N$, $\xi\in\R^d$ and every $\gamma\colon\N\to\R$ with $\sup_{k\in\{0,1,\dots,n-1\}}\|\Theta_{k}^{\xi,\gamma}-\vartheta\|\leq \delta$ and $\sup_{k\in\N}|\gamma_k-\Gamma|\leq \delta$ it holds that 
    \begin{equation}\llabel{evidence 2}
        (\mu_1)^{-1}\|(\gamma_n\varepsilon^{-1}I-\Lambda_n^{\xi,\gamma})(\nabla\smalll)(\Theta_{n-1}^{\xi,\gamma})\|\leq \frac{\varrho-\rho}{8}\nnorm{\Theta_{n-1}^{\xi,\gamma}-\vartheta}.
    \end{equation}}
    \argument{\lref{def: bbL};the assumption that $\smalll\in C^2(\R^d,\R)$;}{that there exists $\delta\in (0,\infty)$ such that for all $\theta\in \R^d$ with $\|\theta-\vartheta\|\leq \delta$ it holds that
    \begin{equation}\llabel{eq17}
    \begin{split}
       & \Gamma\varepsilon^{-1}(\mu_1)^{-1}\|(\nabla\smalll)(\theta)-(\operatorname{Hess}\smalll)(\vartheta)(\theta-\vartheta)\|\\&=\Gamma\varepsilon^{-1}(\mu_1)^{-1}\|(\nabla\smalll)(\theta)-(\nabla\smalll)(\vartheta)-(\operatorname{Hess}\smalll)(\vartheta)(\theta-\vartheta)\|\leq \frac{\varrho-\rho}{4\mu_2}\|\theta-\vartheta\|.
        \end{split}
    \end{equation}}
    \argument{\lref{eq17};\lref{def: mu}}{that there exists $\delta\in (0,\infty)$ such that for every $n\in\N$, $\xi\in\R^d$ and every $\gamma \colon\N\to\R$ with $\|\Theta_{n-1}^{\xi,\gamma}-\vartheta\|\leq \delta$ it holds that
    \begin{equation}\llabel{evidence 3}
        \Gamma\varepsilon^{-1}(\mu_1)^{-1}\|(\nabla\smalll)(\Theta_{n-1}^{\xi,\gamma})-(\operatorname{Hess}\smalll)(\vartheta)(\Theta_{n-1}^{\xi,\gamma}-\vartheta)\|\leq \frac{\varrho-\rho}{4\mu_2}\|\Theta_{n-1}^{\xi,\gamma}-\vartheta\|\leq \frac{\varrho-\rho}{4}\nnorm{\Theta_{n-1}^{\xi,\gamma}-\vartheta}.
    \end{equation}}
    \argument{\lref{eq10};\lref{evidence 1};\lref{evidence 2};\lref{evidence 3};}{that there exist $\delta,\fC\in (0,\infty)$ which satisfy for every $n\in\N$, $\xi\in\R^d$ and every $\gamma\colon\N\to\R$ with $\sup_{k\in\{0,1,\dots,n-1\}}\|\Theta_{k}^{\xi,\gamma}-\vartheta\|\leq \delta$ and $\sup_{k\in\N}|\gamma_k-\Gamma|\leq \delta$ that
    \begin{equation}\llabel{eqmain}
         \begin{split}
             \nnorm{\Theta_n^{\xi,\gamma}-\vartheta}\leq \varrho\nnorm{\Theta_{n-1}^{\xi,\gamma}-\vartheta}.
             \end{split}
    \end{equation}}
    \startnewargseq
    In the following we prove that for all $n\in\N_0$, $\xi\in \R^d$ and all $\gamma\colon\N\to\R$ with $\sup_{k\in\N}|\gamma_k-\Gamma|< \delta$ and $\|\xi-\vartheta\|\leq \delta\mu_1(\mu_2)^{-1}$ it holds that
         \begin{equation}\llabel{need to prove}
             \nnorm{\Theta_n^{\xi,\gamma}-\vartheta}\leq(\mu_1)^{-1}\|\xi-\vartheta\| \varrho^n.
         \end{equation}
         We prove \lref{need to prove} by induction on $n\in\N_0$. For the base case $n=0$ note that \lref{def: Theta} and \lref{def: mu} show that for all $\xi\in \R^d$ and all $\gamma\colon\N\to\R$ with $\|\xi-\vartheta\|\leq \delta\mu_1(\mu_2)^{-1}$ it holds that
         \begin{equation}\llabel{ntp: base 1}
             \nnorm{\Theta_0^{\xi,\gamma}-\vartheta}=\nnorm{\xi-\vartheta}\leq (\mu_1)^{-1}\|\xi-\vartheta\|= (\mu_1)^{-1}\|\xi-\vartheta\|=(\mu_1)^{-1}\|\xi-\vartheta\|\varrho^{0}.
         \end{equation}
         This shows \lref{need to prove} in the base case $n=0$. For the induction step we assume that there exists $n\in\N$ which satisfies for all $k\in\{0,1,\dots,n-1\}$, $\xi\in \R^d$ and all $\gamma\colon\N\to\R$ with $\sup_{m\in\N}|\gamma_m-\Gamma|< \delta$ and $\|\xi-\vartheta\|\leq \delta\mu_1(\mu_2)^{-1}$ that
         \begin{equation}\llabel{ntp: assume}
              \nnorm{\Theta_{k}^{\xi,\gamma}-\vartheta}\leq(\mu_1)^{-1}\|\xi-\vartheta\|\varrho^{k}.
         \end{equation}
         \startnewargseq
         \argument{\lref{def: mu};\lref{ntp: assume};the fact that $0<\varrho<1$}{that for all $k\in\{1,2,\dots,n-1\}$, $\xi\in \R^d$ and all $\gamma\colon\N\to\R$ with $\sup_{m\in\N}|\gamma_m-\Gamma|<\delta$ and $\|\xi-\vartheta\|\leq \delta\mu_1(\mu_2)^{-1}$ it holds that
         \begin{equation}\llabel{ntp: eq1}
         \begin{split}
           \textstyle  \|\Theta_{k}^{\xi,\gamma}-\vartheta\|\leq \mu_2\nnorm{\Theta_{k}^{\xi,\gamma}-\vartheta}\leq \mu_2(\mu_1)^{-1}\varrho^{k}\|\xi-\vartheta\|\leq \delta.
             \end{split}
         \end{equation}}
         \argument{\lref{ntp: eq1};\lref{def: mu};\lref{eqmain};\lref{ntp: assume}}{for all $\xi\in \R^d$ and all $\gamma\colon\N\to\R$ with $\sup_{k\in\N}|\gamma_k-\Gamma|< \delta$  and $\|\xi-\vartheta\|\leq \delta\mu_1(\mu_2)^{-1}$ that
         \begin{equation}\llabel{ntp: eq2}
             \begin{split}
             \nnorm{\Theta_n^{\xi,\gamma}-\vartheta}\leq  \varrho\nnorm{\Theta_{n-1}^{\xi,\gamma}-\vartheta}
            \leq \varrho(\mu_1)^{-1} \|\xi-\vartheta\|\varrho^{n-1}=(\mu_1)^{-1}\|\xi-\vartheta\|\varrho^{n}.
             \end{split}
         \end{equation}}
         \argument{\lref{ntp: eq2};\lref{ntp: assume};induction}{\lref{need to prove}\dott}
         \startnewargseq
         \argument{\lref{def: mu};\lref{need to prove};}{that for all $n\in\N_0$, $\xi\in \R^d$ and all $\gamma\colon\N\to\R$ with $\sup_{k\in\N}|\gamma_k-\Gamma|<\delta$ and $\|\xi-\vartheta\|\leq \delta\mu_1(\mu_2)^{-1}$ it holds that
         \begin{equation}\llabel{eq19}
             \|\Theta_n^{\xi,\gamma}-\vartheta\|\leq \mu_2\nnorm{\Theta_n^{\xi,\gamma}-\vartheta}\leq \mu_2(\mu_1)^{-1}\|\xi-\vartheta\|\varrho^{n}=\mu_2(\mu_1)^{-1}\|\Theta_0^{\xi,\gamma}-\vartheta\|\varrho^{n}.
         \end{equation}}
         \argument{\lref{eq19};\lref{def: bfm};\lref{def: Theta};}{that there exist $\scrC,\eta \in (0,\infty)$ such that for every $\gamma\colon\N\to(\Gamma-\eta,\Gamma+\eta)$, every $\bbV=(\bbV^1,\dots,\bbV^{d})\colon \mathbb{N}_0\to\R^d$, and every $\Theta=(\Theta^1,\dots,\Theta^{d})\colon \mathbb{N}_0\to\R^d$ with the property that for all $n \in \mathbb{N}$, $i\in\{1,2,\dots,d\}$ it holds that 
         \begin{equation}
            \bbV_0=0,\qquad \bbV_n^i= \beta \bbV_{n-1}^i+(1-\beta)|(\nabla_{\theta_i}\smalll)(\Theta_{n-1})|^2,
            \end{equation}
            \begin{equation}
            \|\Theta_0-\vartheta\|< \eta,\qquad\text{and}\qquad  \Theta_n^i=\Theta_{n-1}^i-\gamma_n
          [\varepsilon+(\bbV_n^i)^{1/2}]^{-1}(\nabla_{\theta_i}\smalll)(\Theta_{n-1})
         \end{equation}
          we have for all $n\in\N_0$ that
$    \textstyle
    \|\Theta_n-\vartheta\|\leq \scrC\|\Theta_0-\vartheta\|\varrho^{n}$\dott}
\end{aproof}
\begin{tcolorbox}[colback=white!95!gray,
                  colframe=black,
                  boxrule=0.5pt,
                  sharp corners,
                  enhanced,
                  breakable,
                 ]
\begin{athm}{cor}{main theorem higher bound rmsprop}
Let $d\in\N$, $\chi\in (1,\infty)$, $\beta\in (0,1]$, $\varepsilon\in (0,\infty)$. Then for every $\convexcoe \in (0,\infty)$, $\vartheta \in \R^d$, $\smalll\in \allowbreak C^{2}(\R^d,\mathbb{R})$, $\varrho\in (\frac{ \chi  - 1 }  { \chi  + 1 },1)$ 
there exist $\fC,\Gamma, \delta \in (0,\infty)$ such that for every $\gamma\colon\N\to(\Gamma-\delta,\Gamma+\delta)$, every $\bbV=(\bbV^1,\dots,\bbV^d)\colon \mathbb{N}_0\to\R^d$, and every $\Theta=(\Theta^1,\dots,\Theta^d)\colon \mathbb{N}_0\to\R^d$ with the property that for all $n \in \mathbb{N}$, $i\in\{1,2,\dots,d\}$, $\theta\in\R^d$ it holds that
\begin{equation}
     \textcolor{magenta}{(\nabla\smalll)(\vartheta)=0},\qquad\textcolor{magenta}{\convexcoe\|\theta\|^2\leq \smalll'' (\vartheta)(\theta,\theta)\leq \convexcoe\chi \| \theta \|^2},
 \end{equation} 
 \begin{equation}
            \textcolor{magenta}{\bbV_0=0},\qquad \textcolor{magenta}{\bbV_n^i= \beta \bbV_{n-1}^i+(1-\beta)|(\nabla_{\theta_i}\smalll)(\Theta_{n-1})|^2},
            \end{equation}
            \begin{equation}
            \textcolor{magenta}{\|\Theta_0-\vartheta\|< \delta},\qquad\text{and}\qquad  \textcolor{magenta}{\Theta_n^i=\Theta_{n-1}^i-\gamma_n
          [\varepsilon+(\bbV_n^i)^{1/2}]^{-1}(\nabla_{\theta_i}\smalll)(\Theta_{n-1})}
         \end{equation}
          we have for all $n\in\N_0$ that
$    \textcolor{magenta}{ \|\Theta_n-\vartheta\|\leq \fC\|\Theta_0-\vartheta\|\varrho^n}$.
\end{athm}
\end{tcolorbox}
\begin{aproof}
    \argument{\cref{theo: higher bound rmsprop}}{that for every $\convexcoe \in (0,\infty)$, $\vartheta \in \R^d$, $\smalll\in \allowbreak C^{2}(\R^d,\mathbb{R})$, $\varrho\in (\frac{ \chi  - 1 }  { \chi  + 1 },1)$ there exist $\fC,\Gamma, \delta \in (0,\infty)$ such that for every $\gamma\colon\N\to(\Gamma-\delta,\Gamma+\delta)$, every $\bbV=(\bbV^1,\dots,\bbV^d)\colon \mathbb{N}_0\to\R^d$, and every $\Theta=(\Theta^1,\dots,\Theta^d)\colon \mathbb{N}_0\to\R^d$ with the property that for all $n \in \mathbb{N}$, $i\in\{1,2,\dots,d\}$, $\theta\in\R^d$ it holds that
    \begin{equation}
     (\nabla\smalll)(\vartheta)=0,\qquad\convexcoe\|\theta\|^2\leq \smalll'' (\vartheta)(\theta,\theta)\leq \convexcoe\chi \| \theta \|^2,
 \end{equation} 
 \begin{equation}
           \bbV_0=0,\qquad \bbV_n^i= \beta \bbV_{n-1}^i+(1-\beta)|(\nabla_{\theta_i}\smalll)(\Theta_{n-1})|^2,
            \end{equation}
            \begin{equation}
           \|\Theta_0-\vartheta\|< \delta,\qquad\text{and}\qquad \Theta_n^i=\Theta_{n-1}^i-\gamma_n
          [\varepsilon+(\bbV_n^i)^{1/2}]^{-1}(\nabla_{\theta_i}\smalll)(\Theta_{n-1})
         \end{equation}
          we have for all $n\in\N_0$ that
$    \|\Theta_n-\vartheta\|\leq \fC\|\Theta_0-\vartheta\|\varrho^n$\dott}
\end{aproof}
\subsection{Lower bounds on convergence rates of RMSprop}
\begin{tcolorbox}[colback=white!95!gray,
                  colframe=black,
                  boxrule=0.5pt,
                  sharp corners,
                  enhanced,
                  breakable,
                 ]
\begin{athm}{lemma}{lem: lower bound rmsprop}
      Let $d\in\N\backslash\{1\}$, $\convexcoe,\delta,\gamma,\varepsilon,\fC\in (0,\infty)$, $\chi\in (1,\infty)$, $\beta\in (0,1]$, $\varrho\in (0,\frac{\chi-1}{\chi+1})$ and let $\smalll=(\smalll(\theta))_{\theta=(\theta_1,\dots,\theta_d)\in\R^d}\colon \R^d\to\R$ satisfy for all $\theta=(\theta_1,\dots,\theta_d)\in \R^d$ that $\smalll(\theta)=\convexcoe\chi|\theta_1|^2+\sum_{i=2}^d \convexcoe|\theta_i|^2$. Then there exist $\bbV=(\bbV^1,\dots,\bbV^d)\colon \mathbb{N}_0\to\R^d$ and $\Theta=(\Theta^1,\dots,\Theta^d)\colon \mathbb{N}_0\to\R^d$ such that
\begin{enumerate}[label=(\roman*)]
    \item \llabel{item 1} it holds for all $n \in \mathbb{N}$, $i\in\{1,2,\dots,d\}$ that
          \begin{equation}
            \textcolor{magenta}{\bbV_0=0},\qquad \textcolor{magenta}{\bbV_n^i= \beta \bbV_{n-1}^i+(1-\beta)|(\nabla_{\theta_i}\smalll)(\Theta_{n-1})|^2},
            \end{equation}
            \begin{equation}
            \textcolor{magenta}{\|\Theta_0\|<\delta},\qquad\text{and}\qquad  \textcolor{magenta}{\Theta_n^i=\Theta_{n-1}^i-\gamma
          [\varepsilon+(\bbV_n^i)^{1/2}]^{-1}(\nabla_{\theta_i}\smalll)(\Theta_{n-1})}
         \end{equation}
         and
         \item \llabel{item 2} there exists $n\in\N$ such that $\textcolor{magenta}{ \textstyle\|\Theta_n\|>\fC\|\Theta_0\|\varrho^n}$.
\end{enumerate} 
\end{athm}
\end{tcolorbox}
\begin{aproof}
    Throughout this proof for every $i\in\{1,2,\dots,d\}$ let $h_i\in\R$ satisfy 
    \begin{equation}\llabel{def: h}
        h_i=1+(\chi-1)\mathbbm 1_{\{1\}}(i).
    \end{equation}
    \argument{the fact that $\chi>1$;}{that
    \begin{equation}\llabel{eq1}
        \textstyle \chi-1=|\chi(1-\frac{2r\gamma}{\varepsilon})-(1-\frac{2r\gamma \chi}{\varepsilon})|\leq |\chi(1-\frac{2r\gamma}{\varepsilon})|+|1-\frac{2r\gamma \chi}{\varepsilon}|\leq (\chi+1)\max\bigl\{|1-\frac{2r\gamma}{\varepsilon}|,|1-\frac{2r\gamma \chi}{\varepsilon}|\bigr\}.
    \end{equation}}
    \argument{\lref{eq1};}{that
    \begin{equation}\llabel{eq2}
        \textstyle\max\bigl\{|1-\frac{2r\gamma}{\varepsilon}|,|1-\frac{2r\gamma \chi}{\varepsilon}|\bigr\}\geq \frac{\chi-1}{\chi+1}.
    \end{equation}}
    \argument{\lref{def: h};}{that \llabel{eq3}
        $\{h_1,h_2,\dots,h_d\}=\{1,\chi\}$\dott}
    \argument{\lref{eq3};\lref{eq2};}{that
    \begin{equation}\llabel{eq4}
      \textstyle  \max_{i\in\{1,2,\dots,d\}}|1-\frac{2r\gamma h_i}{\varepsilon}|\geq \frac{\chi-1}{\chi+1}.
    \end{equation}}
    \argument{\lref{eq4};}{that there exists $i\in\{1,2,\dots,d\}$ which satisfies 
    \begin{equation}\llabel{def: i}
        \biggl|1-\frac{2r\gamma h_i}{\varepsilon}\biggr|\geq \frac{\chi-1}{\chi+1}.
    \end{equation}}
    \startnewargseq
    \argument{\lref{def: i};the fact that $\varrho<\frac{\chi-1}{\chi+1}$}{that
    \begin{equation}\llabel{eq5}
      \textstyle  \bigl|1-\frac{2r\gamma h_i}{\varepsilon}\bigr|>\varrho.
    \end{equation}}
    \argument{\lref{eq5};}{that there exists $N\in\N$ which satisfies
    \begin{equation}\llabel{def: N}
        \textstyle  \bigl|1-\frac{2r\gamma h_i}{\varepsilon}\bigr|^N>\fC\varrho^N.
    \end{equation}}
\startnewargseq
In the following let $\bfe=(\bfe_1,\dots,\bfe_d)\in \R^d$ satisfy for all $j\in\{1,2,\dots,d\}$ that
    \begin{equation}\llabel{def: bfe}
        \bfe_{j}=\mathbbm 1_{\{i\}}(j)
    \end{equation}
     and for every $t\in\R$ let $\bbV^{t}=(\bbV^{t,1},\dots,\bbV^{t,d})\colon\N_0\to\R^d$ and $\Theta^{t}=(\Theta^{t,1},\dots,\Theta^{t,d})\colon\N_0\to\R^d$ satisfy for all $n\in\N$, $j\in\{1,2,\dots,d\}$ that
     \begin{equation}\llabel{def: bfm}
      \bbV^{t}_0=0,\qquad \bbV^{t,j}_n= \beta \bbV^{t,j}_{n-1}+(1-\beta)|2rh_j\Theta^{t,j}_{n-1}|^2,
            \end{equation}
            \begin{equation}\llabel{def: Theta}
        \Theta^{t}_0=t\bfe,\qquad\text{and}\qquad \Theta^{t,j}_n=\Theta^{t,j}_{n-1}-\gamma
          [\varepsilon+(\bbV_n^{t,j})^{1/2}]^{-1}2rh_j\Theta_{n-1}^{t,j}
         \end{equation}
         and let $\bfu^t\colon \N_0\to\R$, $\bfv^t\colon\N_0\to\R$ satisfy for all $n\in\N$ that
         \begin{equation}\llabel{def: bfu}
             \bfu_0^t=t,\qquad \bfu_n^t=\bfu_{n-1}^t-\gamma[\varepsilon+(\bfv_n^{t})^{1/2}]^{-1}2rh_i\bfu_{n-1}^{t},
             \end{equation}
             \begin{equation}\llabel{def: bfv}
        \bfv_0^t=0,\qqandqq\bfv_n^t=\beta\bfv_{n-1}^t+(1-\beta)|2rh_i \bfu_{n-1}^t|^2.
         \end{equation}
         In the following we prove that for all $n\in\N_0$ it holds that
         \begin{equation}\llabel{need to prove}
            \textstyle\limsup_{t\to 0}\bigl|\frac{\bfu_n^t}{t}-(1-\frac{2r\gamma h_i}{\varepsilon})^n\bigr|=0\qqandqq \limsup_{t\to 0}|\bfv_n^t|=0.
         \end{equation}
         We prove \lref{need to prove} by induction on $n\in\N_0$. For the base case $n=0$ note that  \lref{def: bfu} and \lref{def: bfv} show that for all $t\in \R\backslash\{0\}$ it holds that
         \begin{equation}\llabel{base: eq1}
             \textstyle\bigl|\frac{\bfu_0^t}{t}-(1-\frac{2r\gamma h_i}{\varepsilon})^0\bigr|=\bigl|\frac{t}{t}-1\bigr|=0\qqandqq \bfv_0^t=0.
         \end{equation}
         This shows \lref{need to prove} in the base case $n=0$. For the induction step we assume that there exists $n\in\N$ which satisfies that
         \begin{equation}\llabel{ntp: assume}
             \textstyle\limsup_{t\to 0}\bigl|\frac{\bfu_{n-1}^t}{t}-(1-\frac{2r\gamma h_i}{\varepsilon})^{n-1}\bigr|=0\qqandqq \limsup_{t\to 0}|\bfv_{n-1}^t|=0. 
         \end{equation}
         \argument{\lref{ntp: assume};}{that
         \begin{equation}\llabel{ntp: eq1}
             \textstyle\limsup_{t\to 0}|\bfu_{n-1}^t|=0.
         \end{equation}}
         \argument{\lref{ntp: eq1};\lref{def: bfv};\lref{ntp: assume}}{that
         \begin{equation}\llabel{ntp: eq2}
           \textstyle  \limsup_{t\to 0}|\bfv_{n}^t|\leq \bigl[\limsup_{t\to 0}|\beta\bfv_{n-1}^t|\bigr]+(1-\beta)\bigl[\limsup_{t\to 0}|2rh_i \bfu_{n-1}^t|^2\bigr] =0+0=0.
         \end{equation}}
         \argument{\lref{ntp: eq2};\lref{def: bfu};\lref{ntp: assume}}{that
         \begin{equation}\llabel{ntp: eq3}
         \begin{split}
        \textstyle     \lim_{t\to 0} \frac{\bfu_n^t}{t}&\textstyle=  \lim_{t\to 0} \frac{\bfu_{n-1}^t-\gamma[\varepsilon+(\bfv_n^{t})^{1/2}]^{-1}2rh_i\bfu_{n-1}^{t}}{t}= \lim_{t\to 0} \frac{\bfu_{n-1}^t(1-\gamma[\varepsilon+(\bfv_n^{t})^{1/2}]^{-1}2rh_i)}{t}\\&\textstyle=\bigl[\lim_{t\to 0} \frac{\bfu_{n-1}^t}{t}\bigr]\bigl(1-\bigl[\lim_{t\to 0} [\varepsilon+(\bfv_n^{t})^{1/2}]^{-1}2r\gamma h_i\bigr]\bigr)\\
        &=\textstyle\bigl[\lim_{t\to 0} \frac{\bfu_{n-1}^t}{t}\bigr]\bigl(1-\bigl[\lim_{t\to 0} \frac{\varepsilon}{\varepsilon+(\bfv_n^{t})^{1/2}}\bigr]\bigl[\lim_{t\to 0} \frac{2r\gamma h_i}{\varepsilon}\bigr]\bigr)\\
        &=\textstyle(1-\frac{2r\gamma h_i}{\varepsilon})^{n-1}(1-\frac{2r\gamma h_i}{\varepsilon})=(1-\frac{2r\gamma h_i}{\varepsilon})^{n}.
        \end{split}
         \end{equation}}
         \argument{\lref{ntp: eq3};\lref{ntp: eq2}}{that
         \begin{equation}\llabel{ntp: eq4}
              \textstyle\limsup_{t\to 0}\bigl|\frac{\bfu_{n}^t}{t}-(1-\frac{2r\gamma h_i}{\varepsilon})^{n}\bigr|=0\qqandqq \limsup_{t\to 0}|\bfv_{n}^t|=0. 
         \end{equation}}
         \argument{\lref{ntp: eq4};\lref{ntp: assume};induction}{\lref{need to prove}\dott}
         \startnewargseq
         In the following we prove that for all $n\in\N_0$, $t\in \R$, $j\in\{1,2,\dots,d\}$ it holds that
         \begin{equation}\llabel{need to prove 2}
           \Theta_n^{t,j}=\bfu_n^t\mathbbm 1_{\{i\}}(j)\qqandqq \bbV_n^{t,j}=\bfv_n^t\mathbbm 1_{\{i\}}(j).
         \end{equation}
         We prove \lref{need to prove 2} by induction on $n\in\N_0$. For the base case $n=0$ note that \lref{def: bfe}, \lref{def: bfm}, \lref{def: Theta}, \lref{def: bfu}, and \lref{def: bfv} show that for all $t\in \R$, $j\in\{1,2,\dots,d\}$ it holds that
         \begin{equation}
             \Theta_0^{t,j}=t\bfe_j=\bfu_0^t\mathbbm 1_{\{i\}}(j)\qqandqq \bbV_0^{t,j}=0=\bfv_0^t\mathbbm 1_{\{i\}}(j).
         \end{equation}
         This proves \lref{need to prove 2} in the base case $n=0$. For the induction step we assume that there exists $n\in\N$ which satisfies for all $t\in\R$, $j\in\{1,2,\dots,d\}$ that
         \begin{equation}\llabel{ntp2: assume}
              \Theta_{n-1}^{t,j}=\bfu_{n-1}^t\mathbbm 1_{\{i\}}(j)\qqandqq \bbV_{n-1}^{t,j}=\bfv_{n-1}^t\mathbbm 1_{\{i\}}(j).
         \end{equation}
         \argument{\lref{def: bfm};\lref{def: bfv};\lref{ntp2: assume}}{that for all  $t\in\R$, $j\in\{1,2,\dots,d\}$ it holds that
         \begin{equation}\llabel{ntp2: eq1}
         \begin{split}
             \bbV_n^{t,j}&=\beta \bbV^{t,j}_{n-1}+(1-\beta)|2rh_j\Theta^{t,j}_{n-1}|^2=\beta\bfv_{n-1}^t\mathbbm 1_{\{i\}}(j)+(1-\beta)|2rh_j\bfu_{n-1}^t\mathbbm 1_{\{i\}}(j)|^2\\
             &=\bigl[\beta\bfv_{n-1}^t+(1-\beta)|2rh_j\bfu_{n-1}^t|^2\bigr]\mathbbm 1_{\{i\}}(j)=\bfv_{n}^t\mathbbm 1_{\{i\}}(j).
             \end{split}
         \end{equation}}
         \argument{\lref{ntp2: eq1};\lref{def: Theta};\lref{def: bfu};\lref{ntp2: assume}}{for all $t\in\R$, $j\in\{1,2,\dots,d\}\backslash\{i\}$ that
         \begin{equation}\llabel{ntp2: eq2}
         \begin{split}
             \Theta_n^{t,i}=\Theta^{t,i}_{n-1}-\gamma
          [\varepsilon+(\bbV_n^{t,i})^{1/2}]^{-1}2rh_i\Theta_{n-1}^{t,i}=\bfu_{n-1}^t-\gamma [\varepsilon+(\bfv_n^{t})^{1/2}]^{-1}2rh_i\bfu_{n-1}^t=\bfu_n^t
          \end{split}
         \end{equation}
         and 
         \begin{equation}\llabel{ntp2: eq3}
             \Theta_n^{t,j}=\Theta^{t,j}_{n-1}-\gamma
          [\varepsilon+(\bbV_n^{t,j})^{1/2}]^{-1}2rh_j\Theta_{n-1}^{t,j}=0-0\gamma
          [\varepsilon+(\bbV_n^{t,j})^{1/2}]^{-1}2rh_j=0.
         \end{equation}}
         \argument{\lref{ntp2: eq3};\lref{ntp2: eq2};}{for all $t\in\R$, $j\in\{1,2,\dots,d\}$ that
         \begin{equation}\llabel{ntp2: eq4}
               \Theta_{n}^{t,j}=\bfu_{n}^t\mathbbm 1_{\{i\}}(j)\qqandqq \bbV_{n}^{t,j}=\bfv_{n}^t\mathbbm 1_{\{i\}}(j).
         \end{equation}}
         \argument{\lref{ntp2: eq4};\lref{ntp2: assume};induction}{\lref{need to prove 2}\dott}
         \startnewargseq
         \argument{\lref{def: N};\lref{need to prove};\lref{need to prove 2};}{that
         \begin{equation}\llabel{eq10}
             \textstyle \lim_{t\to 0}\frac{\|\Theta_N^t\|}{|t|}=\lim_{t\to 0}\frac{\|\bfu_N^t\|}{|t|}=\bigl|(1-\frac{2r\gamma h_i}{\varepsilon})^N\bigr|> \fC\varrho^N.
         \end{equation}}
         \argument{\lref{eq10};\lref{def: Theta};}{that there exists $t\in(-\delta,\delta)$ such that
         \begin{equation}\llabel{eq11}
          \|\Theta_N^t\|>\fC\varrho^N |t|=\fC\|\Theta_0^t\|\varrho^N.
         \end{equation}}
         \argument{\lref{def: h};the assumption that for all $\theta=(\theta_1,\dots,\theta_d)\in \R^d$ it holds that $\smalll(\theta)=\convexcoe\chi|\theta_1|^2+\sum_{i=2}^d \convexcoe|\theta_i|^2$}{that for all $j\in\{1,2,\dots,d\}$, $\theta=(\theta_1,\dots,\theta_d)\in\R^d$ it holds that
         \begin{equation}\llabel{eq12}
           (\nabla_{\theta_j}\smalll)(\theta)=2r h_j\theta_j.  
         \end{equation}}
         \argument{\lref{eq12};\lref{def: bfm};\lref{def: Theta};\lref{eq11}}{\lref{item 1,item 2}\dott}
\end{aproof}
\begin{tcolorbox}[colback=white!95!gray,
                  colframe=black,
                  boxrule=0.5pt,
                  sharp corners,
                  enhanced,
                  breakable,
                 ]
\begin{athm}{cor}{main theorem lower bound rmsprop}
     Let $d\in\N\backslash\{1\}$, $\chi\in (1,\infty)$, $\beta\in (0,1]$, $\varrho\in (0,\frac{\chi - 1 }  { \chi + 1 })$, $\varepsilon\in (0,\infty)$. Then there exist $\convexcoe\in (0,\infty)$, $\vartheta\in\R^d$ and $\smalll\in C^2( \R^d,\R)$ with $(\nabla \smalll)( \vartheta ) = 0$ and $
			\forall\, \theta\in\R^d\colon \convexcoe\|\theta\|^2\leq  \smalll''(\vartheta)(\theta,\theta)\leq \convexcoe\chi \| \theta \|^2
			$ such that for all $\delta,\Gamma,\fC \in (0,\infty)$ there exist $\gamma\colon\N\to (\Gamma-\delta,\Gamma+\delta)$, $\bbV=(\bbV^1,\dots,\bbV^d)\colon \mathbb{N}_0\to\R^d$, and $\Theta=(\Theta^1,\dots,\Theta^d)\colon \mathbb{N}_0\to\R^d$  such that
\begin{enumerate}[label=(\roman*)]
    \item \llabel{item 1} it holds for all $n \in \mathbb{N}$, $i\in\{1,2,\dots,d\}$  that
        \begin{equation}
            \textcolor{magenta}{\bbV_0=0},\qquad \textcolor{magenta}{\bbV_n^i= \beta \bbV_{n-1}^i+(1-\beta)|(\nabla_{\theta_i}\smalll)(\Theta_{n-1})|^2},
            \end{equation}
            \begin{equation}
            \textcolor{magenta}{\|\Theta_0-\vartheta\|< \delta},\qquad\text{and}\qquad  \textcolor{magenta}{\Theta_n^i=\Theta_{n-1}^i-\gamma_n
          [\varepsilon+(\bbV_n^i)^{1/2}]^{-1}(\nabla_{\theta_i}\smalll)(\Theta_{n-1})}
         \end{equation}
         and
         \item \llabel{item 2} there exists $n\in\N$ such that $\textcolor{magenta}{ \textstyle
   \textstyle\|\Theta_n-\vartheta\|>\fC\|\Theta_0-\vartheta\|\varrho^n}$.
\end{enumerate} 
\end{athm}
\end{tcolorbox}
    \begin{aproof}
    Throughout this proof let $\convexcoe\in\R$, $\vartheta\in\R^d$  satisfy $\convexcoe=2$ and $\vartheta=0$ and let $\smalll\colon\R^d\to\R$ satisfy for all $\theta=(\theta_1,\dots,\theta_d)\in\R^d$ that $\smalll(\theta)=\chi|\theta_1|^2+\sum_{i=2}^d|\theta_i|^2$.
\argument{\cref{lem: lower bound rmsprop};}{\lref{item 1,item 2}\dott}
    \end{aproof}
    \subsection{Sharp convergence rates of GD and RMSprop}\label{subsec: sharp rmsprop}
      \begin{samepage}
 \begin{tcolorbox}[colback=white!95!gray,
                  colframe=black,
                  boxrule=0.5pt,
                  sharp corners,
                  enhanced,
                  unbreakable,
                 ]
                 \begin{athm}{theorem}{theo: sharp rate rmsprop}
                 Let $d\in \mathbb{N}\backslash\{1\}$, $\chi\in (1,\infty)$, $\beta\in(0,1]$, $\varrho=\frac{ \chi  - 1 }  {  \chi  + 1 }$, $\varepsilon\in (0,\infty)$ and for every $n\in\N$ let $\varphi_n=(\varphi_n^1,\dots,\varphi_n^d)\colon (\R^d)^n\to\R^d$ satisfy for all $\theta=((\theta_k^i)_{i\in\{1,2,\dots,d\}})_{k\in\{1,2,\dots,n\}}\in (\R^d)^n$, $i\in\{1,2,\dots,d\}$ that
\begin{equation}\llabel{def: Phi}
\textcolor{magenta}{\textstyle\varphi_n^i(\theta)=\bigl[\varepsilon+\bigl((1-\beta)\sum_{k=1}^{n}\beta^{n-k}|\theta_k^i|^2\bigr)^{1/2}\bigr]^{-1}\theta_n^i}.
\end{equation}
Then it holds for all $\rho \in (0,\infty) \backslash \{ \varrho \}$ that the following three statements are equivalent:
             \begin{enumerate}[label=(\roman*)]
                 \item \llabel{item 1} It holds that $\rho > 
                 \varrho$.
                 \item \llabel{item 2} For all $\convexcoe \in (0,\infty)$, $\cL \in C^2( \R^d, \R )$, $\vartheta\in(\nabla\cL)^{-1}(\{0\})$ with $
			\forall\, \theta\in\R^d\colon \convexcoe\|\theta\|^2\leq \cL'' (\vartheta)(\theta,\theta)\leq \convexcoe\chi \| \theta \|^2
			$ there exist $\delta,\Gamma,\fC \in (0,\infty)$ such that for every $\gamma=(\gamma_n)_{n\in\N}\colon\N\to (\Gamma-\delta,\Gamma+\delta)$, every $\bbV=(\bbV^1,\dots,\bbV^d)\colon \mathbb{N}_0\to\R^d$, and every $\Theta=(\Theta^1,\dots,\Theta^d)\colon \mathbb{N}_0\to\R^d$ with the property that 
    for all $n \in \mathbb{N}$, $i\in\{1,2,\dots,d\}$ it holds that
          \begin{equation}
            \textcolor{magenta}{\bbV_0=0},\qquad \textcolor{magenta}{\bbV_n^i= \beta \bbV_{n-1}^i+(1-\beta)|(\nabla_{\theta_i}\cL)(\Theta_{n-1})|^2},
            \end{equation}
            \begin{equation}
            \textcolor{magenta}{\|\Theta_0-\vartheta\|< \delta},\qquad\text{and}\qquad  \textcolor{magenta}{\Theta_n^i=\Theta_{n-1}^i-\gamma_n
          [\varepsilon+(\bbV_n^i)^{1/2}]^{-1}(\nabla_{\theta_i}\cL)(\Theta_{n-1})}
         \end{equation}
         we have for all $n\in\N_0$ that $\textcolor{magenta}{ \textstyle
   \textstyle\|\Theta_n-\vartheta\|\leq\fC\rho^n\|\Theta_0-\vartheta\|}$.
                 \item \llabel{item 3} For all $\convexcoe \in (0,\infty)$, $\cL \in C^2(\R^d, \R )$, $\vartheta\in(\nabla\cL)^{-1}(\{0\})$ with $
			\forall\, \theta\in\R^d\colon \convexcoe\|\theta\|^2\leq \cL'' (\vartheta)(\theta,\theta)\leq \convexcoe\chi \| \theta \|^2
			$ there exist $\delta,\Gamma, \fC \in (0,\infty) $ such that 
                  for every $\gamma=(\gamma_n)_{n\in\N}\colon\N\to (\Gamma-\delta,\Gamma+\delta)$ and every $ \Theta \colon \N_0\to\R^d $
			with $\| \Theta_0 - \vartheta \| < \delta$ and
			\begin{equation}
				 \textcolor{magenta}{\forall \, n \in \N \colon
				\Theta_n
				=
				\Theta_{ n - 1 }
				-
				\gamma_n \varphi_n\bigl(
				( \nabla \cL )( \Theta_0 ) ,
				( \nabla \cL )( \Theta_1 ) ,
				\dots ,
				( \nabla \cL )( \Theta_{ n - 1 } )
				\bigr)}
			\end{equation}
			we have for all $n\in\N_0$ that
			$
			\textcolor{magenta}{\| \Theta_n - \vartheta \|
			\leq \fC\rho^{  n } \|\Theta_0-\vartheta\|}
			$.
             \end{enumerate}
                 \end{athm}
                 \end{tcolorbox}
                 \end{samepage}
    \begin{aproof}
        Throughout this proof for every $\Theta\colon\N_0\to\R^d$ and $\cL=(\cL(\theta))_{\theta=(\theta_1,\dots,\theta_d)\in\R^d}\in C^1(\R^d,\allowbreak\R)$ let $\bbV(\Theta,\cL)=(\bbV_n^i(\Theta,\cL))_{(i,n)\in\{1,2,\dots,d\}\times\N_0}\colon\N_0\to\R^d$ satisfy for all $n\in\N_0$, $i\in\{1,2,\dots,d\}$ that
        \begin{equation}\llabel{def: bfm}
            \bbV_n^i(\Theta,\cL)=(1-\beta)\sum_{k=1}^n\beta^{n-k}|(\nabla_{\theta_i}\cL)(\Theta_{k-1})|^2.
        \end{equation}
        \argument{\cref{main theorem higher bound rmsprop}}{that (\lref{item 1}$\rightarrow$\lref{item 2})\dott}
        \startnewargseq
         \argument{\cref{main theorem lower bound rmsprop}}{that (\lref{item 2}$\rightarrow$\lref{item 1})\dott}
        \startnewargseq
        \argument{\lref{def: bfm};}{that for all $n\in\N$, $\cL=(\cL(\theta))_{\theta=(\theta_1,\dots,\theta_d)\in\R^d}\in C^1(\R^d,\R)$, $i\in\{1,2,\dots,d\}$ and all $\Theta\colon\N_0\to\R^d$ it holds that
        \begin{equation}\llabel{eq1}
            \bbV^i_{n}(\Theta,\cL)=\beta\bbV^i_{n-1}(\Theta,\cL)+(1-\beta)|(\nabla_{\theta_i}\cL)(\Theta_{n-1})|^2.
        \end{equation}}
       \argument{\lref{def: Phi};\lref{def: bfm}}{that for all $n\in\N$, $\cL=(\cL(\theta))_{\theta=(\theta_1,\dots,\theta_d)\in\R^d}\in C^1(\R^d,\R)$, $i\in\{1,2,\dots,d\}$ and all $\Theta\colon\N_0\to\R^d$ it holds that
       \begin{equation}\llabel{eq2}
       \begin{split}
      &  \varphi_n^i\bigl(
				( \nabla \cL )( \Theta_0 ) ,
				( \nabla \cL )( \Theta_1 ) ,
				\dots ,
				( \nabla \cL )( \Theta_{ n - 1 } )
				\bigr)
                \\&= \textstyle \bigl[\varepsilon+\bigl((1-\beta)\sum_{k=1}^{n}\beta^{n-k}|(\nabla_{\theta_i}\cL)(\Theta_{k-1})|^2\bigr)^{1/2}\bigr]^{-1}(\nabla_{\theta_i}\cL)(\Theta_{n-1})\\
                &\textstyle=\bigl[\varepsilon+(\bbV_n^i(\Theta,\cL))^{1/2}\bigr]^{-1}(\nabla_{\theta_i}\cL)(\Theta_{n-1}).
                \end{split}
       \end{equation}}
       \argument{\lref{eq2};\lref{eq1}}{that (\lref{item 2}$\rightarrow$\lref{item 3})\dott}
       \startnewargseq
       \argument{\unskip, \eg, \cite[Lemma 2.7]{DeRoArAd2025}}{that for every $\cL=(\cL(\theta))_{\theta=(\theta_1,\dots,\theta_d)\in\R^d}\in C^2(\R^d,\R)$, every $\bbV=(\bbV^1,\dots,\bbV^d)\colon\N_0\to\R^d$, and every $\Theta\colon\N_0\to\R^d$ with the property that for all $n\in\N$, $i\in\{1,2,\dots,d\}$ it holds that $ \bbV_0=0$ and $\bbV_n^i= \beta\bbV_{n-1}^i+(1-\beta)|(\nabla_{\theta_i}\cL)(\Theta_{n-1})|^2$ we have for all $n\in\N_0$, $i\in\{1,2,\dots,d\}$ that
       \begin{equation}\llabel{eq4}
          \textstyle \bbV_n^i=(1-\beta)\sum_{k=1}^n\beta^{n-k}|(\nabla_{\theta_i}\cL)(\Theta_{k-1})|^2.
       \end{equation}}
       \argument{\lref{eq4};\lref{def: Phi}}{that for every $\cL=(\cL(\theta))_{\theta=(\theta_1,\dots,\theta_d)\in\R^d}\in C^2(\R^d,\R)$, every $\bbV=(\bbV^1,\dots,\bbV^d)\colon\N_0\to\R^d$, and every $\Theta\colon\N_0\to\R^d$ with the property that for all $n\in\N$, $i\in\{1,2,\dots,d\}$ it holds that $ \bbV_0=0$ and $\bbV_n^i= \beta \bbV_{n-1}^i+(1-\beta)|(\nabla_{\theta_i}\cL)(\Theta_{n-1})|^2$ we have for all $n\in\N$, $i\in\{1,2,\dots,d\}$ that
       \begin{equation}\llabel{eq5}
       \begin{split}
           &\varphi^i_n\bigl(
				( \nabla \cL )( \Theta_0 ) ,
				( \nabla \cL )( \Theta_1 ) ,
				\dots ,
				( \nabla \cL )( \Theta_{ n - 1 } )
				\bigr)\\
               &= \textstyle \bigl[\varepsilon+\bigl((1-\beta)\sum_{k=1}^{n}\beta^{n-k}|(\nabla_{\theta_i}\cL)(\Theta_{k-1})|^2\bigr)^{1/2}\bigr]^{-1}(\nabla_{\theta_i}\cL)(\Theta_{n-1})\\
                &\textstyle=\bigl[\varepsilon+(\bbV_n^i)^{1/2}\bigr]^{-1}(\nabla_{\theta_i}\cL)(\Theta_{n-1}).
                \end{split}
       \end{equation}}
       \argument{\lref{eq5};}{that (\lref{item 3}$\rightarrow$\lref{item 2})\dott}
    \end{aproof}
 \begin{tcolorbox}[colback=white!95!gray,
                  colframe=black,
                  boxrule=0.5pt,
                  sharp corners,
                  enhanced,
                  breakable,
                 ]
                 \begin{athm}{cor}{cor: sharp rate rmsprop}
                 Let $d\in \mathbb{N}\backslash\{1\}$, $\chi\in (1,\infty)$, $\beta\in (0,1]$, $\varrho=\frac{ \chi  - 1 }  {  \chi  + 1} $, $\varepsilon\in (0,\infty)$ and for every $n\in\N$ let $\varphi_n=(\varphi_n^1,\dots,\varphi_n^d)\colon (\R^d)^n\to\R^d$ satisfy for all $\theta=((\theta_k^i)_{i\in\{1,2,\dots,d\}})_{k\in\{1,2,\dots,n\}}\in (\R^d)^n$, $i\in\{1,2,\dots,d\}$ that
\begin{equation}\llabel{def: Phi}
\textcolor{magenta}{\textstyle\varphi_n^i(\theta)=\bigl[\varepsilon+\bigl((1-\beta)\sum_{k=1}^{n}\beta^{n-k}|\theta_k^i|^2\bigr)^{1/2}\bigr]^{-1}\theta_n^i}.
\end{equation}
Then $(\varphi_n)_{n\in\N}$ converges with $\chi$-rate $\varrho$ (cf.\ \cref{definition: rate of convergence}).
                 \end{athm}
                 \end{tcolorbox}
                 \begin{aproof}
                     \argument{\cref{theo: sharp rate rmsprop};}{that $(\varphi_n)_{n\in\N}$ converges with $\chi$-rate $\varrho$\dott}
                 \end{aproof}
    \section{Optimal convergence rates of the Adam optimizer}
\label{sec: adam}
In this section we establish in \cref{cor: sharp rate adam} below that for every $d\in\N\backslash\{1\}$, $\beta\in(0,1)$, $\varepsilon\in(0,\infty)$ the $\alpha$-$\beta$-$\varepsilon$-\Adam\ optimizer on $\R^d$ with $\alpha=(\sqrt{\chi}+1)^{-2}(\sqrt{\chi}-1)^2$ converges with $\chi$-rate $(\sqrt{\chi}+1)^{-1}(\sqrt{\chi}-1)$.

    \subsection{Adam optimizer}
    
In the next \cref{definition: Adam} we recall the concept of the \Adam\ optimizer (see \cite{KingmaBa2024_Adam}) using the functions $\Phi_n\colon(\R^d)^n\to\R^d$, $n\in\N$, from \cref{definition: rate of convergence} above; cf., \eg, \cite[Definition~4.1]{DeArAdsharp2025}.
    \begin{tcolorbox}[colback=white!95!gray,
                  colframe=black,
                  boxrule=0.5pt,
                  sharp corners,
                  enhanced,
                  breakable,
                 ]
    \begin{definition}[\Adam]
    \label{definition: Adam}
    Let $ d \in \N $,
    $ \alpha, \beta \in [0,1) $,
    $ \varepsilon \in (0,\infty) $
    and let
    $
      \Phi_n = ( \Phi_n^1, \dots, \Phi_n^{ d } )
      \colon
      \allowbreak
      ( \R^{ d } )^n
      \to \R^{ d }
    $,
    $ n \in \N $, be functions.
    Then we say that
    $ ( \Phi_n )_{ n \in \N } $ is the
    $ \alpha $-$ \beta $-$\varepsilon$-\Adam\ optimizer
    on $ \R^{ d } $ (we say that
    $ ( \Phi_n )_{ n \in \N } $ is
    the $ \alpha $-$ \beta $-$\varepsilon$-\Adam\ optimizer)
    if and only if it holds for all
    $ n \in \N $,
    $
      g =
      (
        ( g_{ i, j } )_{ j \in \{ 1, 2, \dots, d \} }
      )_{ i \in \{ 1, 2, \dots, n \} } \in ( \R^{ d } )^n
    $,
    $ j \in \{ 1, 2, \dots, d \} $
    that
    \begin{equation}
    \textstyle
	  \Phi_n^j( g )
	  =
	  \frac{ ( 1 - \alpha ) }{
	    ( 1 - \alpha^n )
	  }
	  \bigl[
	  \sum\nolimits_{ i = 1 }^n
	  \alpha^{ n - i }
	  g_{ i, j }
	  \bigr]
	  \bigl[
	    \varepsilon
	    +
	    (
	      ( 1 - \beta^n )^{ - 1 }
	      ( 1 - \beta )
	      \sum_{ i = 1 }^n
	      \beta^{ n - i }
	      | g_{ i, j } |^2
	    )^{ \nicefrac{ 1 }{ 2 } }
	  \bigr]^{ - 1 }
	  .
    \end{equation}
	\end{definition}
    \end{tcolorbox}

\subsection{Upper bounds on convergence rates of Adam}
\begin{samepage}
 \begin{tcolorbox}[colback=white!95!gray,
                  colframe=black,
                  boxrule=0.5pt,
                  sharp corners,
                  enhanced,
                  breakable,
                 ]
\begin{athm}{theorem}{theo: higher bound adam}
      Let $d\in \mathbb{N}$,  $\chi\in (1,\infty)$, $\convexcoe,\varepsilon\in (0,\infty)$, $\varrho\in (\frac{ \sqrt{\chi}  - 1 }  { \sqrt{\chi}  + 1 },1)$, $\Gamma=(\convexcoe\sqrt{\chi})^{-1}\varepsilon$, $\alpha=(\frac{\sqrt{\chi}-1}{\sqrt{\chi}+1})^2$, $\beta\in (0,1)$, $\vartheta\in\R^d$ and let $\smalll=(\smalll(\theta))_{\theta=(\theta_1,\dots,\theta_d)\in\R^d}\in \allowbreak C^{2}(\R^d,\mathbb{R})$ satisfy for all $\theta \in \R^d$ that \begin{equation}\llabel{def: bbL}
     \textcolor{magenta}{(\nabla\smalll)(\vartheta)=0}\qqandqq\textcolor{magenta}{\convexcoe\|\theta\|^2\leq \theta^\top(\operatorname{Hess}\smalll) (\vartheta)\theta\leq \convexcoe\chi \| \theta \|^2}.
 \end{equation}  
 Then 
there exist $\fC,\delta \in (0,\infty)$ such that for every $\gamma\colon\N\to(\Gamma-\delta,\Gamma+\delta)$, every $\Theta=(\Theta^1,\dots,\Theta^{d})\colon \mathbb{N}_0\to\R^d$, every $\bfm=(\bfm^1,\dots,\bfm^d)\colon\N_0\to\R^d$, and every $\bbV=(\bbV^1,\dots,\bbV^{d})\colon \mathbb{N}_0\to\R^d$ with the property that for all $n \in \mathbb{N}$, $i\in\{1,2,\dots,d\}$ it holds that 
\begin{equation}
    \textcolor{magenta}{\bfm_0=0},\qquad \textcolor{magenta}{\bfm_n= \alpha \bfm_{n-1}+(1-\alpha)(\nabla\smalll)(\Theta_{n-1})},
\end{equation}
         \begin{equation}
            \textcolor{magenta}{\bbV_0=0},\qquad \textcolor{magenta}{\bbV_n^i= \beta \bbV_{n-1}^i+(1-\beta)|(\nabla_{\theta_i}\smalll)(\Theta_{n-1})|^2},
            \end{equation}
            \begin{equation}
            \textcolor{magenta}{\|\Theta_0-\vartheta\|< \delta},\qquad\text{and}\qquad  \textcolor{magenta}{\Theta_n^i=\Theta_{n-1}^i-\gamma_n
        \textstyle \bigl[\frac{\bfm_{n}^i}{1-\alpha^n}\bigr] \Bigl[\varepsilon+\bigl[\frac{\democrator_{n}^i}{1-\beta^n}\bigr]^{\nicefrac{1}{2}}\Bigr]^{-1}}
         \end{equation}
          we have for all $n\in\N_0$ that
$    \textcolor{magenta}{ \textstyle
    \|\Theta_n-\vartheta\|\leq \fC\|\Theta_0-\vartheta\|\varrho^{n}}$.
\end{athm}
\end{tcolorbox}
\end{samepage}
\begin{aproof}
     Throughout this proof for every $\xi\in\R^d$ and every $\gamma\colon\N\to\R$ let $\Theta^{\xi,\gamma}=(\Theta^{\xi,\gamma,1},\dots,\allowbreak\Theta^{\xi,\gamma,d})\colon\N_0\to\R^d$, $\bfm^{\xi,\gamma}=(\bfm^{\xi,\gamma,1},\dots,\bfm^{\xi,\gamma,d})\colon\N_0\to\R^d$, and  $\bbV^{\xi,\gamma}=(\bbV^{\xi,\gamma,1},\dots,\bbV^{\xi,\gamma,d})\colon\N_0\to\R^d$ satisfy for all $n\in\N$, $i\in\{1,2,\dots,d\}$ that
     \begin{equation}\llabel{def: bfm}
      \bfm^{\xi,\gamma}_0=0,\qquad \bfm^{\xi,\gamma}_n= \alpha \bfm^{\xi,\gamma}_{n-1}+(1-\alpha)(\nabla\smalll)(\Theta^{\xi,\gamma}_{n-1}),
            \end{equation}
              \begin{equation}\llabel{def: bbV}
      \bbV^{\xi,\gamma}_0=0,\qquad \bbV^{\xi,\gamma,i}_n= \beta\bbV^{\xi,\gamma,i}_{n-1}+(1-\beta)|(\nabla_{\theta_i}\smalll)(\Theta^{\xi,\gamma}_{n-1})|^2,
            \end{equation}
            \begin{equation}\llabel{def: Theta}
        \Theta^{\xi,\gamma}_0=\xi,\qquad\text{and}\qquad \Theta^{\xi,\gamma,i}_n=\Theta^{\xi,\gamma,i}_{n-1}-\gamma_n
        \textstyle \bigl[\frac{\bfm_{n}^{\xi,\gamma,i}}{1-\alpha^n}\bigr] \Bigl[\varepsilon+\bigl[\frac{\democrator_{n}^{\xi,\gamma,i}}{1-\beta^n}\bigr]^{\nicefrac{1}{2}}\Bigr]^{-1},
         \end{equation}
          let $I=(I_{i,j})_{(i,j)\in\{1,2,\dots,d\}\times\{1,2,\dots,d\}}\in\R^{d\times d}$ satisfy for all $i,j\in\{1,2,\dots,d\}$ that $I_{i,j}=\mathbbm 1_{\{j\}}(i)$,
       for every $n\in \N$, $\xi\in\R^d$ and every $\gamma\colon\N\to\R$ let $\Lambda^{\xi,\gamma}_n=(\Lambda^{\xi,\gamma,i,j}_n)_{(i,j)\in\{1,2,\dots,d\}\times\{1,2,\dots,d\}}\in\R^{d\times d}$ satisfy for all $i,j\in\{1,2,\dots,d\}$ that
      \begin{equation}\llabel{def: Lambda}
\Lambda^{\xi,\gamma,i,j}_n=\gamma_n[\varepsilon+\textstyle\bigl[\frac{\bbV_n^{\xi,\gamma,i}}{1-\beta^n}\bigr]^{1/2}]^{-1}\mathbbm 1_{\{j\}}(i),
      \end{equation}
         and let $\cA\in \R^{2d\times 2d}$ satisfy 
         \begin{equation}\llabel{def: cA}
           \cA=\begin{pmatrix}
               I-\Gamma\varepsilon^{-1}(1-\alpha)(\operatorname{Hess}\smalll)(\vartheta)&-\Gamma\varepsilon^{-1}\alpha I\\
               (1-\alpha)(\operatorname{Hess}\smalll)(\vartheta) & \alpha I
           \end{pmatrix}.
         \end{equation}  
         \argument{\lref{def: bfm};\lref{def: Theta};\lref{def: Lambda}}{that for every $n\in\N$, $\xi\in\R^d$, $i\in\{1,2,\dots,d\}$ and every $\gamma\colon\N\to\R$ it holds that
      \begin{equation}\llabel{eq1}
      \begin{split}
         & \Theta_n^{\xi,\gamma}-\vartheta=\Theta_{n-1}^{\xi,\gamma}-\vartheta-\frac{\Lambda_n^{\xi,\gamma}\bfm_n^{\xi,\gamma}}{1-\alpha^n}=\Theta_{n-1}^{\xi,\gamma}-\vartheta-\Gamma\varepsilon^{-1}\bfm_n^{\xi,\gamma}-\biggl(\frac{\Lambda_n^{\xi,\gamma}}{1-\alpha^n}-\Gamma\varepsilon^{-1}I\bigg)\bfm_n^{\xi,\gamma}\\
         &=\Theta_{n-1}^{\xi,\gamma}-\vartheta-\Gamma\varepsilon^{-1}(\alpha\bfm_{n-1}^{\xi,\gamma}+(1-\alpha)(\nabla\smalll)(\Theta_{n-1}^{\xi,\gamma}))-\biggl(\frac{\Lambda_n^{\xi,\gamma}}{1-\alpha^n}-\Gamma\varepsilon^{-1}I\bigg)\bfm_n^{\xi,\gamma}\\
         &=\Theta_{n-1}^{\xi,\gamma}-\vartheta-\Gamma\varepsilon^{-1}\alpha\bfm_{n-1}^{\xi,\gamma}-\Gamma\varepsilon^{-1}(1-\alpha)(\operatorname{Hess}\smalll)(\vartheta)(\Theta_{n-1}^{\xi,\gamma}-\vartheta) \\
         &-\Gamma\varepsilon^{-1}(1-\alpha)((\nabla\smalll)(\Theta_{n-1}^{\xi,\gamma})-(\operatorname{Hess}\smalll)(\vartheta)(\Theta_{n-1}^{\xi,\gamma}-\vartheta))   
        -\biggl(\frac{\Lambda_n^{\xi,\gamma}}{1-\alpha^n}-\Gamma\varepsilon^{-1}I\biggr)\bfm_n^{\xi,\gamma}.
          \end{split}
      \end{equation}}
      \argument{\lref{def: bfm};}{that for every $n\in\N$, $\xi\in\R^d$, $i\in\{1,2,\dots,d\}$ and every $\gamma\colon\N\to\R$ it holds that
      \begin{equation}\llabel{eq2}
      \begin{split}
          &\bfm^{\xi,\gamma}_n= \alpha \bfm^{\xi,\gamma}_{n-1}+(1-\alpha)(\nabla\smalll)(\Theta^{\xi,\gamma}_{n-1})\\
          &=\alpha \bfm^{\xi,\gamma}_{n-1}+(1-\alpha)(\operatorname{Hess}\smalll)(\vartheta)(\Theta_{n-1}^{\xi,\gamma}-\vartheta)+(1-\alpha)((\nabla\smalll)(\Theta^{\xi,\gamma}_{n-1})-(\operatorname{Hess}\smalll)(\vartheta)(\Theta_{n-1}^{\xi,\gamma}-\vartheta))
          \end{split}
      \end{equation}}
      \argument{\lref{def: cA};\lref{eq1};\lref{eq2};}{that for every $n\in\N$, $\xi\in\R^d$, $i\in\{1,2,\dots,d\}$ and every $\gamma\colon\N\to\R$ it holds that
      \begin{equation}\llabel{eq3}
      \begin{split}
          &(\Theta_n^{\xi,\gamma}-\vartheta,\bfm^{\xi,\gamma}_n)\\&=\cA(\Theta_{n-1}^{\xi,\gamma}-\vartheta,\bfm^{\xi,\gamma}_{n-1})-\Gamma\varepsilon^{-1}(1-\alpha)((\nabla\smalll)(\Theta_{n-1}^{\xi,\gamma})-(\operatorname{Hess}\smalll)(\vartheta)(\Theta_{n-1}^{\xi,\gamma}-\vartheta),0)\\
          &-\textstyle\bigl(\bigl(\frac{\Lambda_n^{\xi,\gamma}}{1-\alpha^n}-\Gamma\varepsilon^{-1}I\bigr)\bfm_n^{\xi,\gamma},0\bigr)+(1-\alpha)(0,(\nabla\smalll)(\Theta^{\xi,\gamma}_{n-1})-(\operatorname{Hess}\smalll)(\vartheta)(\Theta_{n-1}^{\xi,\gamma}-\vartheta))
          \end{split}
      \end{equation}}
        \argument{the fact that $\alpha=(\frac{\sqrt{\chi}-1}{\sqrt{\chi}+1})^2$}{that
    \begin{equation}\llabel{eqq0}
        \textstyle\chi=(\frac{1+\sqrt{\alpha}}{1-\sqrt{\alpha}})^2.
    \end{equation}}
    \argument{\lref{eqq0};the assumption that $\Gamma=(\convexcoe\sqrt{\chi})^{-1}\varepsilon$;}{that
    \begin{equation}\llabel{eqq1}
    \begin{split}
          &(1+\alpha-\Gamma\varepsilon^{-1}\convexcoe(1-\alpha))^2= (1+\alpha-\chi^{-1/2}(1-\alpha))^2=\biggl(1+\alpha-\frac{(1-\alpha)(1-\sqrt{\alpha})}{1+\sqrt{\alpha}}\biggr)^2=4\alpha\\
          &\biggl(1+\alpha-\frac{(1-\alpha)(1+\sqrt{\alpha})}{1-\sqrt{\alpha}}\biggr)^2=(1+\alpha-\chi^{1/2}(1-\alpha))^2=(1+\alpha-\Gamma\varepsilon^{-1}\convexcoe\chi(1-\alpha))^2.
          \end{split}
    \end{equation}}
    \argument{\lref{eqq1};the fact that $\Gamma\varepsilon^{-1}>0$}{that for all $\lambda\in [\convexcoe,\convexcoe\chi]$ it holds that
    \begin{equation}\llabel{eqq1'}
        (1+\alpha-\Gamma\varepsilon^{-1} \lambda(1-\alpha))^2\leq 4\alpha.
    \end{equation}}
       \argument{\lref{def: bbL};}{that for all $\lambda\in \R$, $\theta\in \bbC^{d}\backslash\{0\}$ with $(\operatorname{Hess} \smalll)(\vartheta)\theta=\lambda\theta$ it holds that
        \begin{equation}\llabel{eq4.1}
\convexcoe\leq \lambda\leq \convexcoe\chi.
        \end{equation}}
        \argument{\lref{eq4.1}; \lref{eqq1'}}
      {that for all $\lambda\in \R$, $\theta\in \bbC^{d}\backslash\{0\}$ with $(\operatorname{Hess} \smalll)(\vartheta)\theta=\lambda\theta$ it holds that
        \begin{equation}\llabel{eq4}
          \textstyle\bigl(  \frac{1+\alpha-\Gamma\varepsilon^{-1}\lambda(1-\alpha)}{2}\bigr)^2+\bigl(  \frac{1+\alpha-\Gamma\varepsilon^{-1}\lambda(1-\alpha)}{2}\bigr)^2-\alpha=\frac{(1+\alpha-\Gamma\varepsilon^{-1}\lambda(1-\alpha))^2}{2}-\alpha\leq \alpha.
        \end{equation}}
\argument{\lref{def: cA};\lref{eq4};\cite[item (i) in Lemma 3.2]{DeArAdsharp2025} (applied with $d\curvearrowleft d$, $\alpha\curvearrowleft\alpha$, $\gamma\curvearrowleft \Gamma\varepsilon^{-1}$, $H\curvearrowleft(\operatorname{Hess}\smalll)(\vartheta)$,  $A\curvearrowleft \cA$ in the notation of \cite[Lemma 3.2]{DeArAdsharp2025})}{that for all $\lambda\in \bbC$, $\theta\in \bbC^{2d}\backslash\{0\}$ with $\cA\theta=\lambda\theta$ it holds that
\begin{equation}\llabel{eq5}
    |\lambda|\leq \sqrt{\alpha}.
\end{equation}}
\argument{\lref{eq5};the fact that $\varrho>\sqrt{\alpha}$;\unskip, \eg, \cite[Lemma 2.3]{DeArAdsharp2025}}{that there exists a norm $\nnorm{\cdot}\colon \R^{2d}\to[0,\infty)$ which satisfies for all $\theta\in\R^{2d}$ that
\begin{equation}\llabel{evidence 1}
    \nnorm{\cA\theta}\leq \frac{\varrho+\sqrt{\alpha}}{2}\nnorm{\theta}.
\end{equation}}
\startnewargseq
  \argument{the fact that $\|\cdot\|$ is a norm; the fact that $\nnorm{\cdot}$ is a norm; the fact that $\R^{2d}$ is finite dimensional vector space;}{that there exist $\mu_1,\mu_2\in (0,\infty)$ which satisfy for all $\theta\in \R^{2d}$ that
         \begin{equation}\llabel{def: mu}
             \mu_1\nnorm{\theta}\leq \|\theta\|\leq \mu_2\nnorm{\theta}.
         \end{equation}}
\startnewargseq
   \argument{\lref{def: bbL};the assumption that $\smalll\in C^2(\R^d,\R)$;the fact that $\varrho>\sqrt{\alpha}$}{that there exists $\delta\in (0,\infty)$ such that for all $\theta\in \R^d$ with $\|\theta-\vartheta\|\leq \delta$ it holds that
    \begin{equation}\llabel{eq7}
    \begin{split}
       & \|(\nabla\smalll)(\theta)-(\operatorname{Hess}\smalll)(\vartheta)(\theta-\vartheta)\|\\&=\|(\nabla\smalll)(\theta)-(\nabla\smalll)(\vartheta)-(\operatorname{Hess}\smalll)(\vartheta)(\theta-\vartheta)\|\leq \textstyle\min\bigl\{\frac{(\varrho-\sqrt{\alpha})\varepsilon\mu_1}{8\mu_2\Gamma(1-\alpha)},\frac{(\varrho-\sqrt{\alpha})\mu_1}{8\mu_2(1-\alpha)}\bigr\}\|\theta-\vartheta\|.
        \end{split}
    \end{equation}}
    \argument{\lref{eq7};\lref{def: mu}}{that there exists $\delta\in (0,\infty)$ such that for every $n\in\N$, $\xi\in\R^d$ and every $\gamma \colon\N\to\R$ with $\|\Theta_{n-1}^{\xi,\gamma}-\vartheta\|\leq \delta$ it holds that
    \begin{equation}\llabel{eq8}
    \begin{split}
       &\Gamma\varepsilon^{-1}(1-\alpha)\nnorm{((\nabla\smalll)(\Theta_{n-1}^{\xi,\gamma})-(\operatorname{Hess}\smalll)(\vartheta)(\Theta_{n-1}^{\xi,\gamma}-\vartheta),0)}\\
       &\leq \Gamma\varepsilon^{-1}(1-\alpha)(\mu_1)^{-1}\|((\nabla\smalll)(\Theta_{n-1}^{\xi,\gamma})-(\operatorname{Hess}\smalll)(\vartheta)(\Theta_{n-1}^{\xi,\gamma}-\vartheta),0)\|\\
       &=\Gamma\varepsilon^{-1}(1-\alpha)(\mu_1)^{-1}\|(\nabla\smalll)(\Theta_{n-1}^{\xi,\gamma})-(\operatorname{Hess}\smalll)(\vartheta)(\Theta_{n-1}^{\xi,\gamma}-\vartheta)\|\\
       &\leq \frac{\varrho-\sqrt{\alpha}}{8\mu_2}\|\Theta_{n-1}^{\xi,\gamma}-\vartheta\|\leq \frac{\varrho-\sqrt{\alpha}}{8\mu_2}\|(\Theta_{n-1}^{\xi,\gamma}-\vartheta,\bfm_{n-1}^{\xi,\gamma})\|\leq \frac{\varrho-\sqrt{\alpha}}{8}\nnorm{(\Theta_{n-1}^{\xi,\gamma}-\vartheta,\bfm_{n-1}^{\xi,\gamma})}
       \end{split}
    \end{equation}
    and 
    \begin{equation}\llabel{eq9}
    \begin{split}
       &(1-\alpha)\nnorm{(0,(\nabla\smalll)(\Theta_{n-1}^{\xi,\gamma})-(\operatorname{Hess}\smalll)(\vartheta)(\Theta_{n-1}^{\xi,\gamma}-\vartheta))}\\
       &\leq (1-\alpha)(\mu_1)^{-1}\|(0,(\nabla\smalll)(\Theta_{n-1}^{\xi,\gamma})-(\operatorname{Hess}\smalll)(\vartheta)(\Theta_{n-1}^{\xi,\gamma}-\vartheta))\|\\
       &=(1-\alpha)(\mu_1)^{-1}\|(\nabla\smalll)(\Theta_{n-1}^{\xi,\gamma})-(\operatorname{Hess}\smalll)(\vartheta)(\Theta_{n-1}^{\xi,\gamma}-\vartheta)\|\\
       &\leq \frac{\varrho-\sqrt{\alpha}}{8\mu_2}\|\Theta_{n-1}^{\xi,\gamma}-\vartheta\|\leq \frac{\varrho-\sqrt{\alpha}}{8\mu_2}\|(\Theta_{n-1}^{\xi,\gamma}-\vartheta,\bfm_{n-1}^{\xi,\gamma})\|\leq \frac{\varrho-\sqrt{\alpha}}{8}\nnorm{(\Theta_{n-1}^{\xi,\gamma}-\vartheta,\bfm_{n-1}^{\xi,\gamma})}.
       \end{split}
    \end{equation}}
    \argument{\lref{eq9};}{that there exists $\delta\in (0,\infty)$ such that for every $n\in\N$, $\xi\in\R^d$ and every $\gamma \colon\N\to\R$ with $\|\Theta_{n-1}^{\xi,\gamma}-\vartheta\|\leq \delta$ it holds that
    \begin{equation}\llabel{evidence 2}
        \begin{split}
            &\Gamma\varepsilon^{-1}(1-\alpha)\nnorm{((\nabla\smalll)(\Theta_{n-1}^{\xi,\gamma})-(\operatorname{Hess}\smalll)(\vartheta)(\Theta_{n-1}^{\xi,\gamma}-\vartheta),0)}\\
            &+(1-\alpha)\nnorm{(0,(\nabla\smalll)(\Theta_{n-1}^{\xi,\gamma})-(\operatorname{Hess}\smalll)(\vartheta)(\Theta_{n-1}^{\xi,\gamma}-\vartheta))}\\
            &\leq \frac{\varrho-\sqrt{\alpha}}{4}\nnorm{(\Theta_{n-1}^{\xi,\gamma}-\vartheta,\bfm_{n-1}^{\xi,\gamma})}.
        \end{split}
    \end{equation}}
    \argument{\lref{def: bbL};the fact that $\smalll\in C^2(\R^d,\R)$}{that there exist $\delta,\fC_1,\fC_2\in (0,\infty)$ such that for all $\theta\in \R^d$ with $\|\theta-\vartheta\|\leq \delta$ it holds that
         \begin{equation}\llabel{eq11}
            \fC_1\|\theta-\vartheta\|\leq \|(\nabla\smalll)(\theta)\|=\|(\nabla\smalll)(\theta)-(\nabla\smalll)(\vartheta)\|\leq \fC_2\|\theta-\vartheta\|.
         \end{equation}}
\argument{\lref{def: bfm};\lref{eq11};\lref{def: mu};}{that there exist $\delta,\fC_1,\fC_2\in (0,\infty)$ such that for every $n\in\N$, $\xi\in\R^d$ and every $\gamma\colon\N\to\R$ with $\|\Theta_{n-1}^{\xi,\gamma}-\vartheta\|\leq \delta$ it holds that
         \begin{equation}\llabel{eq14'}
         \begin{split}
             &\|\bfm_{n}^{\xi,\gamma}\|=\|\alpha\bfm_{n-1}^{\xi,\gamma}+(1-\alpha)(\nabla\smalll)(\Theta_{n-1}^{\xi,\gamma})\|\leq \alpha \|\bfm_{n-1}^{\xi,\gamma}\|+(1-\alpha)\|(\nabla\smalll)(\Theta_{n-1}^{\xi,\gamma})\|\\
             &\leq \fC_1\|\bfm_{n-1}^{\xi,\gamma}\|+\fC_1\|\Theta_{n-1}^{\xi,\gamma}-\vartheta\|\leq 2\fC_1\|(\Theta_{n-1}^{\xi,\gamma}-\vartheta,\bfm_{n-1}^{\xi,\gamma})\|\leq \fC_2\nnorm{(\Theta_{n-1}^{\xi,\gamma}-\vartheta,\bfm_{n-1}^{\xi,\gamma})}.
             \end{split}
         \end{equation}}
  \argument{\lref{def: bbV};\lref{eq11};\unskip, \eg, \cite[Lemma 2.7]{DeRoArAd2025}}{that there exist $\delta,\fC\in (0,\infty)$ such that for every $n\in\N$, $\xi\in\R^d$, $i\in\{1,2,\dots,d\}$ and every $\gamma\colon\N\to\R$ with $\sup_{k\in\{0,1,\dots,n-1\}}\|\Theta_{k}^{\xi,\gamma}-\vartheta\|\leq \delta$ it holds that
  and
    \begin{equation}\llabel{eq14}
        \bbV_n^{\xi,\gamma, i}=(1-\beta)\bigl[\textstyle\sum_{k=1}^n\beta^{n-k}|(\nabla_{\theta_i}\smalll)(\Theta_{k-1}^{\xi,\gamma})|^2 \textstyle \bigr]\leq (1-\beta)\bigl[\textstyle\sum_{k=1}^n\beta^{n-k}\fC\|\Theta_{k-1}^{\xi,\gamma}-\vartheta\|^2 \textstyle \bigr]\leq \fC\delta^2.
    \end{equation}}
    \argument{\lref{def: Lambda};\lref{eq14};}{that there exist $\delta,\fC\in (0,\infty)$ such that for every $n\in\N$, $\xi\in\R^d$ and every $\gamma\colon\N\to\R$ with $\sup_{k\in\{0,1,\dots,n-1\}}\|\Theta_{k}^{\xi,\gamma}-\vartheta\|\leq \delta$ and $\sup_{k\in\N}|\gamma_k|\leq 2\Gamma$ it holds that  
    \begin{equation}\llabel{eq15}
    \begin{split}
       \textstyle\|\gamma_n\varepsilon^{-1} I-\Lambda_n^{\xi,\gamma}\|\leq \sum\limits_{i=1}^d|\gamma_n|\bigl|\bigl(\varepsilon+\bigl(\frac{\bbV_n^{\xi,\gamma,i}}{1-\beta^n}\bigr)^{1/2}\bigr)^{-1}-\varepsilon^{-1}\bigr|\leq \sum\limits_{i=1}^d\bigl(2\Gamma\bigl|\frac{\bbV_n^{\xi,\gamma,i}}{1-\beta}\bigr|^{1/2}\varepsilon^{-2}\bigr)\leq \fC\delta.
        \end{split}
    \end{equation}}
    \argument{\lref{eq15};\lref{def: mu};\lref{eq14'}}{that there exist $\delta,\fC_1,\fC_2\in (0,\infty)$ such that for every $n\in\N$, $\xi\in\R^d$ and every $\gamma\colon\N\to\R$ with $\sup_{k\in\{0,1,\dots,n-1\}}\|\Theta_{k}^{\xi,\gamma}-\vartheta\|\leq \delta$ and $\sup_{k\in\N}|\gamma_k|\leq 2\Gamma$ it holds that 
    \begin{equation}\llabel{eq16}
        \begin{split}
            \nnorm{((\gamma_n\varepsilon^{-1}I-\Lambda_n^{\xi,\gamma})\bfm_{n}^{\xi,\gamma},0)}&\leq \fC_1\|((\gamma_n\varepsilon^{-1}I-\Lambda_n^{\xi,\gamma})\bfm_{n}^{\xi,\gamma},0)\|=\fC_1\|(\gamma_n\varepsilon^{-1}I-\Lambda_n^{\xi,\gamma})\bfm_{n}^{\xi,\gamma}\|\\
           &\leq \fC_1\|\gamma_n\varepsilon^{-1}I-\Lambda_n^{\xi,\gamma}\|\|\bfm_{n}^{\xi,\gamma}\|\leq \fC_2\delta\nnorm{(\Theta_{n-1}^{\xi,\gamma}-\vartheta,\bfm_{n-1}^{\xi,\gamma})}.
        \end{split}
    \end{equation}}
    \argument{\lref{def: mu};\lref{eq14'};}{that there exists $\fC\in (0,\infty)$ such that for every $n\in\N$, $\xi\in\R^d$ and every $\gamma\colon\N\to\R$ with $\|\Theta_{n-1}^{\xi,\gamma}-\vartheta\|\leq 1$ it holds that
    \begin{equation}\llabel{eq17}
       \nnorm{(\bfm_n^{\xi,\gamma},0)}\leq (\mu_1)^{-1}\|(\bfm_n^{\xi,\gamma},0)\|=(\mu_1)^{-1}\|\bfm_n^{\xi,\gamma}\|\leq \fC\nnorm{(\Theta_{n-1}^{\xi,\gamma}-\vartheta,\bfm_{n-1}^{\xi,\gamma})}.
    \end{equation}}
    \argument{\lref{eq17};\lref{eq16};the fact that $\nnorm{\cdot}$ is a norm}{that there exist $\delta,\fC\in (0,\infty)$ such that for every $n\in\N$, $\xi\in\R^d$ and every $\gamma\colon\N\to\R$ with $\sup_{k\in\{0,1,\dots,n-1\}}\|\Theta_{k}^{\xi,\gamma}-\vartheta\|\leq \delta$ and $\sup_{k\in\N}|\gamma_k|\leq 2\Gamma$ it holds that 
    \begin{equation}\llabel{eq18}
    \begin{split}
    &\textstyle\bignnorm{\bigl(\bigl(\frac{\Lambda_n^{\xi,\gamma}}{1-\alpha^n}-\Gamma\varepsilon^{-1}I\bigr)\bfm_n^{\xi,\gamma},0\bigr)}\\
    &\leq\textstyle  \frac{1}{1-\alpha^n}\textstyle\nnorm{((\Lambda_n^{\xi,\gamma}-\gamma_n\varepsilon^{-1}I)\bfm_n^{\xi,\gamma},0)}+ \frac{\alpha^n\Gamma\varepsilon^{-1}}{1-\alpha^n}\textstyle\nnorm{(\bfm_n^{\xi,\gamma},0)}+\frac{1}{1-\alpha^n}\nnorm{(\Gamma\varepsilon^{-1}I-\gamma_n\varepsilon^{-1}I)\bfm_n^{\xi,\gamma},0}\\
        &\leq \textstyle\frac{1}{1-\alpha}\textstyle\nnorm{((\Lambda_n^{\xi,\gamma}-\gamma_n\varepsilon^{-1}I)\bfm_n^{\xi,\gamma},0)}+ \frac{\alpha^n\Gamma\varepsilon^{-1}}{1-\alpha}\textstyle\nnorm{(\bfm_n^{\xi,\gamma},0)}+\frac{|\gamma_n-\Gamma|\varepsilon^{-1}}{1-\alpha}\textstyle\nnorm{(\bfm_n^{\xi,\gamma},0)}\\
        &\leq \fC\delta\textstyle\nnorm{(\Theta_{n-1}^{\xi,\gamma}-\vartheta,\bfm_{n-1}^{\xi,\gamma})}+ \fC\alpha^n\nnorm{(\Theta_{n-1}^{\xi,\gamma}-\vartheta,\bfm_{n-1}^{\xi,\gamma})}+\fC|\gamma_n-\Gamma|\nnorm{(\Theta_{n-1}^{\xi,\gamma}-\vartheta,\bfm_{n-1}^{\xi,\gamma})}.
        \end{split}
    \end{equation}}
    \argument{\lref{eq18};the fact that $\varrho>\sqrt{\alpha}$}{that there exist $\delta\in (0,\infty)$, $N\in\N$ such that for every $n\in\N\cap[N,\infty)$, $\xi\in\R^d$ and every $\gamma\colon\N\to\R$ with $\sup_{k\in\{0,1,\dots,n-1\}}\|\Theta_{k}^{\xi,\gamma}-\vartheta\|\leq \delta$ and $\sup_{k\in\N}|\gamma_k-\Gamma|\leq \delta$ it holds that 
    \begin{equation}\llabel{evidence 3}
       \textstyle\bignnorm{\bigl(\bigl(\frac{\Lambda_n^{\xi,\gamma}}{1-\alpha^n}-\Gamma\varepsilon^{-1}I\bigr)\bfm_n^{\xi,\gamma},0\bigr)}\leq \frac{\varrho-\sqrt{\alpha}}{4}\nnorm{(\Theta_{n-1}^{\xi,\gamma}-\vartheta,\bfm_{n-1}^{\xi,\gamma})}.
    \end{equation}}   
    \argument{\lref{eq3};\lref{evidence 1};\lref{evidence 2};\lref{evidence 3}}{that there exist $\delta\in (0,\min\{\Gamma,1\})$, $N\in\N$ which satisfy for every $n\in\N\cap[N,\infty)$, $\xi\in\R^d$ and every $\gamma\colon\N\to\R$ with $\sup_{k\in\{0,1,\dots,n-1\}}\|\Theta_{k}^{\xi,\gamma}-\vartheta\|\leq \delta$ and $\sup_{k\in\N}|\gamma_k-\Gamma|\leq \delta$ that
    \begin{equation}\llabel{eqmain}
        \nnorm{(\Theta_n^{\xi,\gamma}-\vartheta,\bfm^{\xi,\gamma}_n)}\leq \varrho \nnorm{(\Theta_{n-1}^{\xi,\gamma}-\vartheta,\bfm^{\xi,\gamma}_{n-1})}.
    \end{equation}}
    \startnewargseq
    \argument{\lref{def: Theta};\lref{def: mu};\lref{eq14'};the fact that $\delta<\Gamma$}{that there exists $\fC\in (0,\infty)$ such that for every $n\in\N$, $\xi\in\R^d$ and every $\gamma\colon\N\to\R$ with $\|\Theta_{n-1}^{\xi,\gamma}-\vartheta\|\leq \delta$ and $\sup_{k\in\N}|\gamma_k-\Gamma|\leq \delta$ it holds that
    \begin{equation}\llabel{eqt1}
        \|\Theta_{n}^{\xi,\gamma}-\vartheta\|\leq \|\Theta_{n-1}^{\xi,\gamma}-\vartheta\|+\frac{2\Gamma\varepsilon^{-1}}{1-\alpha}\|\bfm_n^{\xi,\gamma}\| \leq \fC \nnorm{(\Theta_{n-1}^{\xi,\gamma}-\vartheta,\bfm_{n-1}^{\xi,\gamma})}.
    \end{equation}}
    \argument{\lref{eqt1};\lref{eq14'};}{that there exists $\fC\in (0,\infty)$ such that for every $n\in\N$, $\xi\in\R^d$ and every $\gamma\colon\N\to\R$ with $\sup_{k\in\{0,1,\dots,n-1\}}\|\Theta_{k}^{\xi,\gamma}-\vartheta\|\leq \delta$ and $\sup_{k\in\N}|\gamma_k-\Gamma|\leq \delta$ it holds that
    \begin{equation}\llabel{eqt2}
        \nnorm{(\Theta_{n}^{\xi,\gamma}-\vartheta,\bfm_n^{\xi,\gamma})}\leq \fC \nnorm{(\Theta_{n-1}^{\xi,\gamma}-\vartheta,\bfm_{n-1}^{\xi,\gamma})}.
    \end{equation}}
    \argument{\lref{eqt2};\lref{def: mu};induction}{that there exist $\fC_1,\fC_2\in (0,\infty)$ such that for every $n\in\N_0$, $\xi\in\R^d$ and every $\gamma\colon\N\to\R$ with $\sup_{k\in\{0,1,\dots,n-1\}}\|\Theta_{k}^{\xi,\gamma}-\vartheta\|\leq \delta$ and $\sup_{k\in\N}|\gamma_k-\Gamma|\leq \delta$ it holds that
    \begin{equation}\llabel{eqt3}
    \begin{split}
       \|\Theta_n^{\xi,\gamma}-\vartheta\|&\leq \|(\Theta_n^{\xi,\gamma}-\vartheta,\bfm_n^{\xi,\gamma})\|\leq   \mu_2\nnorm{(\Theta_{n}^{\xi,\gamma}-\vartheta,\bfm_n^{\xi,\gamma})}\leq \mu_2(\fC_1)^n\nnorm{(\Theta_{0}^{\xi,\gamma}-\vartheta,\bfm_0^{\xi,\gamma})}\\
       &=\mu_2(\fC_1)^n\nnorm{(\xi-\vartheta,0)}\leq \fC_2(\fC_1)^n\|\xi-\vartheta\|.
       \end{split}
    \end{equation}}
    \argument{\lref{eqt3};}{that there exists $\eta\in (0,\delta)$ which satisfies for every $n\in\{0,1,\dots,N\}$, $\xi\in\R^d$ and every $\gamma\colon\N\to\R$ with $\|\xi-\vartheta\|\leq \eta$ and $\sup_{k\in\N}|\gamma_k-\Gamma|\leq \eta$ that
    \begin{equation}\llabel{def: eta}
        \|\Theta_n^{\xi,\gamma}-\vartheta\|\leq \delta.
    \end{equation}}
    \startnewargseq
    \argument{\lref{eqt3};\lref{def: eta};the fact that $\eta<\delta<1$}{that there exists $\fC\in (1,\infty)$ which satisfies for every $n\in\{0,1,\dots,N\}$, $\xi\in\R^d$ and every $\gamma\colon\N\to\R$ with $\|\xi-\vartheta\|\leq \eta$ and $\sup_{k\in\N}|\gamma_k-\Gamma|\leq \eta$ it holds that
    \begin{equation}\llabel{def: fC}
        \nnorm{(\Theta_n^{\xi,\gamma}-\vartheta,\bfm^{\xi,\gamma}_n)}\leq \fC\|\xi-\vartheta\|\varrho^n.
    \end{equation}}
     In the following we prove that for all $n\in\N\cap[N,\infty)$, $\xi\in \R^d$ and all $\gamma\colon\N\to\R$ with $\sup_{k\in\N}|\gamma_k-\Gamma|< \eta$ and $\|\xi-\vartheta\|\leq  \min\{\eta,\delta(\mu_2\fC)^{-1}\}$ it holds that
         \begin{equation}\llabel{need to prove}
             \nnorm{(\Theta_n^{\xi,\gamma}-\vartheta,\bfm^{\xi,\gamma}_n)}\leq\fC\|\xi-\vartheta\| \varrho^n.
         \end{equation}
         We prove \lref{need to prove} by induction on $n\in\N\cap[N,\infty)$. Note that \lref{def: fC} 
          shows \lref{need to prove} in the base case $n=N$. For the induction step we assume that there exists $n\in\N\cap[N+1,\infty)$ which satisfies for all $k\in\{N,N+1,\dots,n-1\}$, $\xi\in \R^d$ and all $\gamma\colon\N\to\R$ with $\sup_{m\in\N}|\gamma_m-\Gamma|< \eta$ and $\|\xi-\vartheta\|\leq \min\{\eta,\delta(\mu_2\fC)^{-1}\}$ that
         \begin{equation}\llabel{ntp: assume}
              \nnorm{(\Theta_{k}^{\xi,\gamma}-\vartheta,\bfm_k^{\xi,\gamma})}\leq\fC\|\xi-\vartheta\|\varrho^{k}.
         \end{equation}
         \startnewargseq
         \argument{\lref{def: fC};\lref{ntp: assume};}{that for all $k\in\{0,1,\dots,n-1\}$, $\xi\in \R^d$ and all $\gamma\colon\N\to\R$ with $\sup_{m\in\N}|\gamma_m-\Gamma|< \eta$ and $\|\xi-\vartheta\|\leq \min\{\eta,\delta(\mu_2\fC)^{-1}\}$ it holds that
         \begin{equation}\llabel{ntp: assume'}
             \nnorm{(\Theta_{k}^{\xi,\gamma}-\vartheta,\bfm_k^{\xi,\gamma})}\leq\fC\|\xi-\vartheta\|\varrho^{k}.
         \end{equation}}
         \argument{\lref{def: mu};\lref{ntp: assume'};the fact that $0<\varrho<1$}{that for all $k\in\{1,2,\dots,n-1\}$, $\xi\in \R^d$ and all $\gamma\colon\N\to\R$ with $\sup_{m\in\N}|\gamma_m-\Gamma|<\eta$ and $\|\xi-\vartheta\|\leq \min\{\eta,\delta(\mu_2\fC)^{-1}\}$ it holds that
         \begin{equation}\llabel{ntp: eq1}
         \begin{split}
           \textstyle  \|\Theta_{k}^{\xi,\gamma}-\vartheta\| \leq \|(\Theta_{k}^{\xi,\gamma}-\vartheta,\bfm_k^{\xi,\gamma})\|\leq \mu_2\nnorm{(\Theta_{k}^{\xi,\gamma}-\vartheta,\bfm_k^{\xi,\gamma})}\leq \mu_2\fC\varrho^{k}\|\xi-\vartheta\|\leq \mu_2\fC\|\xi-\vartheta\|\leq \delta.
             \end{split}
         \end{equation}}
         \argument{\lref{ntp: eq1};\lref{eqmain};\lref{ntp: assume'}}{for all $\xi\in \R^d$ and all $\gamma\colon\N\to\R$ with $\sup_{k\in\N}|\gamma_k-\Gamma|< \eta$  and $\|\xi-\vartheta\|\leq \min\{\eta,\delta(\mu_2\fC)^{-1}\}$ that
         \begin{equation}\llabel{ntp: eq2}
             \begin{split}
             \nnorm{(\Theta_n^{\xi,\gamma}-\vartheta,\bfm_n^{\xi,\gamma})}\leq  \varrho\nnorm{(\Theta_{n-1}^{\xi,\gamma}-\vartheta,\bfm_{n-1}^{\xi,\gamma})}
            \leq \varrho\fC\|\xi-\vartheta\|\varrho^{n-1}=\fC\|\xi-\vartheta\|\varrho^{n}.
             \end{split}
         \end{equation}}
         \argument{\lref{ntp: eq2};\lref{ntp: assume};induction}{\lref{need to prove}\dott}
         \startnewargseq
         \argument{\lref{def: mu};\lref{need to prove};}{that for all $n\in\N_0$, $\xi\in \R^d$ and all $\gamma\colon\N\to\R$ with $\sup_{k\in\N}|\gamma_k-\Gamma|<\eta$ and $\|\xi-\vartheta\|\leq  \min\{\eta,\delta(\mu_2\fC)^{-1}\}$ it holds that
         \begin{equation}\llabel{eq19}
           \|\Theta_n^{\xi,\gamma}-\vartheta\|\leq   \|(\Theta_n^{\xi,\gamma}-\vartheta,\bfm_n^{\xi,\gamma})\|\leq \mu_2\nnorm{(\Theta_n^{\xi,\gamma}-\vartheta,\bfm_n^{\xi,\gamma})}\leq \mu_2\fC\|\xi-\vartheta\|\varrho^{n}=\mu_2\fC\|\Theta_0^{\xi,\gamma}-\vartheta\|\varrho^{n}.
         \end{equation}}
         \argument{\lref{eq19};\lref{def: bfm};\lref{def: Theta};}{that there exist $\scrC,\scrc \in (0,\infty)$ such that for every $\gamma\colon\N\to(\Gamma-\scrc,\Gamma+\scrc)$, every $\Theta=(\Theta^1,\dots,\Theta^{d})\colon \mathbb{N}_0\to\R^d$, every $\bfm=(\bfm^1,\dots,\bfm^d)\colon\N_0\to\R^d$, and every $\bbV=(\bbV^1,\dots,\bbV^{d})\colon \mathbb{N}_0\to\R^d$ with the property that for all $n \in \mathbb{N}$, $i\in\{1,2,\dots,d\}$ it holds that 
         \begin{equation}
  \bfm_0=0,\qquad \bfm_n= \alpha \bfm_{n-1}+(1-\alpha)(\nabla\smalll)(\Theta_{n-1}),
\end{equation}
         \begin{equation}
           \bbV_0=0,\qquad \bbV_n^i= \beta \bbV_{n-1}^i+(1-\beta)|(\nabla_{\theta_i}\smalll)(\Theta_{n-1})|^2,
            \end{equation}
            \begin{equation}
           \|\Theta_0-\vartheta\|< \scrc,\qquad\text{and}\qquad  \Theta_n^i=\Theta_{n-1}^i-\gamma_n
        \textstyle \bigl[\frac{\bfm_{n}^i}{1-\alpha^n}\bigr] \Bigl[\varepsilon+\bigl[\frac{\democrator_{n}^i}{1-\beta^n}\bigr]^{\nicefrac{1}{2}}\Bigr]^{-1}
        \end{equation}
          we have for all $n\in\N_0$ that
$    \textstyle
    \|\Theta_n-\vartheta\|\leq \scrC\|\Theta_0-\vartheta\|\varrho^{n}$\dott}
\end{aproof}
\begin{tcolorbox}[colback=white!95!gray,
                  colframe=black,
                  boxrule=0.5pt,
                  sharp corners,
                  enhanced,
                  breakable,
                 ]
\begin{athm}{cor}{main theorem higher bound adam}
Let $d\in\N$, $\chi\in (1,\infty)$, $\alpha=(\frac{ \sqrt{\chi}  - 1 }  { \sqrt{\chi}  + 1 })^2$, $\beta\in (0,1)$, $\varepsilon\in (0,\infty)$. Then for every $\convexcoe \in (0,\infty)$, $\vartheta \in \R^d$, $\smalll\in \allowbreak C^{2}(\R^d,\mathbb{R})$, $\varrho\in (\frac{ \sqrt{\chi}  - 1 }  { \sqrt{\chi}  + 1 },1)$ 
there exist $\fC,\Gamma, \delta \in (0,\infty)$ such that for every $\gamma\colon\N\to(\Gamma-\delta,\Gamma+\delta)$, every $\Theta=(\Theta^1,\dots,\Theta^d)\colon \mathbb{N}_0\to\R^d$, every $\bfm=(\bfm^1,\dots,\bfm^d)\colon\N_0\to\R^d$, and every $\bbV=(\bbV^1,\dots,\bbV^d)\colon \mathbb{N}_0\to\R^d$ with the property that for all $n \in \mathbb{N}$, $i\in\{1,2,\dots,d\}$, $\theta\in\R^d$ it holds that
\begin{equation}
    \textcolor{magenta}{(\nabla\smalll)(\vartheta)=0},\qquad\textcolor{magenta}{\convexcoe\|\theta\|^2\leq \theta^\top(\operatorname{Hess}\smalll) (\vartheta)\theta\leq \convexcoe\chi \| \theta \|^2},
\end{equation}
\begin{equation}
    \textcolor{magenta}{\bfm_0=0},\qquad \textcolor{magenta}{\bfm_n= \alpha \bfm_{n-1}+(1-\alpha)(\nabla\smalll)(\Theta_{n-1})},
\end{equation}
         \begin{equation}
            \textcolor{magenta}{\bbV_0=0},\qquad \textcolor{magenta}{\bbV_n^i= \beta \bbV_{n-1}^i+(1-\beta)|(\nabla_{\theta_i}\smalll)(\Theta_{n-1})|^2},
            \end{equation}
            \begin{equation}
            \textcolor{magenta}{\|\Theta_0-\vartheta\|< \delta},\qquad\text{and}\qquad  \textcolor{magenta}{\Theta_n^i=\Theta_{n-1}^i-\gamma_n
        \textstyle \bigl[\frac{\bfm_{n}^i}{1-\alpha^n}\bigr] \Bigl[\varepsilon+\bigl[\frac{\democrator_{n}^i}{1-\beta^n}\bigr]^{\nicefrac{1}{2}}\Bigr]^{-1}}
         \end{equation}
          we have for all $n\in\N_0$ that
$    \textcolor{magenta}{ \|\Theta_n-\vartheta\|\leq \fC\|\Theta_0-\vartheta\|\varrho^n}$.
\end{athm}
\end{tcolorbox}
\begin{aproof}
    \argument{\cref{theo: higher bound adam}}{that for every $\convexcoe \in (0,\infty)$, $\vartheta \in \R^d$, $\smalll\in \allowbreak C^{2}(\R^d,\mathbb{R})$, $\varrho\in (\frac{ \sqrt{\chi}  - 1 }  { \sqrt{\chi}  + 1 },1)$ 
there exist $\fC,\Gamma, \delta \in (0,\infty)$ such that for every $\gamma\colon\N\to(\Gamma-\delta,\Gamma+\delta)$, every $\Theta=(\Theta^1,\dots,\Theta^d)\colon \mathbb{N}_0\to\R^d$, every $\bfm=(\bfm^1,\dots,\bfm^d)\colon\N_0\to\R^d$, and every $\bbV=(\bbV^1,\dots,\bbV^d)\colon \mathbb{N}_0\to\R^d$ with the property that for all $n \in \mathbb{N}$, $i\in\{1,2,\dots,d\}$, $\theta\in\R^d$ it holds that
\begin{equation}
    (\nabla\smalll)(\vartheta)=0,\qquad\convexcoe\|\theta\|^2\leq \theta^\top(\operatorname{Hess}\smalll) (\vartheta)\theta\leq \convexcoe\chi \| \theta \|^2,
\end{equation}
\begin{equation}
\bfm_0=0,\qquad \bfm_n= \alpha \bfm_{n-1}+(1-\alpha)(\nabla\smalll)(\Theta_{n-1}),
\end{equation}
         \begin{equation}
            \bbV_0=0,\qquad \bbV_n^i= \beta \bbV_{n-1}^i+(1-\beta)|(\nabla_{\theta_i}\smalll)(\Theta_{n-1})|^2,
            \end{equation}
            \begin{equation}
            \|\Theta_0-\vartheta\|<\delta,\qquad\text{and}\qquad  \Theta_n^i=\Theta_{n-1}^i-\gamma_n
        \textstyle \bigl[\frac{\bfm_{n}^i}{1-\alpha^n}\bigr] \Bigl[\varepsilon+\bigl[\frac{\democrator_{n}^i}{1-\beta^n}\bigr]^{\nicefrac{1}{2}}\Bigr]^{-1}
         \end{equation}
          we have for all $n\in\N_0$ that
$  \|\Theta_n-\vartheta\|\leq \fC\|\Theta_0-\vartheta\|\varrho^n$\dott}
\end{aproof}
\newpage
\subsection{Lower bounds on convergence rates of Adam}
\begin{tcolorbox}[colback=white!95!gray,
                  colframe=black,
                  boxrule=0.5pt,
                  sharp corners,
                  enhanced,
                  breakable,
                 ]
\begin{athm}{lemma}{lem: lower bound adam}
Let $d\in\N\backslash\{1\}$, $\convexcoe,\delta,\gamma,\varepsilon,\fC\in (0,\infty)$, $\chi\in (1,\infty)$, $\alpha=(\frac{\sqrt{\chi}-1}{\sqrt{\chi}+1})^2$, $\beta\in (0,1)$, $\varrho\in (0,\frac{\sqrt{\chi}-1}{\sqrt{\chi}+1})$ and let $\smalll=(\smalll(\theta))_{\theta=(\theta_1,\dots,\theta_d)\in\R^d}\colon \R^d\to\R$ satisfy for all $\theta=(\theta_1,\dots,\theta_d)\in \R^d$ that $\smalll(\theta)=\convexcoe\chi|\theta_1|^2+\sum_{i=2}^d \convexcoe|\theta_i|^2$. Then there exist $\Theta=(\Theta^1,\dots,\Theta^d)\colon \mathbb{N}_0\to\R^d$, $\bfm=(\bfm^1,\dots,\bfm^d)\colon\N_0\to\R^d$, and  $\bbV=(\bbV^1,\dots,\bbV^d)\colon \mathbb{N}_0\to\R^d$ such that
\begin{enumerate}[label=(\roman*)]
    \item \llabel{item 1} it holds for all $n \in \mathbb{N}$, $i\in\{1,2,\dots,d\}$ that
          \begin{equation}
    \textcolor{magenta}{\bfm_0=0},\qquad \textcolor{magenta}{\bfm_n= \alpha \bfm_{n-1}+(1-\alpha)(\nabla\smalll)(\Theta_{n-1})},
\end{equation}
         \begin{equation}
            \textcolor{magenta}{\bbV_0=0},\qquad \textcolor{magenta}{\bbV_n^i= \beta \bbV_{n-1}^i+(1-\beta)|(\nabla_{\theta_i}\smalll)(\Theta_{n-1})|^2},
            \end{equation}
            \begin{equation}
            \textcolor{magenta}{\|\Theta_0\|< \delta},\qquad\text{and}\qquad  \textcolor{magenta}{\Theta_n^i=\Theta_{n-1}^i-\gamma
        \textstyle \bigl[\frac{\bfm_{n}^i}{1-\alpha^n}\bigr] \Bigl[\varepsilon+\bigl[\frac{\democrator_{n}^i}{1-\beta^n}\bigr]^{\nicefrac{1}{2}}\Bigr]^{-1}}
         \end{equation}
         and
         \item \llabel{item 2} there exists $n\in\N$ such that $\textcolor{magenta}{ \textstyle\|\Theta_n\|>\fC\|\Theta_0\|\varrho^n}$.
\end{enumerate} 
\end{athm}
\end{tcolorbox}
\begin{aproof}
Throughout this proof let $\sgn\colon\R\to\R$ satisfy for all $x\in\R$ that
    \begin{equation}\llabel{def: sgn}
        \sgn(x)=\mathbbm 1_{(0,\infty)}(x)-\mathbbm 1_{(-\infty,0)}(x)
    \end{equation}
    and for every $i\in\{1,2,\dots,d\}$ let $h_i\in\R$ satisfy 
    \begin{equation}\llabel{def: h}
        h_i=1+(\chi-1)\mathbbm 1_{\{1\}}(i).
    \end{equation}
    \argument{the fact that $\alpha=(\frac{\sqrt{\chi}-1}{\sqrt{\chi}+1})^2$}{that
    \begin{equation}\llabel{eqq0}
        \textstyle\chi=(\frac{1+\sqrt{\alpha}}{1-\sqrt{\alpha}})^2.
    \end{equation}}
    \argument{\lref{eqq0};the fact that $\chi>0$}{that
    \begin{equation}\llabel{eqq1}
        (2(1-\alpha)\convexcoe\gamma\varepsilon^{-1}-(1-\sqrt{\alpha})^2)(2(1-\alpha)\convexcoe\gamma\chi\varepsilon^{-1}-(1+\sqrt{\alpha})^2)=\chi(2(1-\alpha)\convexcoe\gamma\varepsilon^{-1}-(1-\sqrt{\alpha})^2)^2\geq 0.
    \end{equation}}
    \argument{\lref{eqq1};the fact that $2(1-\alpha)\convexcoe\chi\gamma\varepsilon^{-1}>2(1-\alpha)\convexcoe\gamma\varepsilon^{-1}>0$}{that 
    \begin{equation}\llabel{eqq2}
        2(1-\alpha)\convexcoe\gamma\varepsilon^{-1}\leq (1-\sqrt{\alpha})^2\qquad \text{or}\qquad  2(1-\alpha)\convexcoe\chi\gamma\varepsilon^{-1}\geq (1+\sqrt{\alpha})^2.
    \end{equation}}
    \argument{\lref{eqq2};}{that
    \begin{equation}\llabel{eqq3}
        1+\alpha- 2(1-\alpha)\convexcoe\gamma\varepsilon^{-1}\geq 2\sqrt{\alpha}\qquad \text{or}\qquad 1+\alpha- 2(1-\alpha)\convexcoe\chi\gamma\varepsilon^{-1}\leq -2\sqrt{\alpha}. 
    \end{equation}}
   \argument{\lref{eqq3};the fact that $\{h_1,h_2,\dots,h_d\}=\{1,\chi\}$}{that there exists  
    $i\in\{1,2,\dots,d\}$ which satisfies that
    \begin{equation}\llabel{def: i}
       |1+\alpha- 2(1-\alpha)rh_i\gamma\varepsilon^{-1}|\geq 2\sqrt{\alpha}.
    \end{equation}}
    \startnewargseq
In the following let $s\in\R$ satisfy 
\begin{equation}\llabel{def: s}
    s=1+\alpha-2(1-\alpha)\convexcoe h_i\gamma\varepsilon^{-1},
\end{equation}
let $\bfe=(\bfe_1,\dots,\bfe_d)\in \R^d$ satisfy for all $j\in\{1,2,\dots,d\}$ that
    \begin{equation}\llabel{def: bfe}
        \bfe_{j}=\mathbbm 1_{\{i\}}(j),
    \end{equation}
     for every $t\in\R$ let $\Theta^{t}=(\Theta^{t,1},\dots,\Theta^{t,d})\colon\N_0\to\R^d$, $\bfm^t=(\bfm^{t,1},\dots,\bfm^{t,d})\colon\N_0\to\R^d$, and $\bbV^{t}=(\bbV^{t,1},\dots,\bbV^{t,d})\colon\N_0\to\R^d$ satisfy for all $n\in\N$, $j\in\{1,2,\dots,d\}$ that
     \begin{equation}\llabel{def: bfm}
      \bfm^{t}_0=0,\qquad \bfm^{t,j}_n= \alpha \bfm^{t,j}_{n-1}+2(1-\alpha)rh_j\Theta^{t,j}_{n-1},
            \end{equation}
            \begin{equation}\llabel{def: bbV}
      \bbV^{t}_0=0,\qquad \bbV^{t,j}_n= \beta \bbV^{t,j}_{n-1}+(1-\beta)|2rh_j\Theta^{t,j}_{n-1}|^2,
            \end{equation}
            \begin{equation}\llabel{def: Theta}
        \Theta^{t}_0=t\bfe,\qquad\text{and}\qquad \Theta^{t,j}_n=\Theta^{t,j}_{n-1}-\gamma
        \textstyle \bigl[\frac{\bfm_{n}^{t,j}}{1-\alpha^n}\bigr] \Bigl[\varepsilon+\bigl[\frac{\bbV_{n}^{t,j}}{1-\beta^n}\bigr]^{\nicefrac{1}{2}}\Bigr]^{-1},
         \end{equation}
         and let $\bfu\colon \N_0\to\R$ and $\bfv\colon\N_0\to\R$ satisfy for all $n\in\N$ that
         \begin{equation}\llabel{def: bfu}
             \bfu_0=1,\qquad \bfu_n=\bfu_{n-1}-\gamma
        \textstyle \bigl[\frac{\bfv_{n}}{1-\alpha^n}\bigr] \varepsilon^{-1},
             \end{equation}
             \begin{equation}\llabel{def: bfv}
        \bfv_0=0,\qqandqq \bfv_n=\alpha\bfv_{n-1}+2(1-\alpha)rh_i \bfu_{n-1}.
         \end{equation}
    \argument{\lref{def: bfu};}{that for all $n\in\N$ it holds that
    \begin{equation}\llabel{eq1}
       \textstyle \bfv_n=-(\bfu_n-\bfu_{n-1})\gamma^{-1}\varepsilon (1-\alpha^n).
    \end{equation}}  
    \argument{\lref{eq1};\lref{def: bfv}}{for all $n\in\N\backslash\{1\}$ that
    \begin{equation}\llabel{eq2'}
       -(\bfu_n-\bfu_{n-1})\gamma^{-1}\varepsilon (1-\alpha^n)=-\alpha(\bfu_{n-1}-\bfu_{n-2})\gamma^{-1}\varepsilon (1-\alpha^{n-1})+2(1-\alpha)rh_i\bfu_{n-1}.
    \end{equation}}
    \argument{\lref{eq2'};}{for all $n\in\N\backslash\{1\}$ that
    \begin{equation}\llabel{eq2}
       (\bfu_n-\bfu_{n-1})(1-\alpha^n)=\alpha(\bfu_{n-1}-\bfu_{n-2}) (1-\alpha^{n-1})-2(1-\alpha)rh_i\gamma\varepsilon^{-1}\bfu_{n-1}.
    \end{equation}}
    \argument{\lref{eq2};\lref{def: s};}{for all $n\in\N\backslash\{1\}$ that
    \begin{equation}\llabel{eq3}
    \begin{split}
        (1-\alpha^{n})\bfu_n&=(1+\alpha-2\alpha^n-2(1-\alpha)rh_i\gamma\varepsilon^{-1})\bfu_{n-1}-\alpha(1-\alpha^{n-1})\bfu_{n-2}\\
        &=(s-2\alpha^n)\bfu_{n-1}-\alpha(1-\alpha^{n-1})\bfu_{n-2}.
        \end{split}
    \end{equation}}
    In the following we prove that for all $n\in\N$ it holds that
    \begin{equation}\llabel{need to prove}
        (\sgn(s))^n\bfu_{n}\geq\sqrt{\alpha}(\sgn(s))^{n-1}\bfu_{n-1}>0.
    \end{equation}
    We prove \lref{need to prove} by induction on $n\in\N$. For the base case $n=1$ note that \lref{def: bfv} shows that \llabel{base: eq1}
        $\bfv_1=2(1-\alpha)rh_i$.
    This, \lref{def: s}, and \lref{def: bfu} establish that
    \begin{equation}\llabel{base: eq2}
        \bfu_1=1-2\gamma\varepsilon^{-1}rh_i=\frac{s-2\alpha}{1-\alpha}.
    \end{equation}
    Combining this with \lref{def: i} and \lref{def: bfu} demonstrates that
    \begin{equation}\llabel{base: eq3}
        \bfu_1\sgn(s)=\frac{|s|-2\alpha\sgn(s)}{1-\alpha}\geq \frac{2\sqrt{\alpha}-2\alpha}{1-\alpha}=\frac{2\sqrt{\alpha}}{1+\sqrt{\alpha}}\geq \sqrt{\alpha}=\sqrt{\alpha}(\sgn(s))^0\bfu_0.
    \end{equation}
    This ensures \lref{need to prove} in the base case $n=1$. For the induction step we assume that there exists $n\in\N\backslash\{1\}$ which satisfies for all $k\in\{1,2,\dots,n-1\}$ that
    \begin{equation}\llabel{ntp: assume}
        (\sgn(s))^k\bfu_{k}\geq\sqrt{\alpha}(\sgn(s))^{k-1}\bfu_{k-1}>0.
    \end{equation}
    \startnewargseq
    \argument{\lref{eq3};\lref{ntp: assume};}{that
    \begin{equation}\llabel{ntp: eq1}
    \begin{split}
       &(1-\alpha^n) (\sgn(s))^n\bfu_n\\
       &=(s-2\alpha^n)(\sgn(s))^{n}\bfu_{n-1}-\alpha(1-\alpha^{n-1})(\sgn(s))^{n}\bfu_{n-2}\\
       &=|s|(\sgn(s))^{n-1}\bfu_{n-1}-2\alpha^n(\sgn(s))^{n}\bfu_{n-1}-\alpha(1-\alpha^{n-1})(\sgn(s))^{n-2}\bfu_{n-2}\\
       &\geq |s|(\sgn(s))^{n-1}\bfu_{n-1}-2\alpha^n(\sgn(s))^{n}\bfu_{n-1}-\sqrt{\alpha}(1-\alpha^{n-1})(\sgn(s))^{n-1}\bfu_{n-1}\\
       &\geq (\sgn(s))^{n-1}\bfu_{n-1}(2\sqrt{\alpha}-
       2\alpha^n-\sqrt{\alpha}(1-\alpha^{n-1}))\\
       &=(\sgn(s))^{n-1}\bfu_{n-1}(\sqrt{\alpha}-2\alpha^n+\alpha^{n-1/2})\geq (\sgn(s))^{n-1}\bfu_{n-1}(\sqrt{\alpha}-\alpha^{n+1/2})\\
       &=\sqrt{\alpha}(\sgn(s))^{n-1}\bfu_{n-1}(1-\alpha^n).
       \end{split}
    \end{equation}}
    \argument{\lref{ntp: eq1};\lref{ntp: assume}}{that
    \begin{equation}\llabel{ntp: eq2}
        (\sgn(s))^n\bfu_n\geq \sqrt{\alpha}(\sgn(s))^{n-1}\bfu_{n-1}>0.
    \end{equation}}
    \argument{\lref{ntp: eq2};\lref{ntp: assume};induction}{\lref{need to prove}\dott}
    \startnewargseq
    \argument{\lref{need to prove};}{that for all $n\in\N$ it holds that \llabel{eq5}
        $|\bfu_n|\geq \sqrt{\alpha}
        |\bfu_{n-1}|$\dott}
    \argument{\lref{def: bfu};\lref{eq5};induction}{that for all $n\in\N$ it holds that
    \begin{equation}\llabel{eq6}
        |\bfu_n|\geq (\sqrt{\alpha})^n|\bfu_0|=(\sqrt{\alpha})^n.
    \end{equation}}
       \argument{\lref{eq6}; the fact that $\varrho<\sqrt{\alpha}$;}{that there exists $N\in\N$ which satisfies
    \begin{equation}\llabel{def: N}
        |\bfu_N|> 2\fC\varrho^N.
    \end{equation}}
    \startnewargseq
     In the following we prove that for all $n\in\N_0$, $t\in \R$, $j\in\{1,2,\dots,d\}$ it holds that
         \begin{equation}\llabel{need to prove 2}
         \limsup_{t\to 0} \Bigl|\frac{\Theta_n^{t,j}}{t}-\bfu_n\mathbbm 1_{\{i\}}(j)\Bigr|=0,\qquad \limsup_{t\to 0} \Bigl|\frac{\bfm_n^{t,j}}{t}-\bfv_n\mathbbm 1_{\{i\}}(j)\Bigr|=0,\qqandqq \limsup_{t\to 0} |\bbV_n^{t,j}|=0.
         \end{equation}
         We prove \lref{need to prove 2} by induction on $n\in\N_0$. For the base case $n=0$ note that \lref{def: bfe}, \lref{def: bfm}, \lref{def: Theta}, \lref{def: bfu}, and \lref{def: bfv} show that for all $t\in \R$, $j\in\{1,2,\dots,d\}$ it holds that
         \begin{equation}
             \Theta_0^{t,j}=t\bfe_j=t\bfu_0\mathbbm 1_{\{i\}}(j),\qquad \bfm_0^{t,j}=0=t\bfv_0\mathbbm 1_{\{i\}}(j),\qqandqq \bbV_0^{t,j}=0.
         \end{equation}
         This proves \lref{need to prove 2} in the base case $n=0$. For the induction step we assume that there exists $n\in\N$ which satisfies for all $t\in\R$, $j\in\{1,2,\dots,d\}$ that
         \begin{equation}\llabel{ntp2: assume}
              \textstyle \limsup\limits_{t\to 0} \Bigl|\frac{\Theta_{n-1}^{t,j}}{t}-\bfu_{n-1}\mathbbm 1_{\{i\}}(j)\Bigr|=0, \quad \limsup\limits_{t\to 0} \Bigl|\frac{\bfm_{n-1}^{t,j}}{t}-\bfv_{n-1}\mathbbm 1_{\{i\}}(j)\Bigr|=0,\qandq\limsup\limits_{t\to 0} |\bbV_{n-1}^{t,j}|=0.
         \end{equation}
         \startnewargseq
         \argument{\lref{ntp2: assume};}{that for all $j\in\{1,2,\dots,d\}$ it holds that
         \begin{equation}\llabel{ntp2: eq0}
             \textstyle \limsup_{t\to 0}|\Theta_{n-1}^{t,j}|=0.
         \end{equation}}
         \argument{\lref{ntp2: eq0};\lref{def: bbV};\lref{ntp2: assume}}{for all  $j\in\{1,2,\dots,d\}$ that
         \begin{equation}\llabel{ntp2: eq1}
         \begin{split}
            \textstyle \limsup_{t\to 0}|\bbV_n^{t,j}|&\leq \textstyle \beta \bigl[\limsup_{t\to 0}\bbV^{t,j}_{n-1}\big]+(1-\beta)\bigl[\limsup_{t\to 0}|2rh_j\Theta^{t,j}_{n-1}|^2\bigr]=0.
             \end{split}
         \end{equation}}
         \argument{\lref{def: bfm};\lref{def: bfv};
         \lref{ntp2: assume}}{for all $j\in\{1,2,\dots,d\}$ that
         \begin{equation}\llabel{ntp2: eq1'}
         \begin{split}
              \textstyle \lim_{t\to 0}\frac{\bfm_{n}^{t,j}}{t}&\textstyle=\alpha\bigl[\lim_{t\to 0}\frac{\bfm_{n-1}^{t,j}}{t}\bigr]+2(1-\alpha)\bigl[\lim_{t\to 0}\frac{rh_j\Theta_{n-1}^{t,j}}{t}\bigr]\\
              &=\alpha\bfv_{n-1}\mathbbm 1_{\{i\}}(j)+2(1-\alpha)\convexcoe h_j\bfu_{n-1}\mathbbm 1_{\{i\}}(j)=\bfv_{n}\mathbbm 1_{\{i\}}(j).
              \end{split}
         \end{equation}}
         \argument{\lref{def: Theta};\lref{def: bfu};\lref{ntp2: assume};;\lref{ntp2: eq1};\lref{ntp2: eq1'}}{for all $j\in\{1,2,\dots,d\}$ that
         \begin{equation}\llabel{ntp2: eq2}
         \begin{split}
             \textstyle\lim_{t\to 0}\frac{\Theta_n^{t,j}}{t}&\textstyle=\bigl[\lim_{t\to 0}\frac{\Theta^{t,j}_{n-1}}{t}\bigr]-\gamma
          \bigl[\textstyle \lim_{t\to 0}\bigl[\varepsilon+\bigl(\frac{\bbV_n^{t,j}}{1-\beta^n})^{1/2}\bigr]^{-1}\bigr]\bigl[\lim_{t\to 0}\frac{\bfm_{n}^{t,j}}{t(1-\alpha^n)}\bigr]\\
          &=\textstyle\bfu_{n-1}\mathbbm 1_{\{i\}}(j)-\gamma\varepsilon^{-1}\bigl[\frac{\bfv_{n}\mathbbm 1_{\{i\}}(j)}{1-\alpha^n}\bigr]=\bfu_{n}\mathbbm 1_{\{i\}}(j).
          \end{split}
         \end{equation}}
         \argument{\lref{ntp2: eq2};
         \lref{ntp2: eq1};\lref{ntp2: eq1'}}{for all $t\in\R$, $j\in\{1,2,\dots,d\}$ that
         \begin{equation}\llabel{ntp2: eq4}
   \textstyle \limsup\limits_{t\to 0} \Bigl|\frac{\Theta_{n}^{t,j}}{t}-\bfu_{n}\mathbbm 1_{\{i\}}(j)\Bigr|=0, \qquad \limsup\limits_{t\to 0} \Bigl|\frac{\bfm_{n}^{t,j}}{t}-\bfv_{n}\mathbbm 1_{\{i\}}(j)\Bigr|=0,\qqandqq\limsup\limits_{t\to 0} |\bbV_{n}^{t,j}|=0.              
         \end{equation}}
         \argument{\lref{ntp2: eq4};\lref{ntp2: assume};induction}{\lref{need to prove 2}\dott}
         \startnewargseq
         \argument{\lref{def: N};\lref{need to prove 2};}{that
         \begin{equation}\llabel{eq10}
             \textstyle \lim_{t\to 0}\frac{\|\Theta_N^t\|}{|t|}=|\bfu_N|> 2\fC\varrho^N.
         \end{equation}}
         \argument{\lref{eq10};\lref{def: Theta};}{that there exists $t\in(-\delta,\delta)$ such that
         \begin{equation}\llabel{eq11}
          \|\Theta_N^t\|>\fC\varrho^N |t|=\fC\|\Theta_0^t\|\varrho^N. 
         \end{equation}}
         \argument{\lref{def: h};the assumption that for all $\theta=(\theta_1,\dots,\theta_d)\in \R^d$ it holds that $\smalll(\theta)=\convexcoe\chi|\theta_1|^2+\sum_{i=2}^d \convexcoe|\theta_i|^2$}{that for all $j\in\{1,2,\dots,d\}$, $\theta=(\theta_1,\dots,\theta_d)\in\R^d$ it holds that
         \begin{equation}\llabel{eq12}
           (\nabla_{\theta_j}\smalll)(\theta)=2r h_j\theta_j.  
         \end{equation}}
         \argument{\lref{eq12};\lref{def: bfm};\lref{def: bbV};\lref{def: Theta};\lref{eq11}}{\lref{item 1,item 2}\dott}
\end{aproof}
\begin{tcolorbox}[colback=white!95!gray,
                  colframe=black,
                  boxrule=0.5pt,
                  sharp corners,
                  enhanced,
                  breakable,
                 ]
\begin{athm}{cor}{main theorem lower bound adam}
     Let $d\in\N\backslash\{1\}$, $\chi\in (1,\infty)$, $\alpha=(\frac{\sqrt{\chi}-1}{\sqrt{\chi}+1})^2$, $\beta\in (0,1)$, $\varrho\in (0,\frac{\sqrt{\chi} - 1 }  { \sqrt{\chi} + 1} )$, $\varepsilon\in (0,\infty)$. Then there exist $\convexcoe\in (0,\infty)$, $\vartheta\in\R^d$, and $\smalll\in C^2( \R^d,\R)$ with $(\nabla \smalll )( \vartheta ) = 0$ and $
			\forall\, \theta\in\R^d\colon \convexcoe\|\theta\|^2\leq \smalll''(\vartheta)(\theta,\theta)\leq \convexcoe\chi \| \theta \|^2
			$ such that for all $\delta,\Gamma,\fC \in (0,\infty)$ there exist $\gamma\colon\N\to (\Gamma-\delta,\Gamma+\delta)$, $\Theta=(\Theta^1,\dots,\Theta^d)\colon \mathbb{N}_0\to\R^d$, $\bfm=(\bfm^1,\dots,\bfm^d)\colon\N_0\to\R^d$, and  $\bbV=(\bbV^1,\dots,\bbV^d)\colon \mathbb{N}_0\to\R^d$ such that
\begin{enumerate}[label=(\roman*)]
    \item \llabel{item 1} it holds for all $n \in \mathbb{N}$, $i\in\{1,2,\dots,d\}$  that
          \begin{equation}
    \textcolor{magenta}{\bfm_0=0},\qquad \textcolor{magenta}{\bfm_n= \alpha \bfm_{n-1}+(1-\alpha)(\nabla\smalll)(\Theta_{n-1})},
\end{equation}
         \begin{equation}
            \textcolor{magenta}{\bbV_0=0},\qquad \textcolor{magenta}{\bbV_n^i= \beta \bbV_{n-1}^i+(1-\beta)|(\nabla_{\theta_i}\smalll)(\Theta_{n-1})|^2},
            \end{equation}
            \begin{equation}
            \textcolor{magenta}{\|\Theta_0-\vartheta\|< \delta},\qquad\text{and}\qquad  \textcolor{magenta}{\Theta_n^i=\Theta_{n-1}^i-\gamma_n
        \textstyle \bigl[\frac{\bfm_{n}^i}{1-\alpha^n}\bigr] \Bigl[\varepsilon+\bigl[\frac{\democrator_{n}^i}{1-\beta^n}\bigr]^{\nicefrac{1}{2}}\Bigr]^{-1}}
         \end{equation}
         and
         \item \llabel{item 2} there exists $n\in\N$ such that $\textcolor{magenta}{ \textstyle
   \textstyle\|\Theta_n-\vartheta\|>\fC\|\Theta_0-\vartheta\|\varrho^n}$.
\end{enumerate} 
\end{athm}
\end{tcolorbox}
    \begin{aproof}
    Throughout this proof let $\convexcoe\in\R$, $\vartheta\in\R^d$ satisfy $\convexcoe=2$ and $\vartheta=0$ and let $\smalll\colon\R^d\to\R$ satisfy for all $\theta=(\theta_1,\dots,\theta_d)\in\R^d$ that $\smalll(\theta)=\chi|\theta_1|^2+\sum_{i=2}^d|\theta_i|^2$.
\argument{\cref{lem: lower bound adam};}{\lref{item 1,item 2}\dott}
    \end{aproof}
    \subsection{Sharp convergence rates of Adam}\label{subsec: sharp adam}
     \begin{samepage}
 \begin{tcolorbox}[colback=white!95!gray,
                  colframe=black,
                  boxrule=0.5pt,
                  sharp corners,
                  enhanced,
                  unbreakable,
                 ]
                 \begin{athm}{theorem}{theo: sharp rate adam}
                 Let $d\in \mathbb{N}\backslash\{1\}$, $\chi\in (1,\infty)$, $\alpha=(\frac{ \sqrt{\chi}  - 1 }  {  \sqrt{\chi}  + 1 })^2$, $\beta\in(0,1)$, $\varrho=\frac{ \sqrt{\chi}  - 1 }  {  \sqrt{\chi}  + 1 }$, $\varepsilon\in (0,\infty)$ and for every $n\in\N$ let $\varphi_n=(\varphi_n^1,\dots,\varphi_n^d)\colon (\R^d)^n\to\R^d$ satisfy for all $\theta=((\theta_k^i)_{i\in\{1,2,\dots,d\}})_{k\in\{1,2,\dots,n\}}\in (\R^d)^n$, $i\in\{1,2,\dots,d\}$ that
\begin{equation}\llabel{def: Phi}
\textcolor{magenta}{\textstyle\varphi_n^i(\theta)=\bigl[\varepsilon+\bigl(\frac{1-\beta}{1-\beta^n}\sum_{k=1}^{n}\beta^{n-k}|\theta_k^i|^2\bigr)^{1/2}\bigr]^{-1}\bigl(\frac{1-\alpha}{1-\alpha^n}\sum_{k=1}^{n}\alpha^{n-k}\theta_k^i\bigr)}.
\end{equation}
Then it holds for all $\rho \in (0,\infty) \backslash \{ \varrho \}$ that the following three statements are equivalent:
             \begin{enumerate}[label=(\roman*)]
                 \item \llabel{item 1} It holds that $\rho > 
                 \varrho$.
                 \item \llabel{item 2} For all $\convexcoe \in (0,\infty)$, $\cL \in C^2( \R^d, \R )$, $\vartheta\in(\nabla\cL)^{-1}(\{0\})$ with $
			\forall\, \theta\in\R^d\colon \convexcoe\|\theta\|^2\leq \cL'' (\vartheta)(\theta,\theta)\leq \convexcoe\chi \| \theta \|^2
			$ there exist $\delta,\Gamma,\fC \in (0,\infty)$ such that for every $\gamma=(\gamma_n)_{n\in\N}\colon\N\to (\Gamma-\delta,\Gamma+\delta)$, every $\bfm=(\bfm^1,\dots,\bfm^d)\colon\N_0\to\R^d$, every $\bbV=(\bbV^1,\dots,\bbV^d)\colon \mathbb{N}_0\to\R^d$, and every $\Theta=(\Theta^1,\dots,\Theta^d)\colon \mathbb{N}_0\to\R^d$ with the property that 
    for all $n \in \mathbb{N}$, $i\in\{1,2,\dots,d\}$ it holds that
         \begin{equation}
    \textcolor{magenta}{\bfm_0=0},\qquad \textcolor{magenta}{\bfm_n= \alpha \bfm_{n-1}+(1-\alpha)(\nabla\cL)(\Theta_{n-1})},
\end{equation}
         \begin{equation}
            \textcolor{magenta}{\bbV_0=0},\qquad \textcolor{magenta}{\bbV_n^i= \beta \bbV_{n-1}^i+(1-\beta)|(\nabla_{\theta_i}\cL)(\Theta_{n-1})|^2},
            \end{equation}
            \begin{equation}
            \textcolor{magenta}{\|\Theta_0-\vartheta\|< \delta},\qquad\text{and}\qquad  \textcolor{magenta}{\Theta_n^i=\Theta_{n-1}^i-\gamma_n
        \textstyle \bigl[\frac{\bfm_{n}^i}{1-\alpha^n}\bigr] \Bigl[\varepsilon+\bigl[\frac{\democrator_{n}^i}{1-\beta^n}\bigr]^{\nicefrac{1}{2}}\Bigr]^{-1}}
         \end{equation}
         we have for all $n\in\N_0$ that $\textcolor{magenta}{ \textstyle
   \textstyle\|\Theta_n-\vartheta\|\leq\fC\rho^n\|\Theta_0-\vartheta\|}$.
                 \item \llabel{item 3} For all $\convexcoe \in (0,\infty)$, $\cL \in C^2(\R^d, \R )$, $\vartheta\in(\nabla\cL)^{-1}(\{0\})$ with $
			\forall\, \theta\in\R^d\colon \convexcoe\|\theta\|^2\leq \cL'' (\vartheta)(\theta,\theta)\leq \convexcoe\chi \| \theta \|^2
			$ there exist $\delta,\Gamma, \fC \in (0,\infty) $ such that 
                  for every $\gamma=(\gamma_n)_{n\in\N}\colon\N\to (\Gamma-\delta,\Gamma+\delta)$ and every $ \Theta \colon \N_0\to\R^d $
			with $\| \Theta_0 - \vartheta \| < \delta$ and
			\begin{equation}
				 \textcolor{magenta}{\forall \, n \in \N \colon
				\Theta_n
				=
				\Theta_{ n - 1 }
				-
				\gamma_n \varphi_n\bigl(
				( \nabla \cL )( \Theta_0 ) ,
				( \nabla \cL )( \Theta_1 ) ,
				\dots ,
				( \nabla \cL )( \Theta_{ n - 1 } )
				\bigr)}
			\end{equation}
			we have for all $n\in\N_0$ that
			$
			\textcolor{magenta}{\| \Theta_n - \vartheta \|
			\leq \fC\rho^{  n } \|\Theta_0-\vartheta\|}
			$.
             \end{enumerate}
                 \end{athm}
                 \end{tcolorbox}
                 \end{samepage}
    \begin{aproof}
        Throughout this proof for every $\Theta\colon\N_0\to\R^d$ and $\cL=(\cL(\theta))_{\theta=(\theta_1,\dots,\theta_d)\in\R^d}\in C^1(\R^d,\allowbreak\R)$ let $\bfm(\Theta,\cL)=(\bfm_n^i(\Theta,\cL))_{(i,n)\in\{1,2,\dots,d\}\times\N_0}\colon\N_0\to\R^d$ and $\bbV(\Theta,\cL)=(\bbV_n^i(\Theta,\cL))_{(i,n)\in\{1,2,\dots,d\}\times\N_0}\allowbreak\colon\N_0\to\R^d$ satisfy for all $n\in\N_0$, $i\in\{1,2,\dots,d\}$ that
        \begin{equation}\llabel{def: bfm}
         \textstyle  \bfm_n(\Theta,\cL)=(1-\alpha)\sum\limits_{k=1}^n\alpha^{n-k}(\nabla\cL)(\Theta_{k-1})\qandq \bbV_n^i(\Theta,\cL)=(1-\beta)\sum\limits_{k=1}^n\beta^{n-k}|(\nabla_{\theta_i}\cL)(\Theta_{k-1})|^2.
        \end{equation}
        \argument{\cref{main theorem higher bound adam}}{that (\lref{item 1}$\rightarrow$\lref{item 2})\dott}
        \startnewargseq
         \argument{\cref{main theorem lower bound adam}}{that (\lref{item 2}$\rightarrow$\lref{item 1})\dott}
        \startnewargseq
        \argument{\lref{def: bfm};}{that for all $n\in\N$, $\cL=(\cL(\theta))_{\theta=(\theta_1,\dots,\theta_d)\in\R^d}\in C^1(\R^d,\R)$, $i\in\{1,2,\dots,d\}$ and all $\Theta\colon\N_0\to\R^d$ it holds that
        \begin{equation}\llabel{eq1'}
            \bfm_{n}(\Theta,\cL)=\alpha\bfm_{n-1}(\Theta,\cL)+(1-\alpha)(\nabla\cL)(\Theta_{n-1})
            \end{equation}
            and
            \begin{equation}\llabel{eq1}
                \bbV^i_{n}(\Theta,\cL)=\beta\bbV^i_{n-1}(\Theta,\cL)+(1-\beta)|(\nabla_{\theta_i}\cL)(\Theta_{n-1})|^2.
            \end{equation}}
       \argument{\lref{def: Phi};\lref{def: bfm}}{that for all $n\in\N$, $\cL=(\cL(\theta))_{\theta=(\theta_1,\dots,\theta_d)\in\R^d}\in C^1(\R^d,\R)$, $i\in\{1,2,\dots,d\}$ and all $\Theta\colon\N_0\to\R^d$ it holds that
       \begin{equation}\llabel{eq2}
       \begin{split}
      &  \varphi_n^i\bigl(
				( \nabla \cL )( \Theta_0 ) ,
				( \nabla \cL )( \Theta_1 ) ,
				\dots ,
				( \nabla \cL )( \Theta_{ n - 1 } )
				\bigr)
                \\&= \textstyle \bigl[\varepsilon+\bigl(\frac{1-\beta}{1-\beta^n}\sum_{k=1}^{n}\beta^{n-k}|(\nabla_{\theta_i}\cL)(\Theta_{k-1})|^2\bigr)^{1/2}\bigr]^{-1}\bigl(\frac{1-\alpha}{1-\alpha^n}\sum_{k=1}^{n}\alpha^{n-k}(\nabla_{\theta_i}\cL)(\Theta_{k-1})\bigr)\\
&\textstyle=\bigl[\varepsilon+\bigl(\frac{\bbV_n^i(\Theta,\cL)}{1-\beta^n}\bigr)^{1/2}\bigr]^{-1}\bigl[\frac{\bfm_n^i(\Theta,\cL)}{1-\alpha^n}\bigr].
                \end{split}
       \end{equation}}
       \argument{\lref{eq2};\lref{eq1'};\lref{eq1}}{that (\lref{item 2}$\rightarrow$\lref{item 3})\dott}
       \startnewargseq
       \argument{\unskip, \eg, \cite[Lemma 2.7]{DeRoArAd2025}}{that for every $\cL=(\cL(\theta))_{\theta=(\theta_1,\dots,\theta_d)\in\R^d}\in C^2(\R^d,\R)$, every $\Theta\colon\N_0\to\R^d$, every $\bfm=(\bfm^1,\dots,\bfm^d)\colon\N_0\to\R^d$, and every $\bbV=(\bbV^1,\dots,\bbV^d)\colon\N_0\to\R^d$ with the property that for all $n\in\N$, $i\in\{1,2,\dots,d\}$ it holds that $ \bfm_0=\bbV_0=0$, $\bfm_n= \alpha\bfm_{n-1}+(1-\alpha)(\nabla\cL)(\Theta_{n-1})$, and $\bbV_n^i= \beta\bbV_{n-1}^i+(1-\beta)|(\nabla_{\theta_i}\cL)(\Theta_{n-1})|^2$ we have for all $n\in\N_0$, $i\in\{1,2,\dots,d\}$ that
       \begin{equation}\llabel{eq4}
       \textstyle \bfm_n=(1-\alpha)\sum_{k=1}^n\alpha^{n-k}(\nabla\cL)(\Theta_{k-1})\qqandqq   \textstyle \bbV_n^i=(1-\beta)\sum_{k=1}^n\beta^{n-k}|(\nabla_{\theta_i}\cL)(\Theta_{k-1})|^2.
       \end{equation}}
       \argument{\lref{eq4};\lref{def: Phi}}{that for every $\cL=(\cL(\theta))_{\theta=(\theta_1,\dots,\theta_d)\in\R^d}\in C^2(\R^d,\R)$, every $\Theta\colon\N_0\to\R^d$, every $\bfm=(\bfm^1,\dots,\bfm^d)\allowbreak\colon\N_0\to\R^d$, and every $\bbV=(\bbV^1,\dots,\bbV^d)
       \allowbreak \colon\N_0\to\R^d$ with the property that for all $n\in\N$, $i\in\{1,2,\dots,d\}$ it holds that $ \bfm_0=\bbV_0=0$, $\bfm_n= \alpha\bfm_{n-1}+(1-\alpha)(\nabla\cL)(\Theta_{n-1})$, and $\bbV_n^i= \beta \bbV_{n-1}^i+(1-\beta)|(\nabla_{\theta_i}\cL)(\Theta_{n-1})|^2$ we have for all $n\in\N$, $i\in\{1,2,\dots,d\}$ that
       \begin{equation}\llabel{eq5}
       \begin{split}
           &\varphi^i_n\bigl(
				( \nabla \cL )( \Theta_0 ) ,
				( \nabla \cL )( \Theta_1 ) ,
				\dots ,
				( \nabla \cL )( \Theta_{ n - 1 } )
				\bigr)
               \\&= \textstyle \bigl[\varepsilon+\bigl(\frac{1-\beta}{1-\beta^n}\sum_{k=1}^{n}\beta^{n-k}|(\nabla_{\theta_i}\cL)(\Theta_{k-1})|^2\bigr)^{1/2}\bigr]^{-1}\bigl(\frac{1-\alpha}{1-\alpha^n}\sum_{k=1}^{n}\alpha^{n-k}(\nabla_{\theta_i}\cL)(\Theta_{k-1})\bigr)\\
&\textstyle=\bigl[\varepsilon+\bigl(\frac{\bbV_n^i}{1-\beta^n}\bigr)^{1/2}\bigr]^{-1}\bigl[\frac{\bfm_n^i}{1-\alpha^n}\bigr].
                \end{split}
       \end{equation}}
       \argument{\lref{eq5};}{that (\lref{item 3}$\rightarrow$\lref{item 2})\dott}
    \end{aproof}
 \begin{tcolorbox}[colback=white!95!gray,
                  colframe=black,
                  boxrule=0.5pt,
                  sharp corners,
                  enhanced,
                  breakable,
                 ]
                 \begin{athm}{cor}{cor: sharp rate adam}
                 Let $d\in \mathbb{N}\backslash\{1\}$, $\chi\in (1,\infty)$, $\alpha=(\frac{ \sqrt{\chi}  - 1 }  {  \sqrt{\chi}  + 1 })^2$, $\beta\in(0,1)$, $\varrho=\frac{ \sqrt{\chi}  - 1 }  {  \sqrt{\chi}  + 1 }$, $\varepsilon\in (0,\infty)$ and for every $n\in\N$ let $\varphi_n=(\varphi_n^1,\dots,\varphi_n^d)\colon (\R^d)^n\to\R^d$ satisfy for all $\theta=((\theta_k^i)_{i\in\{1,2,\dots,d\}})_{k\in\{1,2,\dots,n\}}\in (\R^d)^n$, $i\in\{1,2,\dots,d\}$ that
\begin{equation}\llabel{def: Phi}
\textcolor{magenta}{\textstyle\varphi_n^i(\theta)=\bigl[\varepsilon+\bigl(\frac{1-\beta}{1-\beta^n}\sum_{k=1}^{n}\beta^{n-k}|\theta_k^i|^2\bigr)^{1/2}\bigr]^{-1}\bigl(\frac{1-\alpha}{1-\alpha^n}\sum_{k=1}^{n}\alpha^{n-k}\theta_k^i\bigr)}.
\end{equation}
Then $(\varphi_n)_{n\in\N}$ converges with $\chi$-rate $\varrho$ (cf.\ \cref{definition: rate of convergence}).
                 \end{athm}
                 \end{tcolorbox}
                 \begin{aproof}
                     \argument{\cref{theo: sharp rate adam};}{that $(\varphi_n)_{n\in\N}$ converges with $\chi$-rate $\varrho$\dott}
                 \end{aproof}
\subsubsection*{Acknowledgements}
{\emergencystretch=2em This work has been partially supported by the National Natural Science Foundation of China (NSFC) under grant number W2531010. This work has also been partially funded by the Deutsche Forschungsgemeinschaft (DFG, German Research Foundation) under Germany's Excellence Strategy EXC \mbox{2044-390685587}, Mathematics M\"unster: Dynamics-Geometry-Structure. Most of the specific formulations in the proofs of this work have been created using~\cite{Bennoargumentcommand}.\par}

\subsubsection*{Use of large language models}

\cref{theo: muon higher bound}, \cref{lem: lower bound}, \cref{lem: lower bound rmsprop}, and \cref{lem: lower bound adam} in this work have been developed based on discussions with {\sc GPT-5.6 Sol Pro}. Specifically, the core ideas for the proofs of \cref{theo: muon higher bound}, \cref{lem: lower bound}, \cref{lem: lower bound rmsprop}, and \cref{lem: lower bound adam} are due to {\sc GPT-5.6 Sol Pro}. In addition, {\sc GPT-5.6 Sol Pro} has helped us to reveal several inaccuracies and errors within this manuscript. The entire article (including all statements and all arguments in this article) has been carefully written by the authors and the authors take full responsibility for each of the sentences/statements made in the paper.

\bibliographystyle{acm}
\bibliography{bibfileMUONconvergenceI}

\end{document}